\def\version{August 31, 2007}

\documentclass[12pt]{article}
\usepackage[psamsfonts]{amsfonts}
\usepackage{amsmath,amssymb,amsthm}
\usepackage[dvips]{graphicx}
\usepackage{psfrag}

\def\UseSection{
        \numberwithin{equation}{section}
	\theoremstyle{plain}
        \newtheorem{theorem}    {Theorem}[section]
        \DefineTheorems 
}

\def\DefineTheorems{
	
	\newtheorem{lemma}      [theorem] {Lemma}
	
	\newtheorem{prop}       [theorem] {Proposition}
	
	\newtheorem{cor}        [theorem] {Corollary}
	\newtheorem{ass}        [theorem] {Assumption}

	\theoremstyle{definition}
	\newtheorem{defn}       [theorem] {Definition}
	
	\newtheorem{example}       [theorem] {Example}

	\theoremstyle{definition}
	\newtheorem*{rem}	{Remark}
	
}

\newcommand{\bt}   {\begin{theorem}}
\newcommand{\et}   {\end  {theorem}}
\newcommand{\bl}   {\begin{lemma}}
\newcommand{\el}   {\end  {lemma}}
\newcommand{\bp}   {\begin{prop}}
\newcommand{\ep}   {\end  {prop}}
\newcommand{\bc}   {\begin{cor}}
\newcommand{\ec}   {\end  {cor}}
\newcommand{\bd}   {\begin{defn}}
\newcommand{\ed}   {\end  {defn}}

\newcommand{\ba}   {\begin{array}}
\newcommand{\ea}   {\end  {array}}
\newcommand{\be}   {\begin{enumerate}}
\newcommand{\ee}   {\end  {enumerate}}
\newcommand{\bi}   {\begin{itemize}}
\newcommand{\ei}   {\end  {itemize}}

\def\eq#1\en{\begin{equation}#1\end{equation}}  
\def\eqsplit#1\ensplit{
	\begin{equation}\begin{split}#1\end{split}\end{equation}
	}
\def\eqalign#1\enalign{
	\begin{align}#1\end{align}
	}
\def\eqmul#1\enmul{
	\begin{multline}#1\end{multline}
	}
\newcommand{\eqarrstar} {\begin{eqnarray*}} 
\newcommand{\enarrstar} {\end{eqnarray*}} 
\newcommand{\eqarray}   {\begin{eqnarray}} 
\newcommand{\enarray}   {\end{eqnarray}} 
\newcommand{\nnb}	{\nonumber \\} 

\newcommand{\lbeq}[1]  {\label{e:#1}}
\newcommand{\refeq}[1] {\eqref{e:#1}}    
\newcommand{\lbfg}[1]  {\label{fg: #1}}
\newcommand{\reffg}[1] {\ref{fg: #1}}

\makeatletter
\newcommand{\labelcounter}[2]{{%
	\stepcounter{#1}
	\protected@write\@auxout{}%
	{\string\newlabel{#2}{{\csname the#1\endcsname}{\thepage}}}%
	{\ref{#2}}
	}}
\makeatother

\newcommand{\sss}   { \scriptscriptstyle } 

\newcommand{\Ebold} {{\mathbb E}}

\newcommand{\Nbold} {{\mathbb N}}

\newcommand{\Pbold} {{\mathbb P}}
\newcommand{\Qbold} {{\mathbb Q}}
\newcommand{\Rbold} {{\mathbb R}}

\newcommand{\Zbold} {{\mathbb Z}}

\newcommand{\pvec}  {\boldsymbol{p}}
\newcommand{\qvec}  {\boldsymbol{q}}
\newcommand{\rvec}  {\boldsymbol{r}}

\newcommand{\uvec}  {\boldsymbol{u}}
\newcommand{\vvec}  {\boldsymbol{v}}
\newcommand{\wvec}  {\boldsymbol{w}}
\newcommand{\xvec}  {\boldsymbol{x}}
\newcommand{\yvec}  {\boldsymbol{y}}
\newcommand{\zvec}  {\boldsymbol{z}}
\newcommand{\zerovec} {\boldsymbol{0}}

\newcommand{\Ccal}   {\mathcal{C}} 
 
\newcommand{\Ecal}   {\mathcal{E}} 
\newcommand{\Fcal}   {\mathcal{F}}

\newcommand{\Ncal}   {\mathcal{N}} 
\newcommand{\Ocal}   {\mathcal{O}}

\newcommand{\Rd}    {{ {\Rbold}^d}}
\newcommand{\Zd}    {{ {\Zbold}^d }}

\newcommand{\spose}[1] {{\hbox to 0pt{#1\hss}} }
\newcommand{\ltapprox} {\mathrel{\spose{\lower 3pt\hbox{$\mathchar"218$}}
 \raise 2.0pt\hbox{$\mathchar"13C$}}}
\newcommand{\gtapprox} {\mathrel{\spose{\lower 3pt\hbox{$\mathchar"218$}}
 \raise 2.0pt\hbox{$\mathchar"13E$}}}

\newcommand{\nin}  {{ \not\in }}

\UseSection  
\renewcommand{\to}      {\rightarrow}

\newcommand{\bE}{\Ebold}
\newcommand{\bP}{\Pbold}

\newcommand{\Z}{\Zbold}
\newcommand{\bZ}{\Zbold}
\newcommand{\N}{\Nbold}
\newcommand{\conn}{\longrightarrow}

\newcommand{\nc}        { \conn  {\hspace{-3.0ex} /} \hspace{1.8ex}   }

\newcommand{\nn}{\nonumber}

\newcommand{\smallsup}[1] {{\scriptscriptstyle{({#1}})}}

\newcommand{\al}{\alpha}
\newcommand{\iictau}{\rho}

\newcommand{\tiZ}{{\tilde Z}}

\newcommand{\es}{{\varnothing}}
\newcommand{\vphi}{{\varphi}}

\newcommand{\ol}{\overline}
\newcommand{\tCcal}[1]{{\tilde\Ccal}^\smallsup{#1}}
\newcommand{\fnts}{\scriptsize}

\newcommand{\IIC}{{\hbox{\footnotesize\rm IIC}}}

\def\Reff{R_{\rm eff}}  
\def\eps{\varepsilon}
\def\lam{\lambda}
\def\Gam{\Gamma}
\def\Qi{\Qbold_\infty}
\def\zerovec{{\boldsymbol{0}}}

\title  {
        Random walk on the incipient infinite cluster \\
        for oriented percolation in high dimensions     }

\author{
Martin T.~Barlow\footnote{Department of Mathematics,
University of British Columbia,
Vancouver, BC V6T 1Z2, Canada.
{\tt barlow@math.ubc.ca}},
Antal A.~J\'arai\footnote{Carleton University,
School of Mathematics and Statistics,
1125 Colonel By Drive,
Ottawa, ON K1S 5B6, Canada.
{\tt jarai@math.carleton.ca}},
Takashi Kumagai\footnote{
Department of Mathematics,
Faculty of Science,
Kyoto University, Kyoto 606-8502, Japan.
{\tt kumagai@math.kyoto-u.ac.jp}},
Gordon Slade\footnote{Department of Mathematics,
University of British Columbia,
Vancouver, BC V6T 1Z2, Canada.
{\tt slade@math.ubc.ca}}}

\begin{document}

\date\version

\maketitle

\begin{abstract}
We consider simple random walk on the incipient infinite cluster for the
spread-out model of oriented percolation on $\Zd \times \Z_+$.
In dimensions $d>6$, we obtain
bounds on exit times, transition probabilities, and the range of the random walk,
which establish that the spectral dimension of the incipient infinite
cluster is $\frac 43$, and thereby prove
a version of the Alexander--Orbach conjecture
in this setting.   The proof
divides into two parts.  One part establishes general estimates
for simple random walk on an arbitrary infinite
random graph, given suitable bounds on volume and effective
resistance for the random graph.  A second part then provides these
bounds on volume and effective resistance for the incipient infinite
cluster in dimensions $d>6$, by extending results about critical oriented percolation
obtained previously via the lace expansion.
\end{abstract}

\section{Introduction and main results}

\subsection{Introduction}
\label{sec-intro}

The problem of random walk on a percolation cluster --- the `ant in the labyrinth'
\cite{Genn76} ---
has received much attention both in the physics and the mathematics literature.
Recently, several papers have considered random walk on a supercritical percolation
cluster \cite{Barl04,BB06,MP06,SS04}.  Roughly speaking, supercritical percolation
clusters on $\Zd$ are $d$-dimensional, and these papers
prove, in various ways, that a random walk on a supercritical percolation cluster
behaves in a diffusive fashion similar to a random walk on the entire lattice $\Zd$.

Although a mathematically rigorous understanding of \emph{critical} percolation
clusters is restricted to examples in dimensions $d=2$ and $d>6$, or $d>4$ in the
case of \emph{oriented} percolation, it is generally believed that critical percolation
clusters in dimension $d$ have dimension less than $d$, and that random walk on a
large critical cluster behaves subdiffusively.  Critical percolation clusters are
believed to be finite in all dimensions, and are known to be finite in the oriented
setting \cite{BG90}.
To avoid finite-size issues associated
with random walk on a finite cluster, it is convenient to consider random walk
on the incipient infinite cluster (\IIC), which can be understood as a critical percolation
cluster conditioned to be infinite.  The \IIC\ has been constructed
so far only when $d=2$ \cite{Kest86}, when $d>6$ (in the spread-out case) \cite{HJ04},
and when $d>4$ for oriented percolation on $\Zd\times\Z_+$
(again in the spread-out case) \cite{HHS02}. See \cite{Slad06} for a summary of
the high-dimensional results.
Also, it is not difficult to construct the \IIC\ on a tree
\cite{BK06,Kest86a}.

Random walk on the \IIC\ has been proved to be subdiffusive on $\Z^2$
 \cite{Kest86a} and on a tree \cite{BK06,Kest86a}.
See also \cite{DC06a,DC06b} for related results in the continuum limit.
In this paper, we prove several estimates for random walk on the
\IIC\ for spread-out
oriented percolation on $\Zd \times \Z_+$ in dimensions $d>6$.
These estimates, which show subdiffusive behaviour, establish that the spectral
dimension of the \IIC\ is $\frac 43$, thereby
proving the Alexander--Orbach \cite{AO82}
conjecture in this setting. For random walk on ordinary (unoriented) percolation
for $d < 6$ the Alexander--Orbach conjecture is generally believed
to be false \cite[Section~7.4]{Hugh96}.

The upper critical dimension for oriented percolation is $4$.
Because of this, we initially expected that the spectral dimension
of the \IIC\ would be equal to $\frac 43$ for oriented percolation
in all dimensions $d>4$, but not for $d < 4$.  However, our
methods require that we take $d>6$.  The random walk is allowed to
travel backwards in `time' ({as measured by the oriented
percolation process}), and this allows the
walk to move between vertices that are not connected
{to each other} in the
oriented sense.  It may be that this effect raises the upper
critical dimension for the random walk in the oriented setting to
$d=6$.  Or it may be that our conclusions for the random walk remain true for
all dimensions $d>4$, despite the fact that our methods force us
to assume $d>6$.  This leads to the open question: Do our results
actually apply in all dimensions $d>4$, or does different
behaviour apply for $4 < d \leq 6$?

\subsection{Random walk on graphs
{and in random environments}}\label{sec-rwre}

{Our results on the \IIC\ will be consequences of
more general results on random walks on a family of random graphs.
We now set up our notation for this. }
Let $\Gamma = (G,E)$ be an infinite graph, with vertex set $G$
and edge set $E$.  The edges $e \in E$ are \emph{not} oriented.
We assume that $\Gam$ is connected.
We write $x\sim y$ if $\{x,y \} \in E$, and
assume that $(G, E)$ is locally finite,
i.e., $\mu_y<\infty$ for each $y\in G$,  where $\mu_y$ is the number of
bonds that contain $y$. We extend $\mu$ to a measure on $G$.
Let $X =(X_n, n  \in  \Z_+, P^x, x \in G)$ be the discrete-time
simple random walk on $\Gamma$, i.e.,
{the Markov chain with transition probabilities}
\eq
    P^x(X_1=y) = \frac{1}{\mu_x}, \quad y \sim x.
\en
We define the transition  density (or discrete-time heat kernel) of $X$ by
\eq
     p_n(x,y)= \frac{P^x(X_n=y)}{\mu_y};
\en
we have $p_n(x,y)=p_n(y,x)$.

The natural metric on $\Gamma$, obtained by counting the number of
steps in the shortest path between points, is written $d(x,y)$ for
$x,y\in G$. We write
\eq
     B(x,r)=\{y: d(x,y) <  r\}, \qquad V(x,r)=\mu(B(x,r)), \quad r \in (0,\infty).
\en
Following terminology used for manifolds, we call $V(x,r)$ the \emph{volume}
of the ball $B(x,r)$.
We will assume $G$ contains a marked vertex, which we denote $0$,
and we write
 \eq
   B(R) = B(0,R), \qquad V(R)= V(0,R).
\en
For $A \subset G$, we write
\eq
\lbeq{TAdef}
     T_A =\inf\{n\ge 0: X_n \in A\}, \qquad \tau_A = T_{A^c},
\en
and let
\eq
\tau_R = \tau_{B(0,R)} =\min\{ n \ge 0: X_n \not\in B(0,R)\}.
\en
Let $W_n = \{X_0,X_1,\ldots, X_n\}$ {be the set of vertices hit by
$X$ up to time $n$,} and let
\eq
  S_n = \mu(W_n) = \sum_{x \in W_n} \mu_x.
\en

We write $\Reff(0,B(R)^c)$ for the effective resistance between
$0$ and $B(R)^c$ in the electric network obtained by making each edge
of $\Gam$ a unit resistor -- see \cite{DS84}.
A precise mathematical definition of $\Reff(\cdot, \cdot)$
will be given in Section \ref{sec:rwresults}.

\medskip

We now consider a probability space $(\Omega, \Fcal, \bP)$ carrying a
family of random graphs $\Gam(\omega)=(G(\omega),E(\omega), \omega\in \Omega)$.
We assume
that, for each $\omega \in \Omega$, the graph $\Gam(\omega)$ is infinite,
locally finite and connected, and contains a marked vertex $0\in G$.
We denote balls in $\Gam(\omega)$ by  $B_\omega(x,r)$, their volume
by $V_\omega(x,r)$, and write
\eq
    B(R) = B_\omega(R) = B_\omega(0,R),
    \quad\quad
    V(R) = V_\omega(R) = V_\omega(0,R).
\en
We write
$X=(X_n , n \ge 0, P_\omega^x, x \in G(\omega))$ for the
simple random walk on $\Gam(\omega)$, and denote
by $p_n^\omega(x,y)$ its transition density with respect
to $\mu(\omega)$.
Formally, we  introduce a second
measure space $(\overline \Omega, \overline \Fcal)$,
and define $X$ on the product $\Omega \times \overline \Omega$.
We write
$\overline \omega$ to denote elements of $\overline \Omega$.

The key ingredients in our analysis of the simple random walk are volume
and resistance bounds.  The following defines a set $J(\lambda)$
of values of $R$ for which we have `good' volume and effective resistance
estimates.
The set $J(\lambda)$ depends on the graph $\Gamma$, and thus is a random set
under $\bP$.

\begin{defn}\label{jdef}
\emph{ Let $\Gamma=(G,E)$ be as above.
For $\lambda>1$, let $J(\lambda)$ be the set of
those $R \in [1,\infty]$ such that the following all hold: \\
(1) $V(R) \le \lambda R^2$,\\
(2) $V(R) \ge \lambda^{-1} R^2$,\\
(3) $\Reff (0, B(R)^c)\ge \lambda^{-1} R$.}
\end{defn}

Note that $\Reff(0, B(R)^c) \le R$ {(see Lemma~\ref{reffacts}(c)
in Section~\ref{sec-gg} )}, so there is no need for an upper
bound complementary to Definition~\ref{jdef}(3).  We now make the
following important
assumption concerning the graphs $(\Gam(\omega))$.  This
involves upper and lower bounds on the volume, as well as an estimate
which says that $R$ is likely to be in $J(\lambda)$ for large enough
$\lambda$.

\begin{ass}\label{ass-rwre}
{There exists $R^*\ge 1$ such that the following hold:}\\
(1) There exists $p(\lambda)\ge 0$,
{with $p(\lambda)\le {c_1}{\lam^{-q_0}}$ for some
$q_0, c_1>0$}, 
such that for each $R\ge R^*$,
\begin{equation}\label{plamassump}
\Pbold (R\in J(\lambda))\ge 1-p(\lambda),
\end{equation}
(2) $\Ebold [V(R)]\le c_{2} R^2$, for $R \in [R^*, \infty)$, \\
(3) $\Ebold [1/V(R)]\le c_{3} R^{-2}$ for $R \in [R^*, \infty)$. \\
\end{ass}

\begin{rem}
Assumption~\ref{ass-rwre}(2,3), together with Markov's inequality,
provides upper bounds of the form $c\lambda^{-1}$ for the probability of the complements
of the events in Definition~\ref{jdef}(1,2). 
This creates {some } redundancy
in our formulation, but {we state things this way}
because some of our conclusions for the
random walk rely only on
Assumption~\ref{ass-rwre}(1) and do not require the stronger volume bounds
given by Assumption~\ref{ass-rwre}(2,3).

Note that Assumption~\ref{ass-rwre} only involves statements about the volume
and resistance from one point $0$ in the graph. In general, this kind of
information would not be enough to give much control of the random walk.
However, the graphs considered here have strong recurrence properties,
and are therefore simpler to handle than general graphs.
We use techniques developed in \cite{BCK05, BK06,T89, T01a, T01b}.

We will prove in Theorem \ref{veriass} that Assumption~\ref{ass-rwre} holds for
the \IIC\ {for sufficiently spread-out oriented percolation on $\Z^d \times \Z_+$
when $d>6$}.
As the reader of Sections~\ref{sec:laceexp}--\ref{sec-piv} will see,
obtaining volume and (especially) resistance bounds on the \IIC\ from one base
point is already difficult; it is fortunate that we {do not need to assume more}.
\end{rem}

{We have the following four consequences of
Assumption~\ref{ass-rwre} for random graphs. They
give control, in different ways, of  the quantities
$E^0_\omega \tau_R$, $p_{2n}(0,0)$, 
$d(0,X_n)$}, and {$S_n$}, 
which measure the rate of dispersion of the random
walk $X$ from the base point $0$. 
Some statements in the first proposition involve the {averaged} law
defined by the semi-direct product $P^* = \bP\times  P^0_\omega$.

\begin{theorem} \label{ptight}
Suppose Assumption \ref{ass-rwre}(1) holds. 
Then,
uniformly with respect to $n\ge 1$ and $R\ge 1$,
\begin{align}
\label{pt-a}
 \Pbold (\theta^{-1}\le R^{-3}E^0_\omega\tau_R\le \theta) &\to 1
\quad \text{ as } \theta\to \infty, \\
\label{pt-b}
  \Pbold (\theta^{-1}\le n^{2/3}p_{2n}^\omega(0,0)\le \theta) &\to 1
\quad \text{ as } \theta\to \infty, \\
\label{pt-c}
  P^*( d(0,X_n)n^{-1/3} < \theta) &\to 1
\quad \text{ as } \theta\to \infty. \\
\label{pt-d}
 P^*( \theta^{-1} < (1+d(0,X_n))n^{-1/3} ) &\to 1
\quad \text{ as } \theta\to \infty.
\end{align}
\end{theorem}

\smallskip
Since $P^0_\omega(X_{2n}=0) \approx n^{-2/3}$, we cannot replace
$1+d(0,X_n)$ by $d(0,X_n)$ in (\ref{pt-d}).

\begin{theorem} \label{pmeans}
Suppose Assumption \ref{ass-rwre}(1,2,3) hold. Then there
exists $n^* \ge 1$ (depending only on $R^*$ and the function $p(\cdot)$
in Assumption~\ref{ass-rwre}), and constants $c_i$ such that
\eqalign
\label{e:mmean}
    c_1R^3 \le \Ebold (E^0_\omega\tau_R) \le c_2R^3 & \text{  for all } R\ge 1,
\\
\label{e:pmean}
 c_3 n^{-2/3}\le \Ebold (p_{2n}^\omega(0,0))\le c_4 n^{-2/3} &\text{  for all } n \ge n^*.
\\
\label{e:dmean}
       c_5 n^{1/3}\le \Ebold (E_\omega^0 d(0,X_n)) &\text{  for all } n\ge n^*.
\enalign
\end{theorem}

We do not have an upper bound in (\ref{e:dmean}); this is discussed further
in Example~\ref{ex-p1} below.

\begin{rem} The above two theorems in fact do not require the
polynomial decay of $p(\lambda)$; it is enough to
have $p(\lambda)\to 0$ as $\lambda\to \infty$.
\end{rem}

{
Let $d_s(G)$ be the {\em spectral dimension of} $G$, defined by
\eq
\label{e:dsdef}
 d_s(G) = -2 \lim_{n \to \infty} \frac{\log p_{2n}(x,x)}{\log n},
\en
if this limit exists. Here $x\in G$; it is easy to see that the limit
is independent of the base point $x$.
Note that $d_s(\Z^d)=d$.
}

In (c) below, recall that  $\overline \Omega$ is the second probability
space, on which the random walk $X$ is defined.

\begin{theorem}
\label{thm-rwre} Suppose Assumption \ref{ass-rwre}(1) holds.
Then there exist $\al_1, \al_2, \al_3, \al_4 < \infty$, and a subset $\Omega_0$
with $\bP(\Omega_0)=1$ such that the following statements hold.\\
(a) For each $\omega \in \Omega_0$ and $x \in G(\omega)$ there
exists $N_x(\omega)< \infty$ such that
\eq
 \label{e:logpnlima}
   (\log n)^{-\al_1} n^{-2/3} \le
 p^\omega_{2n}(x,x) \le (\log n)^{\al_1}  n^{-2/3}, \quad  n\ge N_x(\omega).
\en
In particular, $d_s(G) = \frac43$, $\bP$-a.s.,
and the random walk is recurrent. \\
(b)  For each $\omega \in \Omega_0$ and $x \in G(\omega)$ there
exists $R_x(\omega)< \infty$ such that
\eq
 \label{e:logtaulima}
   (\log R)^{-\al_2} R^3 \le E^x_\omega \tau_R
 \le  (\log R)^{\al_2} R^3, \quad  R\ge R_x(\omega).
\en
Hence
\[ \lim_{R \to \infty} \frac{\log  E^x_\omega \tau_R}{\log R} =3.\]
(c)  Let $Y_n= \max_{0\le k \le n} d(0,X_k)$.
For each $\omega \in \Omega_0$ and $x \in G(\omega)$
there exist $N_x(\omega,\overline \omega), R_x(\omega, \overline \omega)$
such that $P^x_\omega( N_x <\infty)=P^x_\omega( R_x <\infty)=1$,
and such that
\begin{align}
 \label{e:ynlim}
 (\log n)^{-\al_3} n^{1/3} \le Y_n(\omega, \overline \omega)
 &\le (\log n)^{\al_3}  n^{1/3},
 \quad n \geq N_x(\omega, \overline \omega), \\
   \label{e:rnlim}
 (\log R)^{-\al_4} R^3 \le \tau_R(\omega, \overline \omega)
 &\le (\log R)^{\al_4}  R^3, \quad\quad
 R \geq R_x(\omega, \overline \omega).
\end{align}
\end{theorem}

{
\begin{rem}
One cannot expect (\ref{e:logpnlima}) or (\ref{e:logtaulima})
to hold with $\al_1=0$ or $\al_2=0$, since it is known that
$\log \log$ fluctuations occur  in the analogous limits
for the \IIC\ on regular trees \cite{BK06}. (This example is discussed
further in Example~\ref{ex-treeIIC}(i) below).
\end{rem}
}

{Let $W_n = \{X_0,X_1,\ldots, X_n\}$ as before and let
$|W_n|$ denote its cardinality.
For a sufficiently regular recurrent graph one expects that
$|W_n| \approx n^{d_s/2}$.
The original formulation} of the Alexander-Orbach conjecture \cite{AO82}
was that, in all dimensions, for the \IIC,
\eq
\label{e:snaoz}
|W_n| \approx n^{2/3},
\en
so that $d_s = \frac 43$ in all dimensions.
As noted already above, the conjecture is now not believed to hold in low
dimensions.
The following theorem shows that a version of the Alexander--Orbach
conjecture does hold for random graphs that
satisfy  Assumption \ref{ass-rwre}(1).
As we will see in Theorem \ref{veriass}, this is the
case for the \IIC\ for sufficiently spread-out oriented percolation on
$\Zd \times \Z_+$ for $d>6$.

\begin{theorem}
\label{thm-snlim} (a)
Suppose Assumption \ref{ass-rwre}(1) holds. Then
there exists a subset $\Omega_0$ with $\bP(\Omega_0)=1$ such that
for each $\omega \in \Omega_0$ and $x \in G(\omega)$,
\eq
\label{e:snlim}
 \lim_{n \to \infty} \frac{\log S_n}{\log n} = \frac23, \quad
 P^x_\omega \text{-a.s.}.
\en
{
(b) Suppose in addition there exists a constant $c_0$ such that
all vertices in $G$ have degree less than $c_0$. Then
\eq \lbeq{Wnlim}
\lim_{n \to \infty} \frac{\log |W_n|}{\log n} = \frac23, \quad
 P^x_\omega \text{-a.s.}
\en
}
\end{theorem}

See Example \ref{ex-treeIIC} for a graph with unbounded degree which
satisfies Assumption \ref{ass-rwre}, but for which \eqref{e:Wnlim}
fails.

\begin{rem} See \cite{KM07} for results which generalise the
above theorems
to the situation where there exist indices $\alpha< \beta$ such that
$V(R)$ is comparable to $R^\alpha$ and $\Reff(0, B(R)^c)$ is comparable
to $R^{\beta-\al}$. Our case is $\alpha=2$, $\beta =3$.
\end{rem}

\subsection{ {The IIC}}
\label{sec-iic}

In this section, we define the oriented percolation model and recall
the construction of the \IIC\ for spread-out oriented percolation on
$\Zd \times \Z_+$ in dimensions $d>4$ \cite{HHS02}.  For simplicity,
we will consider only the most basic example of a spread-out model.
(In the physics literature, oriented percolation is usually called
\emph{directed} percolation; see \cite{JT05}.)

The spread-out oriented percolation model is defined as follows.
Consider the graph with vertices $\Zd \times {\Zbold}_+$ and
directed bonds $((x,n),(y,n+1))$, for $n \geq 0$ and $x,y \in \Zd$
with $0 \leq \|x-y\|_\infty \leq L$.  Here $L$ is a fixed positive integer
and $\| x\|_\infty = \max_{i=1,\ldots, d} |x_i|$ for $x=(x_1,\ldots, x_d)\in \Zd$.
Let $p\in [0,1]$.
We associate to each directed bond $((x,n),(y,n+1))$ an
independent random variable taking the value $1$ with probability
$p$ and $0$ with probability $1-p$.   We say a
bond is {\em occupied}\/ when the corresponding random variable
is $1$, and {\em vacant}\/ when the random variable
is $0$.
Given a configuration of occupied bonds, we say that $(x,n)$ is {\em connected
to} $(y,m)$, and write $(x,n) \conn (y,m)$, if there is an oriented path
from $(x,n)$ to $(y,m)$ consisting of occupied bonds, or
if $(x,n)=(y,m)$.  Let $C(x,n)$ denote the forward cluster of $(x,n)$, i.e.,
$C(x,n) = \{(y,m) : (x,n) \conn (y,m)\}$, and let $|C(x,n)|$ denote
its cardinality.

The joint probability distribution of the bond variables
will be denoted $\Pbold$, with corresponding expectation denoted
$\Ebold$; these depend on $p$ and are defined
on a probability space $(\Omega, \tilde \Fcal, \Pbold)$.
  Let $\theta(p) = \Pbold(|C(0,0)|=\infty)$.
For all dimensions $d \geq 1$ and for all $L \geq 1$,
there is a critical value $p_c=p_c(d,L) \in (0,1)$ such that
$\theta(p) = 0$ for $p \leq p_c$ and $\theta(p)>0$ for $p>p_c$.
In particular, there is no infinite
cluster when $p = p_c$ \cite{BG90,GH02}.
For the remainder of this paper, we fix $p=p_c$, so that
$\Pbold = \Pbold_{p_c}$.

To define the \IIC, some terminology is required.
A {\em cylinder event}\/ is an event that is determined by the occupation
status of a finite set of bonds.  We denote the algebra of cylinder events by
$\Fcal_0$, and define $\Fcal$ to be the $\sigma$-algebra generated by $\Fcal_0$.
The most natural definition of the \IIC\  is as follows.  Let $\{(x,m) \conn n \}$
denote the event that there exists $(y,n)$ such that $(x,m) \conn (y,n)$.
Let
\eq
\lbeq{Qndef}
    \Qbold_n (E) = \Pbold (E | (0,0) \conn n)
    \quad (E \in \Fcal_0)
\en
and define the \IIC\  by
\eq
\lbeq{Qdef}
    \Qbold_\infty (E) = \lim_{n \to \infty} \Qbold_n(E)
    \quad (E \in \Fcal_0),
\en
assuming the limit exists.
{
A possible alternate definition of the \IIC\ is to define
    \eq
    \lbeq{Pndef}
    {\mathbb P}_n(E) = \frac{1}{\tau_n}
    \sum_{x\in \Z^d} {\mathbb P}(E \cap \{ (0,0)\conn (x,n)\})
    \quad (E \in \Fcal_0)
    \en
with $\tau_n = \sum_{x\in \Z^d}{\mathbb P}( (0,0)\conn (x,n))$,
and to let
    \eq
    \lbeq{IICdef}
    {\mathbb P}_{\infty}(E)=\lim_{n\rightarrow \infty} {\mathbb P}_n(E)
    \quad (E \in \Fcal_0),
    \en
assuming the limit exists.
}

Let $d+1>4+1$ and $p = p_c$.
{
It was proved in \cite{HHS02} that
there is an $L_0=L_0(d)$ such that
for $L \geq L_0$ the limit \refeq{IICdef} exists for every cylinder
event $E \in \Fcal_0$.  Moreover, $\Pbold_\infty$ extends to a
probability measure on the $\sigma$-algebra $\Fcal$, and,
writing $\Ccal=C(0,0)$,
$\Ccal$ is ${\mathbb P}_{\infty}$-a.s.\ an infinite cluster.
It was also proved in \cite{HHS02} that if the critical survival probability
$\Pbold((0,0) \conn n)$ is asymptotic to a multiple of $n^{-1}$
as $n \to \infty$, then for $L_0=L_0(d)$ the limit \refeq{Qdef} exists
and defines a probability measure on $\Fcal$, and moreover $\Qbold_\infty
=\Pbold_\infty$ so both constructions yield the same measure.
Subsequently, it was shown in \cite{HHS05a,HHS05b}  that the survival
probability is indeed asymptotic to a multiple of $n^{-1}$ when $d+1>4+1$
and $L \geq L_0(d)$.  We will find both of
the equivalent definitions
\refeq{Qdef} and \refeq{IICdef} to be useful.}

We call $(\Ccal, \Qbold_\infty){=(\Ccal, \Pbold_\infty)}$
the \IIC, and this provides the random
environment for our random walk.  We write $\Ebold_\infty$ for
expectation with respect to $\Qbold_\infty$.
It will be convenient to remove a $\Qbold_\infty$-null set $\Ncal$
from the configuration space $\Omega$, so that for all
$\omega \in \Omega_0=\Omega- \Ncal$ the cluster
$\Ccal(\omega)$ is infinite (and connected).
{
The \IIC\ $\Ccal(\omega)$, $\omega \in \Omega$ under the law
$\Qbold_\infty$ gives a family of random graphs, with marked vertex
$\zerovec=(0,0)$,
so as in Section
\ref{sec-rwre} we can define a random walk
$X=(X_j, j  \in  \Z_+, P^{(x,n)}_\omega, (x,n) \in \Ccal(\omega))$.
Note that although the orientation is used to construct the cluster $\Ccal$,
once $\Ccal$ has been determined
the random walk on $\Ccal$ can move in any direction --- see Figure~\reffg{conn}.
}

\begin{figure}
\begin{center}
\setlength{\unitlength}{0.0075in}
\begin{picture}(215,200)(20,580)
\thicklines

\qbezier(110,580)(110,680)(110,720)
\qbezier(110,600)(140,640)(160,720)
\qbezier(110,650)(140,700)(160,720)

\put(130,612){\makebox(0,0)[lb]{\raisebox{0pt}[0pt][0pt]{${\sss (x,m)}$}}}
\put(70,665){\makebox(0,0)[lb]{\raisebox{0pt}[0pt][0pt]{${\sss (y,n)}$}}}

\put(130,632){\makebox(0,0)[lb]{\raisebox{0pt}[0pt][0pt]{\circle*{8}}}}
\put(120,665){\makebox(0,0)[lb]{\raisebox{0pt}[0pt][0pt]{\circle*{8}}}}
\end{picture}
\end{center}

\caption{\lbfg{conn}
Although the vertex $(x,m)$ is not connected to $(y,n)$, or vice versa, in the sense
of oriented percolation
{(oriented upwards)}, it is nevertheless possible for a random walk
to move from one of these vertices to the other.
}
\end{figure}
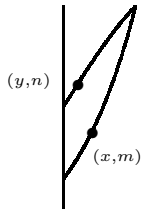

\begin{theorem}
\label{veriass}
 For $d>6$, there is an $L_1=L_1(d) \geq L_0(d)$ such
that for all $L \geq L_1$, Assumption \ref{ass-rwre}(1)--(3) hold
with $q_0=1$ and constants $c_1, c_2, c_3$ independent of $d$ and $L$.
{Consequently, the conclusions of Theorems~\ref{ptight}, \ref{pmeans},
\ref{thm-rwre} and \ref{thm-snlim} all hold for the
random walk on the \IIC.}  {In particular, the Alexander--Orbach conjecture holds
in the form of \refeq{Wnlim}.}

\end{theorem}

As we will see later, the restriction to $d>6$ is required
only for our estimate of the effective resistance.

{
\begin{rem}
Since the constants in Assumption~\ref{ass-rwre} are independent of
$d, L$ for the \IIC\
(provided $d>6$ and $L\ge L_1(d)$), the constants $\al_1, \cdots, \al_4$
in Theorem \ref{thm-rwre} are also independent of $d$ and $L$
when applied to the \IIC.
\end{rem}
}

{The proof of our main results are performed in two principal steps,
corresponding to the results in Section \ref{sec-rwre} and
Theorem~\ref{veriass} respectively.

The results in Section \ref{sec-rwre} are proved in Section
\ref{sec:rwresults}. The first step is to obtain
estimates for a fixed (non-random) graph $\Gam$.
In Section \ref{sec-gg}, using arguments
based on those in \cite{BCK05} and \cite{BK06}, we show that
volume and resistance bounds on $\Gam$ lead to bounds on
transition probabilities and hitting times. Then, in Section
\ref{sec-genrwre} we translate these results into the
random graph context, and prove Theorems~\ref{ptight}--\ref{thm-snlim}.}

{
The second step is the proof of Theorem~\ref{veriass}.
Section \ref{sec:assrw} states three properties of the \IIC\
for critical spread-out oriented percolation in dimensions $d>6$,
and show that these imply Theorem~\ref{veriass}.
These properties are proved in
Sections~\ref{sec:laceexp}--\ref{sec-piv},
using an extension of results of
\cite{HHS02,HHS05a,HS03b}
that were obtained using the lace expansion.
}

{
\subsection{ Further Examples}

We have some other examples of random graphs which satisfy
Assumption~\ref{ass-rwre}.
}

\begin{example} \label{ex-treeIIC}
(i)
Assumption~\ref{ass-rwre} holds for random walk on the \IIC\ for the binomial tree;
see  \cite[Corollary~2.12]{BK06}.  Therefore the conclusions of
Theorems~\ref{ptight}--\ref{thm-snlim}
hold for random walk on this \IIC.
The results of \cite{BK06} go beyond Theorem~\ref{thm-rwre}(a) and (b) in this context,
but Theorem~\ref{thm-rwre}(c) and Theorem~\ref{thm-snlim} here are new.

\medskip\noindent
(ii) It is shown in \cite{AGHS07} that the invasion percolation
cluster on a regular tree is stochastically dominated by the \IIC\ for
the binomial tree.  Consequently, upper bounds on the volume and lower
bounds on the effective resistance of the invasion percolation cluster
follow from the corresponding bounds for the \IIC\
{(using Lemma~\ref{reffacts}(e) in Section~\ref{sec-gg})}.
Assumption~\ref{ass-rwre}(1,2) for the
invasion percolation cluster therefore follows from its counterpart
for the \IIC\ for the binomial tree.
{In addition, the lower bound on the volume in
Assumption~\ref{ass-rwre}(3) is proved for the invasion percolation
cluster in \cite{AGHS07}.   Therefore}
Assumption~\ref{ass-rwre}
holds for the invasion percolation cluster on a regular tree, and
hence simple random walk on the invasion percolation cluster also
obeys the conclusions of Theorems~\ref{ptight}--\ref{thm-snlim}.
See \cite{AGHS07} for further details about {this example}.

\medskip\noindent
{(iii)}
Consider the incipient infinite branching random walk
({\hbox{\footnotesize\rm IIBRW}}), obtained as the limit as
$n\to\infty$ of critical branching random walk (say with binomial
offspring distribution) conditioned to survive to at least $n$
generations \cite[Section~2]{Hofs06}.  We interpret the
{\hbox{\footnotesize\rm IIBRW}} as a random infinite subgraph of
$\Z^d\times \Z_+$.  There is the option of considering either one edge
per particle jump, leading to the occurrence of multiple edges between
vertices, or identifying any such multiple edges as a single edge; we
believe both options will behave similarly in dimensions $d>4$.
Consider simple random walk on the {\hbox{\footnotesize\rm IIBRW}}.
Our volume estimates for the \IIC\ for oriented percolation for $d>4$
will adapt to give similar estimates for the {\hbox{\footnotesize\rm
IIBRW}} for $d>4$.  The effective resistance $\Reff(0,B(R)^c)$ for the
{\hbox{\footnotesize\rm IIBRW}} is lower than it is for the \IIC\ on a
tree, due to cycles in the {\hbox{\footnotesize\rm IIBRW}}.  It is an
interesting open problem to obtain a lower bound on $\Reff(0,B(R)^c)$
for the {\hbox{\footnotesize\rm IIBRW}}, to establish
Assumption~\ref{ass-rwre} and hence its consequences
Theorems~\ref{ptight}--\ref{thm-snlim}
for random walk on the
{\hbox{\footnotesize\rm IIBRW}}.  Our main interest is the question:
Does random walk on the {\hbox{\footnotesize\rm IIBRW}} have the same
behaviour in all dimensions $d>4$, or is there different behaviour for
$4<d\leq 6$ and $d>6$?  An answer would shed light on the question
raised at the end of Section~\ref{sec-intro}.  It would also be of
interest to consider this question in the continuum limit: Brownian
motion on the canonical measure of super-Brownian motion conditioned
to survive for all time (see \cite{Hofs06}).

\medskip\noindent
{(iv) A {non-random} graph $\Gam$ satisfies Assumption~\ref{ass-rwre}
if and only if there exists $\lam$ such that $J(\lam)=[1,\infty)$.
If $\Gam_i$, $1\le i \le n$ are graphs satisfying Assumption~\ref{ass-rwre}
then the graph $\Gam$ obtained by joining the $\Gam_i$ at their
marked vertices also satisfies Assumption~\ref{ass-rwre}.

\medskip\noindent (v)
Consider the {non-random} graph consisting of $\bZ_+$ with for each $n$
a finite subgraph $G_n$ connected by one point in $G_n$ to the vertex
$n$. If $\mu(G_n) \asymp n$ and the diameter of
$G_n$ is $o(n)$ then Assumption~\ref{ass-rwre} holds.
In particular, if we take $G_n$ to be the complete graph with
$r_n = \lfloor n^{1/2} \rfloor$ vertices, then while
$V(R) \asymp R^2$, we have $|B(R)| \asymp R^{3/2}$.  In this case
\eqref{e:snlim} holds, whereas
\eq \lbeq{Wnlim2}
\lim_{n \to \infty} \frac{\log |W_n|}{\log n} = \frac12, \quad
 P^x_\omega \text{-a.s.}
\en
The rough idea behind \refeq{Wnlim2} is as follows.
By \refeq{ynlim}, the distance travelled up to time $n$ is approximately $n^{1/3}$.
The proof of Theorem \ref{thm-snlim} shows that
the random walk will visit a positive fraction of the vertices within this distance,
and there are of order $(n^{1/3})^{3/2}=n^{1/2}$ such vertices,
leading to \refeq{Wnlim2}.}
This shows that some bound on vertex degree is necessary before
one can pass from \eqref{e:snlim} to \eqref{e:Wnlim}.
\end{example}

\medskip

Throughout the paper, we use $c, c'$ to denote strictly positive
finite constants whose values are not significant and may change
from line to line. We write $c_i$ for positive constants whose
values are fixed within theorems and lemmas.

\section{Random walk on a random graph}
\label{sec:rwresults}

{
In this section we prove Theorems~\ref{ptight}--\ref{thm-snlim}.
First, in Section~\ref{sec-gg}, we study the random walk on a fixed graph;
then, in Section~\ref{sec-genrwre} we apply these results to
a family of random graphs satisfying Assumption~\ref{ass-rwre}.
}

\subsection{{ Random walk on a fixed graph }}
\label{sec-gg}

In this section, we fix an infinite locally-finite connected
graph $\Gam = (G,E)$, and will show that
bounds on the quantities $V(R)$ and $\Reff(0, B(R)^c)$ lead to
control of $E^0 \tau_R$, $p_n(0,0)$ and $E^0 d(0,X_n)$.
The results in \cite{BCK05} (see \cite[Theorem~1.3, Lemma~2.2]{BCK05})
cover the case where, for all $x \in G$ and $R\ge 1$,
\begin{equation}\label{vrreg}
c_1 R^2 \le V(x,R) \le c_2 R^2, \quad
c_3 R \le \Reff(x, B(x,R)^c) \le c_4 R.
\end{equation}
Here, we treat the case where we only have information available on
the volume and effective resistance from one fixed point $0$ in the
graph, and only for certain values of $R$.
Our methods are very close to those of
\cite{BCK05}, but the need to keep track of the values of $R$ for
which we make use of the bounds makes the details of the proofs
more complicated.

The following Proposition gives the majority of the bounds on
$\tau_R$, $p_n(0,0)$ and $d(0,X_n)$ that will be used in
Section~\ref{sec-genrwre}.

Recall the definition of $J(\lam)$ from Definition~\ref{jdef}.
In the following proposition, we will take $\lam \ge 1$ and assume
that $R$, and certain multiples of $R$, are in $J(\lam)$.
We then obtain (for example) bounds on $E^0 \tau_R$; these bounds
will involve constants depending on $\lam$.
{For the limit Theorems \ref{thm-rwre} and \ref{thm-snlim}}
we need to know that the
dependence of these constants on $\lam$ is polynomial in $\lam$.
To indicate this, we write $C_i(\lam)$ to denote
positive constants of the form $C_i(\lam) = C_i \lam^{\pm q_i}$,
which will be fixed throughout this section.
The sign accompanying $q_i>0$ is such that statements become weaker as $\lam$
increases.

\begin{prop}\label{rw-main} Let $\lam\ge 1$.
There exist $C_1(\lam),\cdots,C_9(\lam)$
such that the following hold.\\
(a) Suppose that $R \in J(\lam)$. Then
\eq
\label{Etaub}
    E^x \tau_R \le 2 \lam R^3 \quad\text{ for } x\in B(R).\\
\en
{Suppose that $R, R/(4\lam) \in J(\lam)$.}
Then
\begin{align}
\label{e:tmlb}
   E^x \tau_R &\ge C_1(\lam) R^3, \text { for } x \in  B(0, R/{(4\lam)}).
\end{align}
Let $\eps < 1/(4\lam)$ and
$R, \eps R, \eps R/(4\lambda)\in J(\lambda)$. Then
\eq
 \label{trdlb}
 P^y \big( \tau_R \le C_2(\lam) (\eps R)^3 \big) \le C_3(\lam) \eps,
 \quad \text{ for } y \in B(\eps R).
\en
(b)  Suppose that $R \in J(\lam)$. Then
\eq \label{hkgub}
{p_{n}(0,y)+p_{n+1}(0,y)  \le C_4(\lam) n^{-2/3}}  \quad\text{ for }
 y \in B(R) \text{ if } n=2 \lfloor R \rfloor ^3.
\en
Suppose that $R, R/(4\lam) \in J(\lam)$.
Then
\eq \label{e:plb}
 p_{2n}(x,x) \ge C_5(\lam) n^{-2/3}
\quad \quad \text{ for }
\tfrac14 C_1(\lam)R^3 \le   n \le \tfrac12 C_1(\lam) R^3,
 \quad  x \in  B(0, R/{(4\lam)}).
\en
(c) Let $n\ge 1$, $M\ge 1$, and set $R=M n^{1/3}$.
If $R, C_6(\lam)R/M, C_6(\lam)R/(4\lambda M) \in J(\lambda)$, then
\eq  \label{doxub}
P^0\big(n^{-1/3}d(0,X_n)>M\big)\le \frac{C_7(\lam)} {M}. \en
We have $C_7(\lam) \le c \lam^{22/3}$.\\
(d) Let $R=(n/2)^{1/3}$ and $M \ge 1$. If $R, R/M \in J(\lam)$ then
\eq \label{doxlb}
P^0\big( d(0,X_n) < R/M \big) \le \frac{\lam C_4(\lam)}{M^2}.
\en
Also, if  $R, C_8(\lam) R \in J(\lam)$ then
\eq
\lbeq{Edlb}
{   E^0 d(0,X_n) \ge C_9(\lam) n^{1/3}.}
\en
\end{prop}

The overall strategy for the proof of these various inequalities is as follows.
We begin with obtaining bounds on the mean exit time $E^0\tau_R$. Using the
Green function (see \eqref{greendef} below for the definition) we can write
\eq \label{txgreen1}
E^z\tau_B = \sum_{y\in B} g_B(z,y)\mu_y.
\en
Since $g_B(x,x) = \Reff(x, B^c)$ (see \eqref{eq:green2}), this leads
to the upper and lower bounds on  $E^x\tau_R$ for $x$ sufficiently
close to $0$ given in \eqref{Etaub} and \eqref{e:tmlb}.
The final inequality concerning $\tau_R$ is  \eqref{trdlb},
which bounds from above the lower tail of $\tau_R$. (This is equivalent
to bounding from above
the speed at which $X$ can move from its starting point $0$.)
The proof for this takes the bounds in  \eqref{Etaub} and \eqref{e:tmlb}
as its starting point, but also uses a simple inequality relating
effective resistance and hitting probabilities -- see Lemma \ref{hitest}
below.

The next set of inequalities we prove are those for the heat kernel
$p_n(x,y)$. In the continuous time setting these are proved using differential
inequalities which relate the derivative of the heat kernel to its energy.
Unfortunately in discrete time the  differential
inequalities are replaced by rather less intuitive difference equations,
which in addition take a slightly more complicated form.  The estimate
\eqref{hkgub} is proved from an inequality which bounds the heat kernel
just in terms of the volume of balls -- see \eqref{hkub}. Adding information
on $\tau_R$ then enables one to obtain the lower bound \eqref{e:plb}.

The final bounds on $d(0,X_n)$ then follow easily from the bounds on
$\tau_R$ and $p_n(0,x)$.

\subsubsection{Bounds on $\tau_R$}

{We begin by giving a precise definition of effective resistance.
Let $\Ecal$ be the quadratic form given by}
\eq
      \Ecal(f,g)=\tfrac 12\sum_{\substack{x,y\in G \\ x \sim y}}
      (f(x)-f(y))(g(x)-g(y)),
\en
where $x\sim y$ means $\{x,y\}\in E$.  If we regard $\Gamma$ as an
electrical network with a unit resistor on each edge in $E$, then
$\Ecal(f,f)$ is the energy dissipation when the vertices of $G$ are at
a potential $f$.  Set $H^2=\{ f\in \Rbold^{G}: \Ecal(f,f)<\infty\}$.
Let $A,B$ be disjoint subsets of $G$.  The effective resistance
between $A$ and $B$ is defined by:
\begin{equation}
\label{3.3bk}
    \Reff(A,B)^{-1}=\inf\{\Ecal(f,f): f\in H^2, f|_A=1, f|_B=0\}.
\end{equation}
Let $\Reff(x,y)=\Reff(\{x\},\{y\})$, and $\Reff(x,x)=0$.
For general facts on  effective resistance and its connection with
random walks see \cite{AF09,DS84,LP09}.
We recall some basic properties of $\Reff(\cdot,\cdot)$.

\begin{lemma} \label{reffacts}
Let $\Gam=(G,E)$ be an infinite connected graph. \\
(a) $\Reff$ is a metric on $G$. \\
(b) If $A'\subset A$, $B'\subset B$ then
$\Reff(A',B') \ge \Reff(A,B)$. \\
(c) $\Reff(x,y) \le d(x,y)$. \\
(d) If $x, y \in G\setminus A$ then $\Reff(x,A)\le \Reff(x,y) + \Reff(y,A)$.\\
(e) If $\Gamma'=(G',E')$ is a subgraph of $\Gamma$, with effective resistance
$\Reff'$, and if $A'=A\cap G'$ and $B'=B\cap G'$,
then $\Reff'(A',B') \geq \Reff (A,B)$.\\
(f) For all $f\in \Rbold^{G}$ and $x,y\in G$,
\eq\label{resobin}
|f(x)-f(y)|^2\le \Reff(x,y)\Ecal(f,f).
\en
\end{lemma}

\proof
For (a) see \cite[Section~2.3]{Kiga01}. The monotonicity in (b) and (e)
is immediate from the variational definition of $\Reff$.
(c) is easy, and there is a proof in
\cite[Lemma~2.1]{BCK05}.  (d) follows from (a) by
considering the graph $\Gam'$ in which all vertices in $A$ are
connected by short circuits, which reduces $A$ to a single
vertex $a$.  \\
{
(f) If $f(x)=f(y)$ then \eqref{resobin} is immediate. If not,
then set $u(z)=(f(z)-f(y))/(f(x)-f(y))$, so that $u(x)=1$ and
$u(y)=0$. Then by  \eqref{3.3bk}
\[ \Reff(x,y)^{-1} \le \Ecal(u,u)= \Ecal(f,f) |f(x)-f(y)|^{-2}, \]
which gives  \eqref{resobin}. \qed
}

\medskip
{The inequality \eqref{resobin} will play an important
role in obtaining pointwise information on functions
from resistance or energy estimates.

Recall that $T_A$ was defined in \refeq{TAdef} to
be the hitting time of $A \subset G$.
If $A$ and $B$ are disjoint subsets of $G$ and $x \nin A \cup B$, then
(see \cite[Fact 2, p. 226]{BGP03})
\eq \label{bgp-lem}
P^x(T_A < T_B) \le \frac{\Reff(x,B)}{\Reff(x,A)}.
\en
}

{

\begin{lemma}\label{hitest}
Let $\lam\ge 1$ and suppose $R \in J(\lam)$. Let
$0< \eps \le 1/(2\lam)$, and $y \in B(\eps R)$. Then
\begin{align}
P^y ( T_0 < \tau_R ) &\ge 1- \frac{\lam \eps}{1- \eps \lam} \ge 1- 2 \eps \lam ,\\
P^0 ( T_y <  \tau_R ) &\ge  1- \eps \lam .
\end{align}
\end{lemma}

\proof By Lemma \ref{reffacts}(c) $\Reff(y,0)\le d(y,0)$,
while by  Lemma \ref{reffacts}(d) and the definition of
$J(\lam)$,
\[ \Reff(y,B(R)^c) \ge \Reff(0,B(R)^c) - \Reff(0,y)
 \ge \frac{R}{\lam} - \eps R. \]
So by \eqref{bgp-lem}
\[ P^y(\tau_R < T_0) \le  \frac{\Reff(y,0)}{\Reff(y, B(R)^c)}
\le \frac{\eps\lam}{1-\eps \lam}. \]
Similarly,
$P^0(\tau_R < T_y) \le  {\Reff(0,y)}/{\Reff(0, B(R)^c)}
 \le \eps \lam. $
\qed
}

\medskip

The initial steps in bounding $\tau_R$ use the Green kernel for the
random walk $X$, so we now recall its definition.
(These facts about Green functions will only be used in this subsubsection.)
Let $B\subset G$,
$$ L(y,n)= \sum_{k=0}^{n-1} 1_{(X_k=y)}, $$
 and set
\eq \label{greendef}
  g_B(x,y) = \mu_y^{-1} E^x L(y,\tau_B) =
 \mu_y^{-1} \sum_{k=0}^\infty P^x( X_k=y, k<\tau_B).
\en
Then $g_B(x,y)=g_B(y,x)$ and $g_B(x, \cdot)$ is harmonic
on $B\setminus\{x\}$, and zero outside $B$. Using the Markov property at
$T_y$ gives
\eq\label{pxbg}
 g_B(x,y) = P^x(T_y < \tau_B) g_B(y,y).
\en
Summing \eqref{greendef} over $y\in B$ gives
\eq \label{txgreen}
E^z\tau_B = \sum_{y\in B} g_B(z,y)\mu_y.
\en
The final property of $g_B(\cdot,\cdot)$ we will need is that
\eq\label{eq:green2}
\Reff(x,B^c)=g_B(x,x).
\en
One way to see this is to note that $g_B(x, \cdot)$ is the potential
due to a unit current flow from $x$ to $B^c$, so that $g_B(x,x)$ is
the effective resistance from $x$ to $B^c$.
Alternatively, writing
$p^{x}_B (y)= g_B(x,y)/g_B(x,x)$, one can verify that
$p^x_B$ attains the minimum in \eqref{3.3bk}, and that
$\Ecal(p^x_B,p^x_B)= g_B(x,x)^{-1}$.

\smallskip\noindent
\emph{Proof of Proposition \ref{rw-main}(a), \eqref{Etaub}.}
 It is easy to use \eqref{txgreen} to obtain an upper bound
for the exit time from a ball.
By Lemma \ref{reffacts}(d) we have $\Reff(z,B^c)\le 2R$ for any
$z\in B=B(R)$.
So,
\eq \label{tvub}
E^z\tau_B= \sum_{y\in B} g_B(z,y)\mu_y \le \sum_{y\in B} g_B(z,z) \mu_y
  = \Reff(z,B^c) V(R) \le 2\lam R^3,
\en
{which gives (\ref{Etaub})}. \qed

\noindent
\emph{Proof of Proposition \ref{rw-main}(a), \eqref{e:tmlb}}.
Write $B=B(R)$.
To obtain a lower bound for $E^0 \tau_B$ we restrict the
sum in \eqref{txgreen} to a smaller ball $B'=B(R/(4\lam))$,
and use Lemma \ref{hitest} to bound $g_B(0,y)$ from below on
$B'$. If $y \in B'$ then Lemma \ref{hitest} gives
$P^y(T_0 < \tau_B) \ge {\tfrac12}$,
so by \eqref{pxbg} and \eqref{eq:green2}
\[ g_B(0,y) = g_B(0,0)P^y(T_0 < \tau_B) \ge {\tfrac12}  g_B(0,0)
 = \tfrac12 \Reff(0,B^c) \ge \tfrac12 R/\lam.  \]
As $R/(4\lam) \in J(\lam)$ we have $\mu(B') \ge \lam^{-1} (R/(4\lam))^2$,
and therefore we obtain,
 \eq
 E^0 \tau_B \ge \sum_{y\in B'} g_B(0,y)\mu_y
 \ge {\tfrac12} g_B(0,0) \mu(B') \ge c \lam^{-4} R^3.
\en
Then {for $x\in B'$ we have}
$E^x \tau_B \ge P^x(T_0 < \tau_B) E^0 \tau_B$, which gives
\eqref{e:tmlb}.  \qed

\smallskip
The upper and lower bounds on $E^x \tau_R$ lead to a
preliminary inequality on the distribution of $\tau_R$.

\begin{lemma} \label{lem:tprelim}
Suppose that $R, R/(4\lam) \in J(\lam)$.
Let $x \in B(0,R/4\lam)$ and $n \ge 1$. Then
\eq \label{e:tplb}
   P^x( \tau_R > n) \ge \frac{C_1(\lam)R^3 -n}{2\lam R^3}  \quad
\text{ for } n \ge 0. \\
\en
\end{lemma}

\proof By the Markov property,
\eqref{Etaub} and \eqref{e:tmlb},
$$  C_1(\lam) R^3 \le
E^{x}\tau_R \le n+ E^x[1_{\{\tau_R> n\} } E^{X_n}(\tau_R)]
\le n + 2\lam R^3 P^x( \tau_R>n).$$
Rearranging this gives
\eqref{e:tplb}. \qed

\medskip
Setting $n= \delta R^3$ in \eqref{e:tplb} gives
\eq \label{tau-badub}
   P^x( \tau_R \le \delta R^3) \le 1 - \frac{ C_1(\lam) -\delta}{2\lam}.
\en
This inequality has the defect that
the right hand side of \eqref{tau-badub} does not converge to 0 as
$\delta \to 0$.
We will need a better bound in order to control $d(0,X_n)$, and
this is given in \eqref{trdlb}.

\medskip\noindent
\emph{Proof of Proposition \ref{rw-main}(a), \eqref{trdlb}}.
This proof takes a little more work;
we obtain it by a kind of bootstrap from \eqref{e:tplb}
and Lemma \ref{hitest}. The basic point is that, starting at
$y \in B(\eps R)$, $X$ is very likely to visit $0$ before escaping
from $B(R)$. So $X$ will with high probability have made many excursions from
$0$ to $\partial B(\eps R)$ before time $\tau_B$. Thus $\tau_B$ is
stochastically larger than a sum of independent
random variables, each of which, by \eqref{e:tplb}, has a probability
at least $p>0$ of being greater than $cR^3$.
Rather than following this intuition directly and using stochastic
inequalities, it is simpler
to obtain a pair of inequalities \eqref{pytub} and \eqref{pytub2} which
contain the same information.

Let $t_0>0$, and set
$$ q(y)=P^y(\tau_R {\leq} T_0), \qquad a(y)= P^y(\tau_R \le t_0). $$
Then
\begin{align}
 a(y) = P^y(\tau_R \le t_0)
 &= P^y( \tau_R \le t_0, \tau_R {\leq} T_0)
 + P^y( \tau_R \le t_0,  \tau_R > T_0)\nnb
 &\le  P^y( \tau_R \le T_0) +  P^y(T_0< \tau_R,  \tau_R -T_0 \le t_0)\nnb
\label{pytub}
 &\le q(y) + (1-q(y)) a(0)  \le q(y) + a(0),
\end{align}
using the strong Markov property for the second inequality.
Starting $X$ at 0 we have
\begin{align}\label{pytub2}
a(0) = P^0(\tau_R\le t_0) \le
E^0 [1_{\{\tau_{\eps R}\le t_0 \}} P^{X_{\tau_{\eps R}}}( \tau_R \le t_0)]
  \le  P^0( \tau_{\eps R} \le t_0)
\max_{y \in  \partial B(\eps R )} a(y).
\end{align}
Combining  \eqref{pytub} and \eqref{pytub2} gives
\eq \label{a0ub}
  a(0) \le \frac{  \max_{y   \in \partial B(\eps R )} q(y)}
 {  P^0( \tau_{\eps R} > t_0)}.
\en

Note that as $J(\lam)$ is defined to be a subset of
$[1,\infty)$, the condition that $\eps R/(4\lam) \in J(\lam)$
implies that $R \ge 4\lam/\eps$. Since $\eps < 1/(4\lam)$,
$\eps R +1 \le {2\eps R<} R/(2\lam)$, and Lemma~\ref{hitest}
(used with $2\eps$) gives
\eq
\lbeq{qy}
 q(y) \le \frac{2\eps \lam}{1-2\eps \lam} \le 4 \eps \lam.
\en
Let
$t_0= \tfrac12 C_1(\lam) (\eps R)^3$; then using \eqref{e:tplb} for the
ball $B(\eps R)$ we obtain
$$ P^0( \tau_{\eps R} > t_0) \ge \frac{C_1(\lam)}{4\lam}; $$
combining this with \refeq{qy}, (\ref{a0ub}) and (\ref{pytub}) completes
the proof of \eqref{trdlb}. \qed

\subsubsection{Heat kernel bounds}

We now turn to the heat kernel bounds in Proposition \ref{rw-main}(b).
Our first result Proposition \ref{bka}
follows from \cite[Lemmas 1.1, 1.2 and 3.10]{BCK05},
but as the proof is short we give it here.
To deal with issues related to the possible
bipartite structure of the graph it proves helpful to consider
$p_n(x,y)+p_{n+1}(x,y)$.
The main result of the proposition below is the inequality
\eqref{hkub}, which gives an upper bound for $p_n(x,x)$
just in terms of the volume.
The proof of the analogous inequality in continuous time
is a bit easier -- see \cite[Theorem 4.1]{BK06}.

\begin{prop}\label{bka}
Let $x_0\in G$ and $f_n(y)=p_n(x_0,y)+p_{n+1}(x_0,y)$. \\
(a) We have
\eq\label{speinew}
\Ecal(f_n,f_n)\le \frac {2}{n} f_{2 \lfloor n/2 \rfloor}(x_0).
\en \\
 (b) We have
\eq\label{speine2w2}
| f_n(y)-f_n(x_0) |^2
 \le  \frac {2}{n} d(x_0,y)  {f_{ 2 \lfloor n/2 \rfloor} (x_0)} .
 \en
(c) Let $r \in [1,\infty)$ and $n=2 \lfloor r \rfloor ^3$. Then
\eq \label{hkub}
f_{n}(x_0) \le c_1 n^{-2/3} ( 1 \vee (r^2/V(x_0,r)).
\en
\end{prop}

\proof (a)
It is easy to check that
$$\Ecal(f_n,f_n) = f_{2n}(x_0)- f_{2n+2}(x_0).$$
The spectral decomposition
(see for example, {Chapter~3~(32) of \cite{AF09}})
gives that $k\to f_{2k}(x_0)-f_{2k+2}(x_0)$ is non-increasing. Thus
\begin{eqnarray*}
 n \left(f_{2n}(x_0)-f_{2n+2}(x_0)\right)
 &\le & (2 \lfloor n/2 \rfloor+1)
\left(f_{4 \lfloor n/2 \rfloor}(x_0)-f_{4 \lfloor n/2 \rfloor+2}
(x_0)\right)\\
&\le &  2 \sum_{i= \lfloor n/2 \rfloor}^{2 \lfloor n/2 \rfloor}
(f_{2i}(x_0)-f_{2i+2}(x_0))\le 2f_{2 \lfloor n/2 \rfloor}(x_0),
\end{eqnarray*}
and (\ref{speinew}) is obtained.
\newline 
(b) Using Lemma \ref{reffacts}(c),(f),
\[ |f_n(y)-f_n(x)|^2 \le \Reff(x,y) \Ecal(f_n,f_n) \le
 d(x,y)  \Ecal(f_n,f_n). \]
We then use
(\ref{speinew}) to bound $\Ecal(f_n,f_n)$. \\
(c)
Choose $x_*\in B(x_0,r)$ such
that $f_{n}(x_*) = \min_{x\in B(x_0,r)} f_{n}(x)$.
Then
$$ f_{n}(x_*) V(x_0,r) \le \sum_{x\in B(x_0,r)} f_{n}(x) \mu_x
\le \sum_{x\in G} p_{n}(x_0,x) \mu_x +\sum_{x\in G}
p_{{n}+1}(x_0,x) \mu_x  \le 2, $$
so  that $f_{n}(x_*) \le 2/V(x_0,r)$.
{Since $n$ is even,} by (\ref{speine2w2}) we have
\[ f_{n}(x_0)^2 \le  2\left(f_{n}(x_*)^2 + |f_{n}(x_0)-f_{n}(x_*)|^2 \right)
\le \frac{8}{V(x_0,r)^2} + \frac{c r f_{n}(x_0)}{n}.\]
Using $a+b \leq 2(a \vee b)$, we see that
$f_{n}(x_0) \le  ({c'}/{V(x_0,r)}) \vee ({c' r}/{n})$. \qed

\begin{rem}
In fact, (\ref{speinew}) can be sharpened to give
$\Ecal(f_n,f_n)\le c_1n^{-1} p_{2 \lfloor n/2 \rfloor}(x_0,x_0),$ --
see \cite[Lemma~3.10]{BCK05}, but we do not need this.
\end{rem}

\newpage\noindent
\emph{Proof of Proposition \ref{rw-main}(b).}
Let $f_n(y)=p_n(0,y)+p_{n+1}(0,y)$.
As $R \in J(\lam)$, $R^2/V(R) \le \lam$, so by
Proposition~\ref{bka}{(c)},
\eq \label{pnub-a}
 f_n(0) \le c_1 \lam n^{-2/3}.
\en
By Proposition~\ref{bka}{(b)}, if $n$ is even
\eq
 f_{n}(y) \le f_{n}(0) + |f_{n}(y)-f_{n}(0)| \le
 f_{n}(0)+ (2 d(0,y) n^{-1} {f_n(0)} )^{1/2}
\le c \lam n^{-2/3},
\en
which proves \eqref{hkgub}.

To prove the lower bound \eqref{e:plb} we use
Lemma \ref{lem:tprelim}.
For sufficiently small $n$ this bounds from above the probability that
$X$ has left $B$ by time $n$, and so bounds from below
$P^0(X_n \in B)$. This leads
easily to a lower bound on $p_{2n}(x,x)$. Here are the details.
Let $n \le \tfrac12 C_1(\lam){R^3}$. Then using \eqref{e:tplb}
\eq\label{e:pqwes}
 P^x( X_n \in B) \ge  P^x(\tau_B>n) \ge \tfrac{1}{{4}} \lam^{-1} C_1(\lam).
\en
By Chapman--Kolmogorov and Cauchy--Schwarz
$$ P^x( X_n \in B)^2 = ( \sum_{y \in B} p_n(x,y) \mu_y)^2
\le \mu(B) \sum_{y \in B} p_n(x,y)^2 \mu_y  \le p_{2n}(x,x) \lam R^2, $$
and using (\ref{e:pqwes}) gives (\ref{e:plb}).
\qed

\bigskip
\subsubsection{ Bounds on $d(0,X_n)$}

The main work for these bounds has already been done
in the proofs of Proposition \ref{rw-main}(a) and (b),
and in particular the proof of \eqref{trdlb}.

\medskip\noindent \emph{Proof of Proposition \ref{rw-main}}.
(c)
The proof of \eqref{doxub} follows from \eqref{trdlb}
after suitable checking, since
\eq \label{dandtau}
 P^0(d(0,X_n) n^{-1/3}> M) =
 P^0(d(0,X_n)  > R)  \le P^0(\tau_{R}\le n).
\en
We now fill in the details.
Define $\eps$ by the relation $n= C_2(\lam) (\eps R)^3$; so that
$\eps= C_6(\lam)/M$. Using \eqref{dandtau} and \eqref{trdlb},
\eq
  P^0(d(0,X_n) n^{-1/3}> M) \le
P^0(\tau_{R}\le  C_2(\lam) (\eps R)^3)
\le C_3(\lam) \eps \le  \frac{C_7(\lam)}{M},
 \en
which proves \eqref{doxub}. Tracking the powers of $\lam$ gives that
$C_7(\lam) \le \lam^{22/3}$.

\smallskip \noindent
(d)  We can bound the probability that $X$ is in a ball
$B'$ by the volume of the ball and the maximum of the
heat kernel on the ball.
By \eqref{hkgub},  writing $B'=B(0,R/M)\subset B(0,R)$
and $f_n(0,y)=p_n(0,y)+p_{n+1}(0,y)$,
\eq
  P^0( X_n\in B')= \sum_{y\in B'} p_n(0,y)\mu_y \le
   \sum_{y\in B'} f_n(0,y) \mu_y \le V(R/M) C_4(\lam) R^{-2}\le \lam C_4(\lam)/M^2,
\en
proving \eqref{doxlb}.
The final inequality in (d) now follows easily, since all we need is that
$d(0,X_n)$ is greater than $cn^{1/3}$ with positive probability.
Let $M=C_8(\lam)$ satisfy
$M^2 = 2 \lam C_4(\lam)$. Then using (c),
$P^0( d(0,X_n) < R/M ) \le  \tfrac12$,
so $ E^0 d(0,X_n)\ge \tfrac12 R/M$. \qed

\bigskip

We do not have an upper bound on $E^0d(0,X_n)$ to complement the lower bound
of Proposition~\ref{rw-main}(d), which uses
volume and resistance bounds from a single base point, i.e.,
bounds on $V(0,R)$ and $\Reff(0, B(R)^c)$.
 Suppose that $J(\lam)=[1,\infty)$ for some $\lam\ge 1$,
and let $Z_n =  n^{-1/3} d(0,X_n)$.
Then we are able to bound  $E^0 Z_n^p$ for
$p< 1$, since \eqref{doxub} gives
\eqalign
 E^0 [Z_n^p] &\le  \sum_{m=1}^\infty
(2^{m+1})^p P^0\big( 2^m\le n^{-1/3}d(0,X_n)< 2^{m+1}\big)\nnb \nonumber
&\le  \sum_{m=1}^\infty
(2^{m+1})^p P^0\big(n^{-1/3}d(0,X_n)\ge 2^m\big)\le
 c_1\sum_{m=1}^\infty 2^{m(p-1)}= c_2<\infty.
\enalign
On the other hand the following example indicates that, under our hypotheses,
we cannot expect to have a uniform bound on $E^0 (Z_n^p)$ when $p>1$.
{We sketch this argument below.}

\begin{example}
\label{ex-p1}
Let $\Gam$ be the subgraph of $\Z^2$ with vertex set $G=G_0 \cup G_1$,
where $G_0= \{ (n,0), n \in \Z\}$, and
$G_1= \{ (n,m): 0\le m\le n \}$. Let the edges be
$\{ (n,0),(n+1,0)\}$, for $n \in \Z$, and
$\{(n,m),(n,m+1)\}$ if $n \ge 1$ and $0\le m \le n-1$.
Thus $\Gam$ consists of $\Z_-$ and a comb-type graph of vertical
branches with base $\Z_+$.
Write $0$ for $(0,0)$. It is easily checked that
$V(0,R) \asymp R^2$, and  $\Reff(0, B(0,R)^c) \ge R/4$. Thus
there exists $\lam_0<\infty$ such that $J(\lam_0)=[1,\infty)$.
Let
\[ H(a,b) = \{ (n,m) \in G: a\le n \le b \}. \]

Let $X_n$ be the simple random walk on $\Gam$. If we time-change
out the excursions of $X$ away from $\Z$ then we obtain a simple
random walk $Y_n$ on $\Z$. Now let $R\ge 1$, and
$r=R^{2/3} \in \Z$. Let
$A=H(-r,r)$.  Since $B(0,r/2) \subset A \subset B(0,2r)$,
Proposition~\ref{rw-main}(a) implies that
$E^0 \tau_A \approx r^3 \approx R^2$.
Since $X$ only moves horizontally when it is on the $x$-axis,
$P^0( X_{\tau_A} = (-r,0))=1/2$. If $X_{\tau_A}=(-r,0)$
then the probability that $X$ reaches $H(-\infty,-R)$ before returning to
$0$ is $r/R \approx R^{-1/3}$;
also, if $X$ does this then the time
taken to do so will be of order $R^2$.

These arguments lead us to expect that if $n = R^2$ then
\eq
 \label{e:pxineq}
   P^0\big( X_n \in H(-\infty, -R/2)\big) \ge c R^{-1/3}.
\en
Given (\ref{e:pxineq}), it follows from Markov's inequality that
\[ E^0 Z_n^p \ge n^{-p/3} (R/2)^p P^0\big(  X_n \in
  H(-\infty, -R/2)\big)\ge c n^{(p-1)/6}, \]
and the lower bound diverges if $p>1$.
 This concludes Example~\ref{ex-p1}.
\end{example}

\subsection{Results for random graphs}
\label{sec-genrwre}

We now consider a family of random graphs, as described in Section~\ref{sec-rwre},
and prove Theorems~\ref{ptight}--\ref{thm-snlim}.
{Most of the hard work has been done in the previous section, where
we obtained bounds for a fixed graph $\Gam$.  }

We begin by obtaining tightness of the quantities
$R^{-3} E^0  \tau_R$, $n^{2/3} p_{2n}(0,0)$,  and \break
$n^{-1/3}d(0,X_n)$. We recall the definition of the function $p(\lam)$
in  Assumption~\ref{ass-rwre}(1), and that $p(\lam)\le c_0 \lam^{-q_0}$.

\medskip\noindent
\emph{Proof of Theorem \ref{ptight}}.
The basic idea here is straightforward. For each of the quantities we
are interested in, the estimates in Proposition \ref{rw-main}
tell us that provided the environment is `good' at the scale $R$
(that is, more precisely, that $c_i R \in J(\lam)$ for suitable $c_i$)
then the quantity takes the value we want. The bounds we get
will only hold if $R$ or $n$ is large enough, but it is easy to handle
the small values of $R$ or $n$.

We begin with (\ref{pt-a}).
Let $\eps>0$. Choose $\lam\ge 1$ such that $2p(\lam)< \eps$.
Let $R/(4 \lam) \ge R^*$,
and set $F_1 =\{ R, R/(4\lam) \in J(\lam)\}$.
Then, by
Assumption~\ref{ass-rwre}(1), $\bP(F_1)\ge 1-2p(\lam)$.
For $\omega \in F_1$, by Proposition~\ref{rw-main}(a), there exists $c_1<\infty$,
$q_1 \ge 0$ such that
\eq
 \label{e:ttz}
 (c_1 \lam^{q_1})^{-1} \le  R^{-3} E^x_\omega \tau_R \le c_1\lam^{q_1}
 \text{ for } x \in B(R/(4 \lam)).
\en
{ So, if $\theta \ge  c_1 \lam^{q_1}$ then for $R \in [4\lam R^*, \infty)$,}
\eq
\label{e:tta}
\bP\big(  \theta^{-1} \le  R^{-3} E^0_\omega \tau_R \le \theta  \big)
\ge \bP(F_1) \ge 1 -2 p(\lam) \ge 1 - \eps.
\en
{
Let $R_0 \ge 1$. Since $0 < \sup_{1\le r\le R_0} r^{-3} E^0_\omega \tau_r <\infty$,
we have
\[ \lim_{\theta \to \infty}
\Pbold (\theta^{-1}\le r^{-3}E^0_\omega\tau_r\le \theta)=1
\quad \text{ uniformly for } r \in [1,R_0].  \]
Combining this with \eqref{e:tta} gives \eqref{pt-a}.}

A similar argument enables us to handle
the cases of small $n$ in \eqref{pt-b}--\eqref{pt-d}, {and we do not
provide further details on this point below.}

For (\ref{pt-b}) let
$n \ge 1$, $\lam \ge 1$, and let $R_0$, $R_1$ be defined by
$n= \frac12 C_1(\lam) R_1^3=2R_0^3$.
Let $F_2 =\{ R_0, R_1, R_1/(4\lam) \in J(\lam)\}$.
Suppose  that $R_0$ and $R_1/(4\lam)$
are both greater than $R^*$; then $\bP(F_2) \ge 1-3p(\lam)$.
If $\omega \in F_2$ then by Proposition~\ref{rw-main}(b)
\[ (c_2 \lam^{q_2})^{-1} \le n^{2/3} p_{2n}^\omega(0,0) \le c_2\lam^{q_2}. \]
So,
\eq
\lbeq{ttaaa}
     \bP\big(  (c_2 \lam^{q_2})^{-1} \le n^{2/3}p_{2n}^\omega
     (0,0) \le c_2 \lam^{q_2} \big)
 \ge \bP(F_2) \ge 1-3p(\lam),
\en
{proving \eqref{pt-b}.}

We now prove (\ref{pt-c}). Let $n \ge 1$ and $\lam\ge 1$.
Let $M=\lam^8$ and set
\eq \nn
 R_0= M n^{1/3}, \quad R_1= C_6(\lam)n^{1/3}, \quad
 R_2= C_6(\lam)n^{1/3}/(4\lam),
\en
$F_3 =\{ R_0, R_1, R_2 \in J(\lam)\}$.
If $n$ is large enough so that
$R_i\ge R^*$ for $0\le i \le 2$ then by
\eqref{doxub}, if $\omega \in F_3$ then
\[
    P^0_\omega\big( n^{-1/3} d(0,X_n)> \lam^8 \big)\le
    \frac{C_7(\lam)}{\lam^8}\le \frac{c\lam^{22/3}}{\lam^8}
    =\frac{c}{\lam^{2/3}}.
\]
Taking $\theta=\lam^8$, we have
\eq
\label{e:dtta}
 P^* \big( n^{-1/3}d(0,X_n)> \theta \big) \le \bP(F_3^c) +
 \bE \big(  P^0_\omega( n^{-1/3}d(0,X_n)> \lam^8 ) 1_{F_3}\big)
\le 3p(\theta^{1/8}) + c_3 \theta^{-1/12},
\en
{ and \eqref{pt-c} follows.}

{
Finally, we prove (\ref{pt-d}).
Let $R= (n/2)^{1/3}$, $M\ge 1$. If $R, R/M \in J(\lam)$ then
by \eqref{doxlb}
\eq \label{doxlb2}
 P^0_\omega\big( n^{-1/3} d(0,X_n) < 2^{-1/3}M^{-1}\big)
 \le \frac{\lam C_4(\lam)}{M^2}.
\en
Given  $\eps>0$ choose $\lam$ so that $p(\lam)<\eps$ and
$M$ so that $ \lam C_4(\lam)/M^2 < \eps$.
Let $F_4=\{ R, R/M \in J(\lam)\}$. Then
\eqref{doxlb2} holds for $\omega \in F_4$, so
taking expectations with respect to $\bP$
\begin{align*}
    P^*\big( n^{-1/3}(1+ d(0,X_n)) < 2^{-1/3}M^{-1}\big) &\le
  P^*\big( n^{-1/3} d(0,X_n) < 2^{-1/3}M^{-1}\big) \\
  &= \bE P^0_\omega\big( n^{-1/3} d(0,X_n) < 2^{-1/3}M^{-1}\big) \\
 &\le \bP( F_4^c) + \eps < 3\eps.
\end{align*}
This deals with the case of large $n$; for small $n$
we just use $1+d(0,X_n)\ge 1$. \qed
}

\bigskip\noindent
\emph{Proof of Theorem~\ref{pmeans}.}
We begin with the upper bounds in (\ref{e:mmean})--(\ref{e:pmean}).
Here all we need do is to use the bounds on $\bE V(R)$
and $\bE (1/V(R))$ given by  Assumption~\ref{ass-rwre}(2), together with
the bounds on $E^0\tau_R$ and $p_{2n}(0,0)$ obtained above.

By \eqref{tvub} and Assumption~\ref{ass-rwre}(2),
\[ \Ebold (E^0_\omega \tau_R) \le \Ebold (2 RV(R)) \le c R^3, \]
provided $R \ge R^*$. If $R \le R^*$ then since
$\tau_R \le \tau_R^*$ we obtain
the upper bound in \eqref{e:mmean} by adjusting the constant
$c_2$.
Also, by Proposition~\ref{bka}{(c)}, if $r=(n/2)^{1/3}$ then using
Assumption~\ref{ass-rwre}(3)
\[ {\Ebold p_{2n}^\omega(0,0)}
 \le c n^{-2/3} \Ebold (1 + r^2/V(r) )\le c' n^{-2/3}, \]
again provided  $r \ge R^*$.

For each of the lower bounds, it is sufficient to find a set $F\subset \Omega$
of `good' graphs with $\bP(F)\ge c>0$ such that, for all $\omega \in F$
we have suitable lower bounds on $E^0_\omega \tau_R$, $p^\omega_{2n}(0,0)$
or $E^0_\omega d(0,X_n)$.
We assume that $R\geq 1$ is large enough so that $R/(4\lam_0)\ge R^*$,
where $\lam_0$ is
chosen large enough that $p(\lambda_0)< 1/8$.
Again, we obtain the lower bound in \eqref{e:mmean} {for small $R$}
using the fact that $\bE(E^0_\omega \tau_R)\ge 1$ and adjusting
the constant $c_1$.

Let
$F =\{ R, R/(4\lam_0) \in J(\lam_0)\}$. Then $\bP(F)\ge \frac34$,
and for $\omega\in F$, by (\ref{e:tmlb}), $E^0_\omega \tau_R
\ge c_1(\lam_0) R^3$. So,
\[ \bE(E^0_\omega\tau_R)\ge \bE(E^0_\omega\tau_R 1_F)
\ge c_1(\lam_0) R^3 \bP(F) \ge c_2(\lam_0) R^3. \]
Given $n \in \N$, choose $R$ so that
$n=\frac12 C_1(\lam_0) R^3$.
{Then there exists $n^*$ (depending on $\lam_0$ and $R^*$) such that
$n \ge n^*$ implies that $R/(4\lam_0) \ge R^*$. }
Let $F$ be as above. Then using \eqref{e:plb} to bound
$p_{2n}(0,0)$ from below,
\[ \bE p_{2n}^\omega(0,0) \ge \bP(F)c_3(\lam_0) n^{-2/3} \ge c_4(\lam_0)
n^{-2/3}, \]
giving the lower bound in (\ref{e:pmean}).

A similar argument uses \refeq{Edlb} to conclude (\ref{e:dmean}). \qed

\medskip\noindent
\emph{Proof of Theorem~\ref{thm-rwre}.}
These results will follow from the bounds already obtained in Proposition~\ref{rw-main}
and in the proof of Theorem~\ref{ptight} by a straightforward Borel--Cantelli
argument.

We will take $\Omega_0 = \Omega_a \cap \Omega_b \cap \Omega_c$
where the sets $\Omega_*$ are defined in the proofs of (a), (b) and (c).
Recall that by  Assumption~\ref{ass-rwre}(1),
$p(\lam) = \bP( R \not\in J(\lam))\le c_0 \lam^{-q_0}$.

\noindent (a) We begin with the case $x=0$, and
write $w(n) = p^\omega_{2n}(0,0)$. By \refeq{ttaaa} we have
\[ \bP( {(c_1 \lam^{q_1})^{-1}} < n^{2/3} w_n \le c_1 \lam^{-q_1})
 \ge 1-3 p(\lam). \]
Let $n_k= \lfloor e^k \rfloor $ and $\lam_k= k^{2/q_0}$. Then,
since $\sum p(\lam_k) < \infty$, by Borel--Cantelli
there exists $K_0(\omega)$ with $\bP(K_0<\infty)=1$ such that
$c_1^{-1} k^{- 2 q_1/q_0} \le n_k^{2/3} w(n_k)
 \le c_1  k^{2 q_1/q_0}$
for all $k\ge K_0(\omega)$. Let $\Omega_a =\{ K_0 < \infty\}$.
For $k \ge K_0$ we therefore have
\[  c_2^{-1} (\log n_k)^{- 2 q_1/q_0} n_k^{-2/3} \le w(n_k)
 \le c_2 (\log n_ k)^{2 q_1/q_0} n_k^{-2/3}, \]
so that (\ref{e:logpnlima}) holds for the subsequence $n_k$.
The spectral decomposition
(see for example {\cite{AF09}}) 
gives that $p^\omega_{2n}(0,0)$ is monotone
decreasing in $n$.
So, if $n >N_0= e^{K_0} +1$, let $k \ge K_0$ be such that
$n_k \le n < n_{k+1}$.
 Then
\[ w(n) \le w(n_k) \le c_2 (\log n_k)^{2 q_1/q_0} n_k^{-2/3}
 \le 2e^{2/3} c_2  (\log n)^{2 q_1/q_0} n^{-2/3}. \]
Similarly $w(n) \ge w(n_{k+1}) \ge c_3 n^{-2/3} (\log n)^{-2q_1/q_0}$.
Taking $q_2 > 2q_1/q_0$, so that the constants $c_2, c_3$ can be
absorbed into the $\log n$ term, we obtain
\eq
 \label{p2nlim}
  (\log n)^{-q_2} n^{-2/3} \le  p^\omega_{2n}(0,0)
 \le  (\log n)^{q_2}  n^{-2/3} \quad \text { for all }  n\ge N_0(\omega).
\en
That $\lim_n \log p_{2n}^\omega(0,0)/\log n = -2/3$, $\bP$-a.s. is then
immediate. Since $\sum_n p^\omega_{2n}(0,0)=\infty$, $X$ is recurrent.

If $x, y \in \Ccal(\omega)$ and $k=d_\omega(x,y)$, then
the Chapman--Kolmogorov equations give that
\[ p^\omega_{2n}(x,x)(p^\omega_{k}(x,y)
\mu_x(\omega))^2 \le p^\omega_{2n+2k}(y,y), \]
{
and using this it is easy to obtain (\ref{e:logpnlima}) from
\eqref{p2nlim}. }

\medbreak
\noindent (b) Let $R_n = e^n$ and $\lam_n = n^{2/q_0}$.
Let $F_n =\{ R_n, R_n/(4 \lam_n) \in J(\lam_n) \}$.
Then (provided $R_n/(4\lam_n) \ge 1$) we have
$\bP(F_n^c) \le 2 p(\lam_n) \le 2 n^{-2}$.
So, by Borel--Cantelli, if $\Omega_b =\liminf F_n$, then
$\bP(\Omega_b)=1$. Hence there exists $M_0$ with $M_0(\omega)<\infty$
on $\Omega_b$, and such that $\omega \in F_n$ for all $n \ge M_0(\omega)$.

Now fix $\omega \in \Omega_b$, and let $x \in \Ccal(\omega)$.
Write $F(R) = E^x_\omega \tau_R$. By (\ref{e:ttz}) there exist
constants $c_4$, $q_4$ such that
\eq \label{mrnbound}
 (c_4 \lam_n^{q_4})^{-1} \le R_n^{-3} F(R_n) \le c_4 \lam_n^{q_4} .
\en
provided $n \ge M_0(\omega)$ and $n$ is also large enough so that $x \in B(R_n/(4 \lam_n))$.
Writing $M_x(\omega)$ for the smallest such $n$,
\[ c_4^{-1} (\log R_n)^{-2 q_4/q_0} R_n^{3}
  \le F(R_n) \le  c_4 (\log R_n)^{2 q_4/q_0} R_n^{3} ,
\quad \text{ for all } n \ge M_x(\omega). \]
As $F(R)$ is monotonic, the same argument as in (a) enables
us to replace $F(R_n)$ by $F(R)$, for all $R\ge R_x= 1+ e^{M_x}$.
Taking $\al_2 > 2q_4/q_0$ we obtain (\ref{e:logtaulima}).

\medbreak
\noindent (c) Recall that $Y_n =\max_{0\le k \le n} d(0,X_k)$.
We begin by noting that
\eq
 \label{ytaurel}
   \{ Y_n \geq R \} =\{ \tau_R \le n\}.
\en
Using this, (\ref{e:ynlim}) follows easily from (\ref{e:rnlim}).

It remains to prove (\ref{e:rnlim}). Since $\tau_R$ is monotone
in $R$, as in (b) it is enough to prove the result
for the subsequence $R_n=e^n$.

The estimates in (b) give the upper bound. In
fact, if $\omega \in \Omega_b$, and
$n \ge M_x(\omega)$, then
by (\ref{mrnbound})
\[ P^x_\omega( \tau_{R_n} \ge n^2 c_4 \lam_n^{q_4} R_n^3) \le
  \frac{F(R_n)}{ n^2 c_4 \lam_n^{q_4} R_n^3} \leq n^{-2}. \]
So, by Borel--Cantelli (with respect to the law $P^x_\omega$),
there exists $N'_x(\omega, \overline \omega)$ with
\[ P^x_\omega(N'_x< \infty)
 =P^x_\omega( \{ \overline \omega: N'_x(\omega, \overline \omega)<\infty\})=1 \]
such that
\[ \tau_{R_n} \le c_5 (\log R_n)^{q_5} R_n^3, \quad
\text{ for all } n \ge N'_x. \]

For the lower bound, write $C_2(\lam)=c_6 \lam^{-q_6}$,
$C_3(\lam)=c_7 \lam^{q_7}$. Let
$\lam_n= n^{2/q_0}$, and $\eps_n= n^{-2} \lam_n^{-q_6 - q_7}$.
Set $G_n=\{ R_n, \eps_n R_n, \eps_n R_n/(4 \lam_n) \in J(\lam_n)\}$.
Then, for $n$ sufficiently large so that $\eps_n R_n/(4 \lam_n) \ge 1$,
we have $\bP(G_n^c) \le 3 p(\lam_n) \le 3 c_0 n^{-2}$.
Let $\Omega_c =\Omega_b \cap (\liminf G_n)$; then by Borel--Cantelli
$\bP(\Omega_c)=1$ and there exists $M_1$ with $M_1(\omega)<\infty$
for $\omega \in \Omega_c$ such that $\omega \in G_n$ whenever
$n \ge M_1(\omega)$.
By \eqref{trdlb}, if $n\ge M_1$ and
$x \in B(\eps_n R_n)$ then
\eq
 P^x_\omega( \tau_{R_n} \le c_6 \lam_n^{-q_6} \eps_n^3 R_n^3)
 \le c_7 \lam_n^{q_7} \eps_n \le c_7 n^{-2}.
\en
So, using Borel--Cantelli, we deduce that (for some $q_8$)
\[ \tau_{R_n} \ge  c_6 \lam^{-q_6} \eps_n^3 R_n^3
 \ge n^{-q_8} R_n^3 = (\log R_n)^{-q_8} R_n^3, \]
for all $n \ge N''_x(\omega, \overline \omega)$.
This completes the proof of (\ref{e:rnlim}).
\qed

\medskip \noindent
\emph{Proof of Theorem~\ref{thm-snlim}.}
(a) We first consider the case $x=0$.
{The upper bound on $\log S_n/\log n$ follows easily from
the bounds on $\tau_R$ and $V(R)$, {as follows}.
A Borel--Cantelli argument similar to those above implies that
\eq
 V(R) \le R^2 (\log R)^{c}
\en
for all sufficiently large $R$.
Recall that $Y_n = \max_{0\le k \le n} d(0,X_n)$. We have
$W_n \subset B(Y_n)$, so $S_n \le  V(Y_n)$. So, for sufficiently
large $n$, using \eqref{e:ynlim},
\eq
 S_n \le V( (\log n)^{\alpha_3} n^{1/3}) \le n^{2/3} (\log n)^{c'},
\en
proving the upper bound in \eqref{e:snlim}.

For the lower bound, we need to show that a positive proportion
of the points in $B(Y_n)$ have been hit by time $n$, and for this
we use Lemma \ref{hitest}. }

Choose $q_1\ge 1$, $q_2\ge 1$  so that we can write
$C_2(\lam)= c_1 \lam^{-q_1}$ and $C_3(\lam) = c_2 \lam^{q_2}$.
Let $R_k=e^k$, and $\lam_k=k^{q_3}$ where $q_3\ge 2$ is chosen large
enough so that $\sum p(\lam_k) < \infty$. Let
$\eps_k = c_2^{-1} \lam_k^{-q_2} k^{-q_3}$.
Set
\[ F_k = \{ R_k, \eps_k R_k, \eps_k R_k/ 4 \lam_k \in J(\lam_k) \}. \]

Write
$\xi(x,R) = 1_{\{ T_x > \tau_R\}}$.
If $R \in J(\lam)$ and $\eps< 1/2\lam$ then by Lemma \ref{hitest},
\[ P^0_\omega( \xi(x,R) =1 ) \le \eps \lam, \quad
\text{ for } x \in B(\eps R). \]
Set
\[ Z_k = V(\eps_k R_k)^{-1} \sum_{x \in B(\eps_k R_k)}
{\xi(x,R_k)\mu_x}; \]
{
this is the proportion of points in $B(\eps_k R_k)$ which are
not hit by time $\tau_{R_k}$. }
Then if $\omega \in F_k$,
\[ P^0_\omega( Z_k \ge \tfrac12) \le 2 E^0_\omega Z_k
\le 2 \eps_k \lam_k \le k^{-q_3}. \]
Let $m(k) =  k^{q_3} \lam_k R_k^3$. Then  if $\omega \in F_k$,
by \eqref{Etaub},
\[  P^0_\omega( \tau_{R_k} \ge m(k) ) \le 2 \lam_k R_k^3 m(k)^{-1}
 = 2 k^{-q_3}. \]
Thus
\[ P^*( F_k^c \cup \{ Z_k \ge \tfrac12\} \cup\{  \tau_{R_k} \ge m(k)\} )
\le 3p(\lam_k) + 3 k^{-q_3}, \]
so by Borel--Cantelli, $P^*$-a.s.\ there exists a $k_0(\omega,\overline\omega)<\infty$
such that, for all $k\ge k_0$, $F_k$ holds, $\tau_{R_k} \le m(k)$,
and $Z_k \le 1/2$. So, for $k\ge k_0$,
\[ S_{m(k)} \ge
S_{\tau_{R_k}} {=} \sum_{x \in B(\eps_k R_k)} {(1-\xi(x,R_k))\mu_x}=
  V(\eps_k R_k) (1-Z_k) \ge \tfrac12 \lam_k^{-1} (\eps_k R_k)^2. \]
Let $n$ be large enough so that $m(k) \le n < m(k+1)$ for some
$k \ge k_0$. Then
\[ \frac{\log S_n}{\log n} \ge \frac{ \log S_{m(k)}}{ \log m(k+1)}
  \ge \frac{ 2k - c \log k}{3(k+1) + c' \log (k+1)}, \]
and the lower bound  in (\ref{e:snlim}) follows.
This proves  (\ref{e:snlim}) when $x=0$.

Now let
\[ \Omega_0 = \{ \omega: G(\omega) \text{ is recurrent and }
 P^0_\omega( \lim_n (\log S_n/\log n) = \tfrac23)=1 \}. \]
We have $\bP(\Omega_0)=1$. If $\omega \in \Omega_0$, and
$x \in G(\omega)$ then $X$ hits $0$ with $P^x_\omega$--probability
1. Since the limit does not depend on the initial segment
$X_0, \dots , X_{T_0}$, we obtain  (\ref{e:snlim}).
{

\noindent (b) We have $|W_n| \le S_n \le c_0 |W_n|$, so
(\ref{e:Wnlim}) is immediate from (\ref{e:snlim}).
\qed

\begin{rem}
Note that the constants $c_i$ in Theorem \ref{pmeans} and
$\al_i$ in Theorem \ref{thm-rwre} depend only on the constants
$c_1, c_2, c_3, q_0$ in Assumption~\ref{ass-rwre}.
\end{rem}
}

\section{Verification of Assumption~\ref{ass-rwre} for the IIC}
\label{sec:assrw}

In Section~\ref{sec-threeprop}, we state three propositions
which give estimates for the volume and effective resistance
for the \IIC.
Propositions~\ref{prop-iicvol}--\ref{prop-volbd}, which pertain
to the volume growth of $\Ccal$, are proved in Section~\ref{sec:laceexp}.
Proposition~\ref{prop:piv}, which will be used to estimate the
effective resistance, is proved in Section~\ref{sec-piv}.
In Section~\ref{sec-verass}, we use the three propositions
to verify Assumption~\ref{ass-rwre} for the \IIC, and complete
the proof of our main result Theorem~\ref{veriass}.

\subsection{Three propositions}
\label{sec-threeprop}

We will use the following notation for the \IIC.
Let $U(R)= \{ (x,n): n \ge R \}$,
$B(R) = \{(x,n) \in \Ccal: 0\le n < R\}$, and
$\partial B(R) = \{ (x,R): (x,R) \in \Ccal\}$.
We note that, using the graph distance $d$ on
$\Ccal$, $B(R)$ is just the ball $B(\zerovec,R)$, and
$\partial B(R)$ is its exterior boundary.
Let
\eq
\lbeq{ZRdef}
    Z_R = b_0 R^{-2} V(R),
\en
{where $b_0$ is a constant that will be specified
below \refeq{Z_R}.  The constant $b_0$ has limit $\frac 12$
as $L \to \infty$.}

\begin{prop}
\label{prop-iicvol}
Let $d>4$ and $L \geq L_0$.
Under the \IIC\  measure, the random variables
$Z_R$ converge in distribution to a strictly positive
limit $Z$, whose distribution is independent of $d$ and $L$.
{Also}, all moments converge, i.e.,
$\Ebold_\infty Z_R^l \to \Ebold Z^l$ for each $l \in \N$.
In particular,
\eq
\nnb
  c_1(d) R^2 \le \bE_\infty V(R) \le c_2(d) R^2 {,\,\quad R \ge 1}.
\en
Moreover, $c_1$ and $c_2$ do not depend on $d$, if
we further require that $L \ge L_1$, for some
$L_1 = L_1(d)$.
\end{prop}

\begin{rem}
We do not need the full strength of Proposition~\ref{prop-iicvol} to establish
   Assumption~\ref{ass-rwre} for the \IIC. However, since the scaling limit of $V(R)$
   is also of independent interest, we will prove the stronger result, and, moreover,
   identify the limiting random variable $Z$ in terms of super-Brownian motion.
\end{rem}

\begin{prop}
\label{prop-volbd}
Let $d>4$ and $L \geq L_0$.
\eq
\label{e:volbd}
  \Qi ( V(R) R^{-2} < \lambda )
  \le c_1(d) \exp \{ - c_2(d) \lambda^{-1/2} \},\,\quad R \ge 1.
\en
Moreover, $c_1$ and $c_2$ do not depend on $d$, if we further
require that $L \ge L_1$, for some $L_1 = L_1(d)$.
\end{prop}

The third proposition
gives an estimate on the expected number of edges at level $n-1$
that need to be cut in order to disconnect $0$ from level $R$.
We say that $(x,n), (x',n') \in \Ccal$ are \emph{RW-connected},
if there is a path, not necessarily oriented, in $\Ccal$ from
$(x,n)$ to $(x',n')$. We reserve the term \emph{connected}
to mean oriented connection, that is $(x,n) \conn (x',n')$.
Let
\eq
\lbeq{Dndef}
  D(n)
  = \left\{ e = ((w,n-1),(x,n)) \subset \Ccal : \
    \parbox{2in}{$(x,n)$ is RW-connected to level $R$ by a
    path in $\Ccal \cap U(n)$} \right\},
    \quad 0 < n \leq R.
\en
It follows from the definition that all edges in $D(n)$ need to
be cut in order to RW-disconnect $0$ from level $R$. Also, cutting
all the edges in $D(n)$ RW-disconnects $0$ from $B(R)^c$,
since for any RW-path from $0$ to $B(R)^c$ the last
crossing of level $n$ occurs at an edge in $D(n)$.

\begin{prop}
\label{prop:piv}
{Let} $d > 6$. There exists $L_1 = L_1(d) \ge L_0(d)$ such that
for $L \ge L_1$, {$R\geq 1$ and $0 < a < 1$,}
\eq
  \bE_\infty ( |D(n)| ) \le c_1(a), \qquad 0 < n \le {\lfloor aR \rfloor } .
\en
The constant $c_1(a)$ is independent of the dimension $d$
{and also of $L$}.
\end{prop}

\noindent
{\bf Remark.}
Proposition~\ref{prop:piv} is the only place where we need $d>6$ rather than
$d>4$.

\subsection{Verification of Assumption~\ref{ass-rwre} for the IIC}
\label{sec-verass}

We begin with a lemma that relates $|D(n)|$ and the effective resistance.

\begin{lemma}
For oriented percolation in any dimension $d \geq 1$,
\label{lem:rbounda}
\eq\label{e:reffdn}
  \Reff(0, \partial B(R)) \ge \sum_{n=1}^{{R}} \frac{1}{|D(n)|}.
\en
\end{lemma}

\proof
We have that $\Reff(0, \partial B(R))$ is the minimum
energy dissipation of a unit current from $\zerovec$
to $\partial B(R)$ -- see \cite[p. 63]{DS84}.
Let $I$ be such a unit current.
Fix $1\le n \le R$, let $k= |D(n)|$, and
let $J_1, \dots J_k$ be the currents in the bonds in $D(n)$.
Then since all the current must flow through the edges in $D(n)$,
we have $\sum_{i=1}^k |J_i| \ge 1$. Hence the
energy dissipation for $I$ in the bonds in $D(n)$, which is
$\sum_{i=1}^k |J_i|^2$, is greater than $1/k = |D(n)|^{-1}$.
Summing then gives \eqref{e:reffdn}. \qed

\medskip
Now we combine Proposition~\ref{prop:piv} and Lemma~\ref{lem:rbounda}
to show that it is unlikely that the effective resistance
$\Reff(0,  \partial B(R))$ is less than a small multiple of $R$.

\begin{prop}
\label{prop:rpbound}
{There is a constant $c$ such that for $d>6$, $L \geq L_1$,
$R \geq 2$ and $\epsilon >0$,}
\eq
 \label{e:refflb}
 \Qi(  \Reff(0,  \partial B(R))\le \eps R)  \le c \eps.
\en
\end{prop}

\proof
{Let $R \geq 2$.
Fix $\frac 12 < a<1$ and let $r=\lfloor aR \rfloor$;
note that $r \geq 1$.}
By Lemma~\ref{lem:rbounda} and the Cauchy--Schwarz inequality,
\eq \label{reffdn}
 \Reff(0,  \partial B(R))^{-1} \le
 \Big(\sum_{n=1}^r |D(n)|^{-1}\Big)^{-1}
\le r^{-2} \sum_{n=1}^r |D(n)|.
\en
Therefore, by Proposition~\ref{prop:piv}, Markov's inequality
and \eqref{reffdn}
\begin{align*}
 \Qi(  \Reff(0,  \partial B(R))\le \eps R)
 &= \Qi( \Reff(0,  \partial B(R))^{-1} \ge \eps^{-1} R^{-1} ) \\
 &\le \eps R \bE_\infty \big(  \Reff(0,  \partial B(R))^{-1}\big) \\
 &\le {\eps  R r^{-2}} \bE_\infty(  \sum_{n=1}^r |D(n)|) \le
 {\eps R r^{-1} c_1(a) \le 2a^{-1}c_1(a) \eps}.
\end{align*}
\qed

\noindent\emph{{Proof of Theorem \ref{veriass}.}} Let $W_R = V(R)/R^2$.  By
Proposition~\ref{prop-iicvol} we have (2) and
\eq
  \label{e:wub}
 \Qi( W_R \ge \lam)  \le \lam^{-1} {\bE_\infty W_R} \le {c}{\lam^{-1}}.
\en
Also, Proposition~\ref{prop-volbd} gives
\eq
\label{e:wlb}
 \Qi( W_R < \lam^{-1})  \le c \exp ( - c' \lam^{1/2}),
\en
and (3) is then immediate after integration.
{The combination of  (\ref{e:wub})--(\ref{e:wlb}) and (\ref{e:refflb})
(with $\eps = \lam^{-1}$),
together with the fact that each of the bounds is less than $c\lam^{-1}$
for large $\lam$, gives (1)
with $q_0=1$ and $R^*=2$.
The fact that all constants here are independent of $d,L$ implies
that the constants in Assumption~\ref{ass-rwre} share this independence.}
\qed

\section{IIC volume estimates:  Proof of Propositions~\ref{prop-iicvol}--\ref{prop-volbd}}
\label{sec:laceexp}

In Section~\ref{sec-iicvol} we prove Proposition~\ref{prop-iicvol}, and in
Section~\ref{sec-volbd} we prove Proposition~\ref{prop-volbd}.
{The proofs make use of results from several previous papers involving
the lace expansion;
these results are gathered
together and slightly extended in Section~\ref{sec-lacesum}.}

We assume throughout that $d>4$ and that $L$ is large; these
assumptions will often not be mentioned explicitly in the following.
Throughout:

\smallskip \noindent
{\emph{$\beta = L^{-d}$, $K$ denotes a constant that only depends
on $d$, and $\bar{K}$ denotes an absolute constant.}

\smallskip \noindent
The values of the constants $K$ and $\bar K$ may change from one occurrence
to the next.}

\subsection{Preliminaries}
\label{sec-lacesum}

In this section, we recall and slightly extend various results from
\cite{Hofs06,HHS02,HS02,HS03b}.  These results isolate the necessary
ingredients from other papers that will be used in the proof
of Propositions~\ref{prop-iicvol}--\ref{prop-volbd}.

\subsubsection{Critical oriented percolation $r$-point functions}

The critical oriented percolation
two-point function $\tau_n(x)$
is defined by
\eq
    \tau_n(x) = \Pbold_{p_c}( (0,0) \conn (x,n)).
\en
Let $\tau_n = \sum_{x \in \Zd}\tau_n(x)$.
By \cite[Theorem 1.1]{HS03b},
\eqalign
\label{e:taubound}
  &\sup_{x \in \Zd} \tau_n(x)
  \le K \beta (n+1)^{-d/2}, \quad n \ge 1, \\
\label{e:taun}
  &\tau_n
  = A ( 1 + \Ocal( n^{(4-d)/2} ) ), \quad
  \text{as $n \to \infty$,}
\enalign
where $|A-1|\leq K\beta$.  {The estimate
\cite[(4.2)]{HS02} shows that the error term in \eqref{e:taun}
is bounded by $K \beta n^{(4-d)/2}$ (note that
$f_n(0,z_c)$ of \cite{HS02} corresponds to our $\tau_n$).}
Hence for $L \ge L_1 = L_1(d)$, we have
\eq
\label{e:taunabs}
  \bar{K}^{-1}
  \le A
  \le \bar{K}, \qquad
  |\tau_n - A| \le \bar{K} n^{(4-d)/2},\quad n \ge 1, \qquad
  \bar{K}^{-1} \leq \tau_n \le \bar{K}, \quad n \ge 0.
\en
Also, noting that $\tau_1$ is called $p_c$ in \cite{HS03b}, we see
from \cite[Eqn.~(1.12)]{HS03b} that $|\tau_1-1|\leq K\beta \leq \bar K$
for $L \geq L_1(d)$ sufficiently large.

For all $r \geq 2$, the critical oriented percolation
$r$-point function $\tau_n^{\smallsup{r}}(x)$
is defined by
\eq
    \tau_{n_1,\ldots,n_{r-1}}^{\smallsup{r}}(x_1,\ldots, x_{r-1})
    = \Pbold_{p_c}( (0,0)\conn (x_i,n_i) \mbox{~for all~} i=1,\ldots,r-1),
\en
with $x_i \in \Zd$, $n_i \in \Z_+$.  The asymptotic behaviour of the Fourier
transforms of the $r$-point functions is given in \cite[Theorem~1.2]{HS03b}.
A very special case of \cite[Theorem~1.2]{HS03b} is that
there is a $\delta >0$ such that for $t_1,t_2>0$,
\eq
\lbeq{Vstar}
    \sum_{x_1,x_2\in \Z^d}
    \tau_{\lfloor nt_1\rfloor ,\lfloor nt_{2}\rfloor}^{\smallsup{3}}(x_1,x_{2})
    = nV^* A^3\left[ t_1 \wedge t_2 +O(n^{-\delta}) \right]
\en
as $n \to \infty$ (see \cite[(1.22)]{HS03b}).
The \emph{vertex factor} $V^*$  is
written $V$ in \cite{HS03b} but written $V^*$ here to avoid confusion with
the volume.  The vertex factor is a constant with
$|V^*-1| \leq K\beta$, and we assume that $L_1$ has been chosen so that
$\bar{K}^{-1} \le V^* \le \bar{K}$.

\subsubsection{The \IIC\ $r$-point functions}
\label{sec-IICr}

Let $\vec{y}=(y_1,\ldots,y_{r-1})$ and $\vec{m}=(m_1,\ldots,m_{r-1})$
with $y_i \in \Zd$, $m_i \in \Z_+$.  For $r \geq 2$,
the \IIC\  $r$-point function is defined by
\eq
    \iictau^\smallsup{r}_{\vec{m}}(\vec{y})
    =
    \Qi ((0,0)\conn (y_i,m_i) \mbox{~for all~} i=1,\ldots,r-1).
\en
Let
\eq
\label{e:iictau}
    \hat{\iictau}_{\vec{m}}^\smallsup{r}
    =
    \sum_{y_1,\ldots, y_{r-1}\in \Zd} \iictau^\smallsup{r}_{\vec{m}}(\vec{y}).
\en
Let $A$ be the constant
of (\ref{e:taun}), and let $V^*$ be the vertex factor of \refeq{Vstar}.
Let $r \geq 2$, $\vec{t} = (t_1,\ldots,t_{r-1}) \in (0,1]^{r-1}$,
and for a positive integer $m$,
let $m \vec t$\, be the vector with components $\lfloor m t_i \rfloor$.
It is immediate from \cite[(5.15)]{HHS02} (with $\vec k = \vec 0$) that
{for $r \geq 2$},
\eq
\lbeq{rholim}
    \lim_{m \to \infty}
    \frac{1}{(mA^2V^*)^{r-1}}
    \hat{\iictau}_{m\vec{t}}^\smallsup{r}
    =
    \hat{M}_{1,\vec{t}}^\smallsup{r},
\en
where the limit $\hat{M}_{1,\vec{t}}^\smallsup{r}$ is defined recursively as follows
(see \cite[Section~4.2]{HHS02}).

For $r=1$, we have simply
    \eq
    \lbeq{fdSBM}
    \hat{M}^{\smallsup{1}}_{s}
    =
    1.
    \en
For $r>2$ and $\bar{s}=(s_1,\ldots, s_r)$
with each $s_i > 0$,
the $\hat{M}^{\smallsup{r}}_{\bar{s}}$
are given recursively by
    \eq
    \lbeq{recmeas}
        \hat{M}_{\bar{s}}^{\smallsup{r}}
        =
        \int_0^{\underline{s}} ds \;
      \hat{M}_{s}^{\smallsup{1}}
        \sum_{I \subset J_1 : |I| \geq 1}
        \hat{M}_{\bar{s}_{I}-s}^{\smallsup{i}}
        \hat{M}_{\bar{s}_{J\backslash I}-s}^{\smallsup{r-i}}
        ,
    \en
where $i=|I|$, $J=\{1, \ldots, l\}, J_1=J\backslash \{1\}$,
$\underline{s}=\min_i s_i$, $\vec{s}_I$ denotes the vector consisting
of the components $s_i$ of $\vec{s}$ with $i \in I$, and
$\vec{s}_{I}-s$ denotes subtraction of $s$ from each component
of $\vec{s}_I$.
The explicit solution to the recursive formula
\refeq{recmeas} can be found, e.g., in \cite[(1.25)]{HS03b}.  In particular,
$\hat{M}_{s_1,s_2}^{\smallsup{2}} = s_1 \wedge s_2$.
It is shown in \cite[Lemma~4.2]{HHS02}
that for $r \geq 1$ and $t >0$,
    \eq
    \lbeq{Mform}
    \hat{M}_{t,\ldots, t}^{\smallsup{r}} = t^{r-1} 2^{-(r-1)} r!.
    \en
To this we add the following elementary fact.

\begin{lemma}
\label{lem-Ma}
For $r \geq 1$, $\hat{M}_{s_1,\ldots, s_r}^\smallsup{r}$
is nondecreasing in each $s_i$.
\end{lemma}

\proof
The proof is by induction on $r$.  For $r=1$, $\hat{M}_{s_1}^\smallsup{1}
=1$ by \refeq{fdSBM}, which is nondecreasing.  Assume the result holds
for all $j \leq r$.  Then it holds also for $r+1$ by \refeq{recmeas},
since increasing an $s_i$ can only increase the integrand (by the
induction hypothesis) or the domain
of integration in \refeq{recmeas}.
\qed

\subsubsection{Super-Brownian motion}

As discussed in \cite[Section~4]{HHS02}, the quantity $\hat{M}_{\bar s}^\smallsup{r}$
appearing in \refeq{rholim}
is the $r^{\rm th}$ moment of
the canonical measure $\N$ of super-Brownian motion $X_t$, namely
\eq
    \hat{M}_{s_1,\ldots, s_r}^\smallsup{r}
    =
    \N \big(X_{s_1}(\Rd) \cdots X_{s_r}(\Rd)\big).
\en
For an introduction to the canonical measure, see \cite[Chapter~17]{Slad06}.

Let $Y_t$ denote the canonical measure of super-Brownian motion conditioned
to survive for all time (see \cite{Hofs06}).
Let
\eq
\lbeq{Zdef}
    Z = \int_0^1 dt \; Y_t(\Rd),
\en
so that $Z$ is a positive random variable.  It is clear that the distribution
of  $Z$ does not depend on $L$.  It also does not depend on $d$, since it is
equal to the
mass up to time $1$ of the continuum random tree conditioned to survive forever.
The moments of $Z$ are given, for integers $l \geq 1$, by
\eq
\lbeq{YtM}
    \Ebold Z^l
    =
    \int_0^1 dt_1 \cdots \int_0^1 dt_l
    \hat{M}_{1,\vec{t}}^\smallsup{l+1}
\en
(see \cite[Section~3.4]{Hofs06}).
We will use the fact that $Z$ has an exponential moment.
This follows from
\eqalign
\lbeq{Zmgf}
    \Ebold Z^l
    &
    \leq
    \int_0^1 dt_1 \cdots \int_0^1 dt_l
    \hat{M}_{1,1,\dots, 1}^\smallsup{l+1}
    =
    2^{-l}(l+1)!,
\enalign
where we have used \refeq{YtM}, Lemma~\ref{lem-Ma} and \refeq{Mform}.

\subsubsection{Rate of convergence to the \IIC}

{
For the proof of Proposition~\ref{prop-volbd}, we will need an estimate
for the rate of convergence of $\Pbold_n$ to $\Pbold_\infty$
(recall the definitions
from \refeq{Pndef}--\refeq{IICdef}).}
Let $\Ecal_m$ denote the set of cylinder events measurable with
respect to the set of edges up to level $m-1$.
In \cite[Eqn.~(2.19)]{HHS02},
the following representation was obtained for $\bP_n(E)$,
$E \in \Ecal_m$:
\eq
\lbeq{Pnrep}
  \bP_n (E)
  = \frac1{\tau_n} \left[ \sum_{l=m}^{n-1}
    \vphi_l(E) \tau_1 \tau_{n-l-1} + \vphi_n(E) \right],
\en
where $\vphi_l(E)$ is a function arising in the lace expansion.
The factor $\tau_1$ was called $p_c$ in \cite{HHS02}.
By \cite[Lemma 2.2]{HHS02}, $\vphi_l$ satisfies
\eq
\label{e:phibnd1}
  |\vphi_l(E)|
  \le K \beta m (l - m + 1)^{-d/2}, \qquad l \ge m + 1.
\en
However, a very slight modification of
the proof of \cite[Lemma 2.2]{HHS02}
actually shows that
\eq
\label{e:phibnd2}
  |\vphi_l(E)|
  \le K \beta (l - m + 1)^{(2-d)/2},
  \quad l \ge m \ge 1
\en
(replace the upper bound $Km(l-m+1)^{-d/2}$ on $\sum_{a=0}^{m-1}(l-a)^{-d/2}$
used in \cite[(2.33),(2.35)]{HHS02} by the more careful upper bound
$K(l-m+1)^{(2-d)/2}$),
and we will use this variant.
The \IIC\  measure is given in \cite[Eqn.~(2.29)]{HHS02} as
\eq
\label{e:IICexpr}
  \bP (E)
  = \sum_{l=m}^\infty \tau_1 \vphi_l(E), \qquad E \in \Ecal_m.
\en

The following lemma bounds the rate at
which the measure $\bP_{2m}$ converges to $\Pbold_\infty$.

\begin{lemma}
\label{lem:IICapprox}
Let $d > 4$. For $E \in \Ecal_m$,
\eq
\label{e:IICerror}
  \left| \bP_{2m}(E)
  - \bP_\infty ( E ) \right| = \Ocal( (m+1)^{(4-d)/2} )
\en
where the constant in the error term is uniform in $E$ and
$L \ge L_0$. The error term can be guaranteed to be uniform
in $d$ as well, by further requiring that $L \ge L_1$ for
some $L_1 = L_1(d)$.
\end{lemma}

\begin{proof}
By the triangle inequality,
\eq
\lbeq{P2mtri}
    \left| \bP_{2m}(E) - \bP_\infty(E) \right|
    \leq
    \left| \bP_{2m}(E) - \sum_{l=m}^{2m} \tau_1 \vphi_l(E) \right|
    + \left| \bP_\infty(E)- \sum_{l=m}^{2m} \tau_1 \vphi_l(E) \right|.
\en
For the second term on the right-hand side, we use \refeq{IICexpr}
and \eqref{e:phibnd2} to obtain
\eq
\label{e:IICtail}
  \left| \bP_\infty (E) - \sum_{l=m}^{2m} \tau_1 \vphi_l(E) \right|
  \le \sum_{l=2m+1}^\infty \tau_1 |\vphi_l(E)|
  \le K \beta \sum_{l=2m+1}^\infty (l - m + 1)^{(2-d)/2}
  \le K \beta m^{(4-d)/2}.
\en
For the first term on the right-hand side
of \refeq{P2mtri}, we use \refeq{Pnrep}
to obtain
\eqsplit
\label{e:omittau}
  \left| \bP_{2m} (E) - \sum_{l=m}^{2m} \tau_1 \vphi_l(E) \right|
  &\le \sum_{l=m}^{2m-1} \tau_1 |\vphi_l(E)|
     \left| \frac{\tau_{2m-l-1}}{\tau_{2m}} - 1 \right|
     + |\vphi_{2m}(E)| \left| \frac{1}{\tau_{2m}} - \tau_1 \right| .
\ensplit
By \refeq{phibnd2}, the last term is bounded by
$K \beta m^{(2-d)/2}$. To bound the sum, we split it into the
cases $m \le l < 3m/2$ and $3m/2 \le l \le 2m-1$. In the
first case, we use \refeq{taun} to obtain
$|(\tau_{2m-l-1} / \tau_{2m}) - 1| \le K \beta m^{(4-d)/2}$.
Then inserting the bound \refeq{phibnd2} and summing over
$l$, we obtain a bound $K \beta m^{(4-d)/2}$ for the first
case. In the second case, we bound
$|\tau_{2m-l-1} / \tau_{2m} - 1| \le K$. Inserting the bound
on $\vphi_l$, and summing over $l$, we obtain a bound
$K \beta m^{(4-d)/2}$ for the second case.  Thus, in either case,
\refeq{omittau} is bounded by $K \beta m^{(4-d)/2}$.
For $L \ge L_1$ this bound is at most $\bar{K} m^{(4-d)/2}$.
With \refeq{P2mtri}--\refeq{IICtail}, this proves \refeq{IICerror}.
\end{proof}

\subsection{Volume convergence: Proof of Proposition~\ref{prop-iicvol}}
\label{sec-iicvol}

{In this section, we prove Proposition~\ref{prop-iicvol}.
We now choose $b_0 =(2 \tau_1 A^2 V^* R^{2})^{-1}$
in \refeq{ZRdef}, so} that $Z_R$ is defined by
\eq
\lbeq{Z_R}
    Z_R
    = (2 \tau_1 A^2 V^* R^{2})^{-1}  V(R).
\en
{As pointed out in Section~\ref{sec-lacesum},}
the constants $\tau_1, A, V^*$ all have limit $1$ as $L \to \infty$.
Let
\eq
\lbeq{tiZ_R}
    \tiZ_R
    = (A^2 V^* R^{2})^{-1} |B(R)|.
\en
Thus $\tiZ_R$ is defined in terms
of the vertices in $B(R)$, whereas
$Z_R$ is defined in terms of the edges.
Recall the random variable $Z$ defined in \refeq{Zdef}.
We use \eqref{e:rholim} to prove that $\lim_{R\to\infty}\Ebold \tiZ_R^l = \Ebold Z^l$
for all $l \geq 1$, and then adapt this to $Z_R$.

Let $l \geq 1$.  By definition,
\eqalign
\lbeq{ZRl}
    \Ebold \tiZ_R^l &= \frac{1}{(A^2V^* R^{2})^l}
    \sum_{n_1=0}^{R-1} \cdots \sum_{n_l=0}^{R-1}
    \sum_{x_1 \in \Zd} \cdots \sum_{x_l \in \Zd}
    \rho^\smallsup{l+1}_{n_1,\ldots,n_l}(x_1,\ldots,x_l)
    \nnb
    &= \frac{1}{R}\sum_{n_1=0}^{R-1} \cdots
    \frac{1}{R} \sum_{n_l=0}^{R-1}
    \frac{1}{(A^2V^* R)^l}
    \hat{\rho}^\smallsup{l+1}_{\vec{t}R} ,
\enalign
where $\vec{t} = (n_1R^{-1},\ldots, n_lR^{-1})$.
The summand on the right hand side is bounded by a constant,
by standard tree-graph inequalities \cite{AN84}
(see \cite[Section~5.1]{HHS02} for the details when $l=1$).
Therefore, by \eqref{e:rholim}, the dominated convergence theorem,
{and \refeq{YtM}},
\eq
    \lim_{R \to \infty}
    \Ebold \tiZ_R^l
    =
    \int_0^1 dt_1 \cdots \int_0^1 dt_l
    \hat{M}_{1,\vec{t}}^\smallsup{l+1} = {\Ebold Z^l}.
\en

The next lemma implies that it is also the case that $\lim_{R \to \infty}\Ebold Z_R^l
= \Ebold Z^l$ for all $l \geq 1$.

\begin{lemma}
\label{lem-ZtiZ}
For all $l \ge 1$ {and $R \geq 3$},
\eq
\label{e:ZRcompare}
  (1-2/R)^{2l}\Ebold \tiZ_{R-2}^l
  \le \Ebold Z_R^l
  \le \Ebold \tiZ_{R-1}^l  + c(d,L,l)R^{-1}.
\en
\end{lemma}

{Since $Z$ was shown in \refeq{Zmgf} to has a moment
generating function with  radius
of convergence at least $2$, the convergence of moments established
in
Lemma~\ref{lem-ZtiZ} implies that $Z_R$ converges weakly to $Z$
(see \cite[Theorem~30.2]{Bill95}).}.
Note that for
$L \ge L_1$, the constants $A$, $V^*$ and $\tau_1$ satisfy bounds
independent of $d$, hence $c_1$ and $c_2$ in
Proposition~\ref{prop-iicvol} do not depend on $d$.
This completes the proof of Proposition~\ref{prop-iicvol},
subject to Lemma~\ref{lem-ZtiZ}.

\medskip \noindent
\emph{Proof of Lemma~\ref{lem-ZtiZ}.}
For $l \geq 1$, we define
\eq
\nnb
    \sigma^\smallsup{l+1}_{\vec{m}}(\vec{x},\vec{y})
    =
    \Qi ((0,0)\conn (x_i,m_i) \conn (y_i,m_i+1)
     \mbox{~for all~} i=1,\ldots,l).
\en
Note that
\eq
\label{e:numedges}
  2 | \text{edges in $B(R-1)$} |
  \le \sum_{(x,m) \in B(R)} \mu_{(x,m)} = V(R)
  \le 2 | \text{edges in $B(R)$} |,
\en
since edges on the boundary of $B(R)$ are counted once in $V(R)$,
while other edges
are counted twice.
Therefore
\eq
\label{e:edgesR}
    \Ebold Z_R^l
    \ge \frac{1}{(\tau_1 A^2 V^* R^{2})^l}
      \sum_{n_1=0}^{R-2} \cdots \sum_{n_l=0}^{R-2}
      \sum_{x_1,y_1 \in \Zd} \cdots \sum_{x_l,y_l \in \Zd}
      \sigma^\smallsup{l+1}_{n_1,\ldots,n_l}(x_1,\ldots,x_l,y_1,\ldots,y_l),
\en
with a corresponding upper bound if the summations over the $n_i$'s
extend to $R-1$.

\smallskip \noindent \emph{Lower bound.}
The Harris--FKG inequality \cite{FKG71,Grim99} implies that
for increasing events $A$ and $B$ we have
$\Qbold_n ( A \cap B ) \ge \Qbold_n(A) \Pbold(B)$. If $A$
and $B$ are cylinder events, then by passing to the limit, we have
$\Qi(A \cap B) \ge \Qi(A) \Pbold(B)$.
Hence
\eq
\label{e:FKG}
  \sigma^\smallsup{l+1}_{\vec{n}}(\vec{x},\vec{y})
  \ge \iictau^\smallsup{l+1}_{\vec{n}}(\vec{x})
      \prod_{i = 1}^l \tau_1(y_i - x_i).
\en
With \refeq{ZRl}, this gives $\Ebold Z_R^l \ge [(R-2)/R]^{2l}\Ebold \tiZ_{R-2}^l$.

\smallskip \noindent \emph{Upper bound.}
Let
\eq
\nnb
  A_{\vec{m}}(\vec{x})
  = \{ (0,0) \conn \infty,\,
    (0,0) \conn (x_i,m_i),\, i=1,\ldots,l \}.
\en
Let $F_{\vec{m}}(\vec{x},\vec{y})$ denote the event that
the following $l+1$ events occur on disjoint
sets of edges:
\eq
\label{e:events}
  A_{\vec{m}}(\vec{x}),\,
     \{ (x_1,m_1) \conn (y_1,m_1+1) \}, \dots,
     \{ (x_l,m_l) \conn (y_l,m_l+1) \}.
\en
Then
\eq
\label{e:disjoint}
  \sigma^\smallsup{l+1}_{\vec{m}}(\vec{x},\vec{y})
  \le \Qi ( F_{\vec{m}}(\vec{x},\vec{y}) )
    + \Qi ( A_{\vec{m}}(\vec{x}){\cap_{i=1}^l\{(x_i,m_i) \conn (y_i,m_i+1)\}} \setminus
    F_{\vec{m}}(\vec{x},\vec{y}) ).
\en

{The BK inequality implies that
for increasing events $A$ and $B$ that depend on only finitely
many edges we have
$\Pbold ( A \circ B ) \le \Pbold (A) \Pbold(B)$, where $A\circ B$ denotes
disjoint occurrence \cite{BK85,Grim99}.
We will bound the first term by passing to the limit in
the BK inequality.}
{ Let
\eq
\nnb
  A_{\vec{m},n}(\vec{x})
  = \{ (0,0) \conn n,\,
    (0,0) \conn (x_i, m_i),\, i=1,\ldots, l \},
\en
and define $F_{\vec{m},n}(\vec{x}, \vec{y})$ analogously, by
replacing $A_{\vec{m}}(\vec{x})$ in \eqref{e:events} by $A_{\vec{m},n}(\vec{x})$.
Then each event in the definition of $F_{\vec{m},n}(\vec{x}, \vec{y})$ only
depends on finitely many edges, hence by BK,
\eq
\nnb
  \Pbold (F_{\vec{m},n}(\vec{x}, \vec{y}))
  \le \Pbold (A_{\vec{m},n}(\vec{x})) \prod_{i=1}^l \tau_1(y_i - x_i).
\en
Dividing both sides by $\Pbold ( (0,0) \conn n )$ and letting
$n \to \infty$, we get}
\eq
    \Qi ( F_{\vec{m}}(\vec{x},\vec{y}) ) \leq
  \Qi ( A_{\vec{m}}(\vec{x}) ) \,
    \prod_{i=1}^l \tau_1 (y_i - x_i)
  = \iictau^\smallsup{l+1}_{\vec{m}}(\vec{x})\,
    \prod_{i=1}^l \tau_1 (y_i - x_i).
\en
The sum of this bound over $\vec{x}$ and $\vec{y}$ is
$\hat{\rho}^\smallsup{l+1}_{\vec{m}} \tau_1^l$.
With \refeq{ZRl}, this gives a contribution
$\Ebold \tiZ_{R-1}^l$
to the upper bound version of
\refeq{edgesR}.

We claim that on the event
$A_{\vec{m}}(\vec{x}) {\cap_{i=1}^l\{(x_i,m_i) \conn (y_i,m_i+1)\}}
\setminus F_{\vec{m}}(\vec{x},\vec{y})$,
there exists $1 \le i \le l$ such that either
$(x_i,m_i) \conn (x_j,m_j)$ for some $j \not= i$, or
$(x_i,m_i) \conn \infty$. To see this, we may assume that
all the $(x_i,m_i)$'s are different, otherwise there is
nothing to prove. Under this assumption, the last
$l$ events in \eqref{e:events} occur disjointly. As in a
tree-graph bound \cite{AN84}, choose a set of disjoint paths
showing that $A_{\vec{m}}(\vec{x})$ occurs. Then at least
one of the paths uses an edge $((x_i,m_i),(y_i,m_i+1)$,
otherwise $F_{\vec{m}}(\vec{x},\vec{y})$ would occur.
This path includes a connection $(x_i,m_i) \conn (x_j,m_j)$
or $(x_i,m_i) \conn \infty$, proving the claim.

By the claim, the second term on the right hand side
of \eqref{e:disjoint} is at most
\eq
\label{e:onelesspath}
  \sum_{1 \le i \le l} \left[ \sum_{j \not= i}
    \Qi ( A_{\vec{m}}(\vec{x}),\, (x_i,m_i) \conn (x_j,m_j) )
    + \Qi ( A_{\vec{m}}(\vec{x}),\, (x_i,m_i) \conn \infty ) \right].
\en
Each term in \eqref{e:onelesspath} can be bounded using a
tree-graph inequality where the number of internal vertices
in the tree-graph bound is $l-1$, one less than it would be
for $\iictau^\smallsup{l+1}$. This implies that the sum of
\eqref{e:onelesspath} over $\vec{x}$ and $\vec{y}$
{inside $B(R)$} is bounded
by $c(d,L,l) R^{l-1}$. It follows that
\eq
\nnb
   \Ebold Z_R^l
   \le \Ebold \tiZ_{R-1}^l  + c(d,L,l) R^{-1},
\en
which gives the desired upper bound and completes
the proof of \eqref{e:ZRcompare}.
\qed

\subsection{Volume estimate: Proof of Proposition~\ref{prop-volbd}}
\label{sec-volbd}

In this section, we prove Proposition~\ref{prop-volbd}.
{Recall the definitions of ${\mathbb P}_n$ and ${\mathbb P}_\infty$
from \refeq{Pndef}--\refeq{IICdef}.}
It is enough to show that
we can find constants $R_0(d), c_1(d), c_2(d), c_3(d)$
such that for $R \ge R_0$ and $\lambda \le c_3$ we have
\eq
\lbeq{volbd1}
  \Pbold_\infty ( V(R) R^{-2} < \lambda )
  \le c_1 \exp \{ - c_2 \lambda^{-1/2} \}.
\en
Indeed, the restrictions on $\lambda$ and $R$ can be removed by
adjusting the constant $c_1$ as follows. First, for
$\lambda > c_3$, if $c_1 > \exp \{ c_2 (c_3)^{-1/2} \}$, the right
hand side of \eqref{e:volbd1} is larger than $1$. As for
$R < R_0$, due to the (deterministic) inequality $V(R) \ge R$,
we have $V(R) R^{-2} \ge R\,R^{-2} > R_0^{-1}$.
Therefore, if $\lambda < R_0^{-1}$, the left hand side of
\eqref{e:volbd1} is $0$. For $\lambda \ge R_0^{-1}$, it
is enough to require that $c_1 > \exp \{ c_2 R_0^{1/2} \}$.
Finally, note that if initially $R_0$, $c_1$, $c_2$, $c_3$ are
independent of $d$, then so is the adjusted $c_1$.

We begin with a simple consequence of Proposition~\ref{prop-iicvol}.

\begin{cor}
Given $\eps > 0$, there exists
$\lambda_0 = \lambda_0 (\eps,d)$, such that
\eq
\label{e:tightness}
  \Qi( V(R) R^{-2} < \lambda_0 ) < \eps, \qquad R \ge 1.
\en
For $L \ge L_1$, $\lambda_0$ can be chosen independent of $d$.
\end{cor}

\begin{proof}
This follows from Proposition~\ref{prop-iicvol}
and the fact that $Z$ is strictly positive.
\end{proof}

Let $c = c(d) = \sup_{m \ge 1} \tau_m$.
According to \eqref{e:tightness}, 
there is a constant $c_3 = c_3(d)$ such that
\eq
\label{e:simple}
  \bP_\infty ( V(R) < 4 c_3 (R+1)^2 )
  < \frac1{3 c}, \qquad R \ge 1.
\en
We fix $m_0 = m_0(d)$ such that for $m \ge m_0$
the error term on the right-hand side of \eqref{e:IICerror}
is at most $(3 c)^{-1}$.  Let $R_0 = 16 c_3 m_0^2$.
Fix $\lambda \le c_3$ and $R \ge R_0$.
We will prove that \refeq{volbd1} holds
for $\lambda$ and $R$ with the choice of $c_3$ made and with
$c_1 = 1$ and $c_2 = \frac 12\log(3/2) c_3^{1/2}$.

There is nothing to prove if $\lambda < R_0 / R^2$, since, in this case
\eq
  \bP_\infty ( V(R) R^{-2} < \lambda )
  \le \bP_\infty ( V(R) < R_0 )
  \le \bP_\infty ( V(R) < R )
  = 0
\en
and \refeq{volbd1} holds trivially.
Hence, without loss of generality, we assume that
\eq
\label{e:mu}
  \frac{16 c_3 m_0^2}{R^2}
  = \frac{R_0}{R^2}
  \le \lambda \le c_3.
\en

To estimate $\bP_\infty ( V(R) < \lambda R^2 )$, we
subdivide the time interval $[0,R]$ into blocks that
provide roughly independent contributions to the volume,
and apply \eqref{e:simple} in each block.
The number of blocks is
$S = \lfloor ( c_3 / \lambda )^{1/2} \rfloor$, which is at least $1$
by \eqref{e:mu}. The length of a block is $2m$, with
$m = \lfloor R / 2S \rfloor$.  Note that $m \ge m_0$, since
\eq
  \frac{R}{2S}
  \ge \frac{R}{2 (c_3 / \lambda)^{1/2}}
  \ge
  {\frac{R_0^{1/2}}{2c_3^{1/2}}} =2 m_0
  > 1,
\en
and hence
\eq
  m
  = \left\lfloor \frac{R}{2S} \right\rfloor
  \ge \frac{R}{4S}
  \ge \frac{R}{4 (c_3 / \lambda)^{1/2}}
  \ge m_0.
\en
Set $n_i = i(2m)$, $i = 0, \dots, S$,
so that the $i$-th block starts at level $n_{i-1}$ and
ends at level $n_i$.

By \refeq{IICdef},
\eq
\lbeq{IICVR}
  \bP_\infty ( V(R) < \lambda R^2 )
  = \lim_{N \to \infty} \frac1{\tau_N}
    \sum_{x \in \Zd} \bP_{p_c} ( V(R) < \lambda R^2,\, (0,0) \conn (x,N) ).
\en
The path $(0,0) \conn (x,N)$ on the right-hand side passes
through the levels $n_1, \dots, n_S$, and hence there exist
$0 = x_0, x_1, \dots, x_S \in \Zd$ such that
\eq
\nnb
  (0,0) \conn (x_1,n_1) \conn \cdots \conn (x_S, n_S) \conn (x,N).
\en
We write $\xvec_i = (x_i,n_i)$ for $i=0,\ldots, S$, and write
$\xvec = (x,N)$. It follows that
\eqsplit
\label{e:sumoverxs}
  &\bP_{p_c} ( V(R) < \lambda R^2,\, (0,0) \conn (x,N) ) \\
  &\qquad = \bP_{p_c} \left( \bigcup_{x_1,\dots,x_S \in \Zd}
    \{ V(R) < \lambda R^2,\, \xvec_{i-1}  \conn \xvec_i,\,
    i = 1,\dots,S \} \cap \{ \xvec_S \conn \xvec \} \right) \\
  &\qquad \le \sum_{x_1,\dots,x_S \in \Zd}
    \bP_{p_c} ( V(R) < \lambda R^2,\, \xvec_{i-1}  \conn \xvec_i,\,
    i = 1,\dots,S,\, \xvec_S \conn \xvec ).
\ensplit

Let
\eq
    \Ccal(\yvec; n) = C(\yvec) \cap (\Zd \times \{0,1,\ldots, n\}).
\en
On the event on the right-hand side
of \refeq{sumoverxs}, $\xvec_{i-1}$ is
contained in $B(R)$, and hence
$\Ccal(\xvec_{i-1};n_{i-1}+m) \subset B(R)$. Denote
$V_i = \mu(\Ccal(\xvec_{i-1};n_{i-1}+m))$.
Then on the event in the right-hand side of \eqref{e:sumoverxs},
since $\lambda \le c_3 / S^2$ by the choice of $S$, we have
\eq
  V_i
  \le V(R)
  < \lambda R^2
  \le \frac{c_3}{S^2} R^2
  = 4 c_3 \left(\frac{R}{2S}\right)^2
  \le 4 c_3 (m+1)^2.
\en
Hence, the right-hand side of \eqref{e:sumoverxs} is at most
\eq
\label{e:eachVi}
  \sum_{x_1,\dots,x_S \in \Zd}
    \bP_{p_c} \left( \bigcap_{i=1}^S
    \{ V_i < 4 c_3 (m+1)^2,\,
    \xvec_{i-1}  \conn \xvec_i \}
    \cap \{ \xvec_S \conn \xvec \} \right).
\en
The $S+1$ events in \eqref{e:eachVi} depend on disjoint sets
of bonds, so the probability factors as
\eq
\label{e:prodVi}
  \sum_{x_1,\dots,x_S \in \Zd}
  \bP_{p_c} ( \xvec_S \conn \xvec  )
  \prod_{i=1}^S
    \bP_{p_c} ( V_i < 4 c_3 (m+1)^2,\,
    \xvec_{i-1}  \conn \xvec_i )     .
\en

We insert this into \eqref{e:sumoverxs}, and use \refeq{IICVR},
\refeq{taun} and \refeq{Pndef} to obtain
\eqsplit
  \bP_{\infty} ( V(R) < \lambda R^2 )
  &\le
    \prod_{i=1}^S \left( \sum_{x_i \in \Zd}
    \bP_{p_c} ( V_i < 4 c_3 (m+1)^2,\,
    \xvec_{i-1}  \conn \xvec_i ) \right)
    \limsup_{N \to \infty} \frac{\tau_{N-n_S}}{\tau_N}
     \\
  &= \left[ \tau_{2m}
    \bP_{2m} ( V(m) < 4 c_3 (m+1)^2 ) \right]^S.
\ensplit
By Lemma \ref{lem:IICapprox}, the right-hand side equals
\eq
  \tau_{2m}^S \left[ \bP_\infty ( V(m) < 4 c_3 (m+1)^2 )
  + \Ocal ((m+1)^{(4-d)/2}) \right]^S.
\en
By the choice of $m_0$ and \eqref{e:simple}, both terms inside
the square brackets are at most $(3 c)^{-1}$.
Since
\eq
\nnb
  S
  = {\left\lfloor ( c_3/\lambda )^{1/2} \right\rfloor}
  \ge \frac 12 (c_3/\lambda )^{1/2},
\en
it follows from our choice of $c$ that
\eq
  \bP_\infty ( V(R) < \lambda R^2 )
  \le \tau_{2m}^S \left( \frac{2}{3 c} \right)^S
  \le \left( \frac{2}{3} \right)^S
  \le \exp \{ - {\textstyle \frac 12} \log(3/2)  c_3^{1/2} \lambda^{-1/2} \}.
\en
The choice
$c_2 = \frac 12\log(3/2) c_3^{1/2}$ gives \refeq{volbd1}.
Noting that for $L \ge L_1$, $c$, $c_3$ and $m_0$ (and hence all
further constants chosen) are independent of $d$, this completes the
proof of Proposition~\ref{prop-volbd}.

\section{IIC resistance estimates:  Proof of Proposition~\ref{prop:piv}}
\label{sec-piv}

In this section we prove Proposition~\ref{prop:piv}.
{Throughout, we use $\xvec, \yvec, \dots$ to denote
space-time vertices in
$\Zd \times \Z_+$, we denote the spatial component
of a vertex $\xvec$ by $x$, and we write $|\xvec|=n$ when $\xvec = (x,n)$.}
According to \refeq{Dndef},
\eq
  D(n)
  = \left\{ e = (\wvec,\xvec) \subset \Ccal : \
    \parbox{2.2in}{$|\xvec|=n$, $\xvec$ is RW-connected to level $R$
    by a path in $\Ccal \cap U(n)$} \right\},
    \quad 0 < n \leq R.
\en
Our goal is to prove that {for $d>6$, $L$ sufficiently large
and $0<a<1$,}
\eq
\lbeq{EDgoal}
  \bE_\infty ( |D(n)| )
  \le {c_1(a), \quad 0 < n \leq \lfloor aR \rfloor .}
\en
Writing $\yvec = (y,N)$, by \refeq{IICdef} and \refeq{taun} we have
\eqalign
\lbeq{EDlim}
  \bE_\infty |D(n)|
  &= \sum_{w,x \in \Zd} \bP_\infty
     \left[ (\wvec,\xvec) \in D(n) \right] \nonumber \\
  &= \frac{1}{A} \lim_{N \to \infty} \sum_{w,x,y \in \Zd}
     \bP_{p_c} \left[ (\wvec,\xvec) \in D(n),\,
     \zerovec \conn \yvec \right].
\enalign
Hence we will focus on the event
$\{ (\wvec,\xvec) \in D(n),\, \zerovec \conn \yvec \}$, for
fixed $n$, $\wvec = (w,n-1)$, $\xvec = (x,n)$ and $\yvec = (y,N)$.

\begin{figure}
\begin{center}
\setlength{\unitlength}{0.0095in}
\begin{picture}(215,200)(20,580)
\thicklines

\qbezier(110,580)(110,600)(110,620)
\qbezier(110,620)(140,670)(110,720)
\qbezier(110,620)( 80,670)(110,720)

\qbezier(120,700)(145,740)(150,780)

\thinlines
\qbezier(75,670)(110,670)(145,670)
\qbezier(75,760)(110,760)(165,760)

\put(80,680){\makebox(0,0)[lb]{\raisebox{0pt}[0pt][0pt]{$\xvec$}}}
\put(76,650){\makebox(0,0)[lb]{\raisebox{0pt}[0pt][0pt]{$\wvec$}}}
\put(60,665){\makebox(0,0)[lb]{\raisebox{0pt}[0pt][0pt]{$n$}}}
\put(160,780){\makebox(0,0)[lb]{\raisebox{0pt}[0pt][0pt]{$\yvec$}}}
\put(60,755){\makebox(0,0)[lb]{\raisebox{0pt}[0pt][0pt]{$R$}}}
\put(120,580){\makebox(0,0)[lb]{\raisebox{0pt}[0pt][0pt]{$\zerovec$}}}

\put(120,700){\makebox(0,0)[lb]{\raisebox{0pt}[0pt][0pt]{\circle*{6}}}}
\put(110,720){\makebox(0,0)[lb]{\raisebox{0pt}[0pt][0pt]{\circle*{6}}}}
\put(110,620){\makebox(0,0)[lb]{\raisebox{0pt}[0pt][0pt]{\circle*{6}}}}
\put(95,660){\makebox(0,0)[lb]{\raisebox{0pt}[0pt][0pt]{\circle*{6}}}}
\put(110,580){\makebox(0,0)[lb]{\raisebox{0pt}[0pt][0pt]{\circle*{6}}}}
\put(95,670){\makebox(0,0)[lb]{\raisebox{0pt}[0pt][0pt]{\circle*{6}}}}
\put(150,780){\makebox(0,0)[lb]{\raisebox{0pt}[0pt][0pt]{\circle*{6}}}}
\end{picture}
\end{center}

\caption{\lbfg{quick} The configuration bounded in
\refeq{quicksum}. The vertices {$\wvec = (w,n-1)$, $\xvec = (x,n)$,  $\yvec =
(y,N)$ are summed over $w,x,y \in \Zd$,} and the three unlabelled
vertices are summed over space and time.}
\end{figure}

\begin{rem}
For a quick indication of why we need to assume $d>6$, consider the configuration
in Figure~\reffg{quick}, which contributes to the right-hand
side of \refeq{EDlim}.  {Using the fact that $\tau_n$ is
bounded by a constant by \refeq{taunabs}, and using \refeq{taubound}
(see also \refeq{space-conv} below),
the configuration in Figure~\reffg{quick} can be bounded above using the BK
inequality by}
\eq
\lbeq{quicksum}
    c \sum_{l=n}^\infty \sum_{k=n}^l \sum_{j=0}^n (l-j+1)^{-d/2}
    \leq
    c \sum_{l=n}^\infty  \sum_{k=n}^l (l-n+1)^{(2-d)/2}
    \leq
    c \sum_{l=n}^\infty (l-n+1)^{(4-d)/2}
    =
    c\sum_{m=1}^\infty m^{(4-d)/2},
\en
where $j,k,l$ are the time coordinates of the unlabelled vertices, from bottom to top.
{Here, the connection from the lower unlabelled vertex to the upper unlabelled vertex
via $\wvec$ and $\xvec$ contributes $K(l-j+1)^{-d/2}$, and the other connections
all contribute constants.}
The right-hand side is bounded only for $d>6$.  Our complete proof of \refeq{EDgoal}
is more involved since we must estimate the contributions to \refeq{EDlim}
due also to more complex zigzag random walk paths.
\end{rem}

In Section~\ref{sec-il}, we prove Lemma~\ref{lem:intersect}, which
{explores the geometry of
the event $\{ (\wvec,\xvec) \in D(n),\, \zerovec \conn \yvec \}$}.
Then, in Section~\ref{sec-AJ},
we apply Lemma~\ref{lem:intersect} to construct events
$A_J(n,\wvec,\xvec,\yvec)$, $J \ge 0$,
such that
\eq
\label{e:Dnincl}
   \{ (\wvec,\xvec) \in D(n),\, \zerovec \conn \yvec \}
   \subset \bigcup_{J = 0}^\infty A_J(n,\wvec,\xvec,\yvec).
\en
In Section~\ref{sec-diag}, the BK inequality \cite{BK85} is used to obtain
a diagrammatic bound for the probability of the event $A_J(n,\wvec,\xvec,\yvec)$.
Finally, in Section~\ref{ssec:estimates}, we estimate the diagrams in this
diagrammatic bound, to prove \refeq{EDgoal} and hence Proposition~\ref{prop:piv}.
{The need to
restrict to $d>6$, rather than $d>4$, occurs only in
our last lemma, Lemma~\ref{lem:psiJs}.}

\subsection{An intersection lemma}
\label{sec-il}

We will need the existence
of certain intersections within the cluster $\Ccal$ that are implied
by the presence of a random walk path from $\xvec$ to $R$.  These
intersections are isolated in the following lemma.
The following notation will be convenient:
\eq
\nnb
  \tCcal{\pvec,\qvec}
  = \{ \vvec : \text{$\zerovec \conn \vvec$ disjointly from the edge
    $(\pvec,\qvec)$} \}, \quad (\pvec,\qvec) \subset \Ccal.
\en
Also, we write $\ol{\yvec_1 \yvec_2}$ for an occupied oriented
path $\yvec_1 \conn \yvec_2$.  Such paths are in general not unique,
but context will often identify a unique path for consideration.

{
\begin{figure}
\psfrag{x}{\fnts\ \ $\xvec$}
\psfrag{w}{\fnts\ \ $\wvec$}
\psfrag{0}{\fnts\ \ $\zerovec$}
\psfrag{p}{\fnts\ \ $\pvec$}
\psfrag{q}{\fnts\ \ $\qvec$}
\psfrag{ppr}{\fnts\ \ $\pvec'$}
\psfrag{r}{\fnts\ \ $\rvec$}
\psfrag{z}{\fnts\ \ $\zvec$}
\psfrag{A}{\fnts\ \ $A$}
\psfrag{B}{\fnts\ \ $B$}
\psfrag{R}{\fnts\ \ $R$}
\psfrag{n}{\fnts\ \ $n$}
\center{\includegraphics[scale = 0.6]{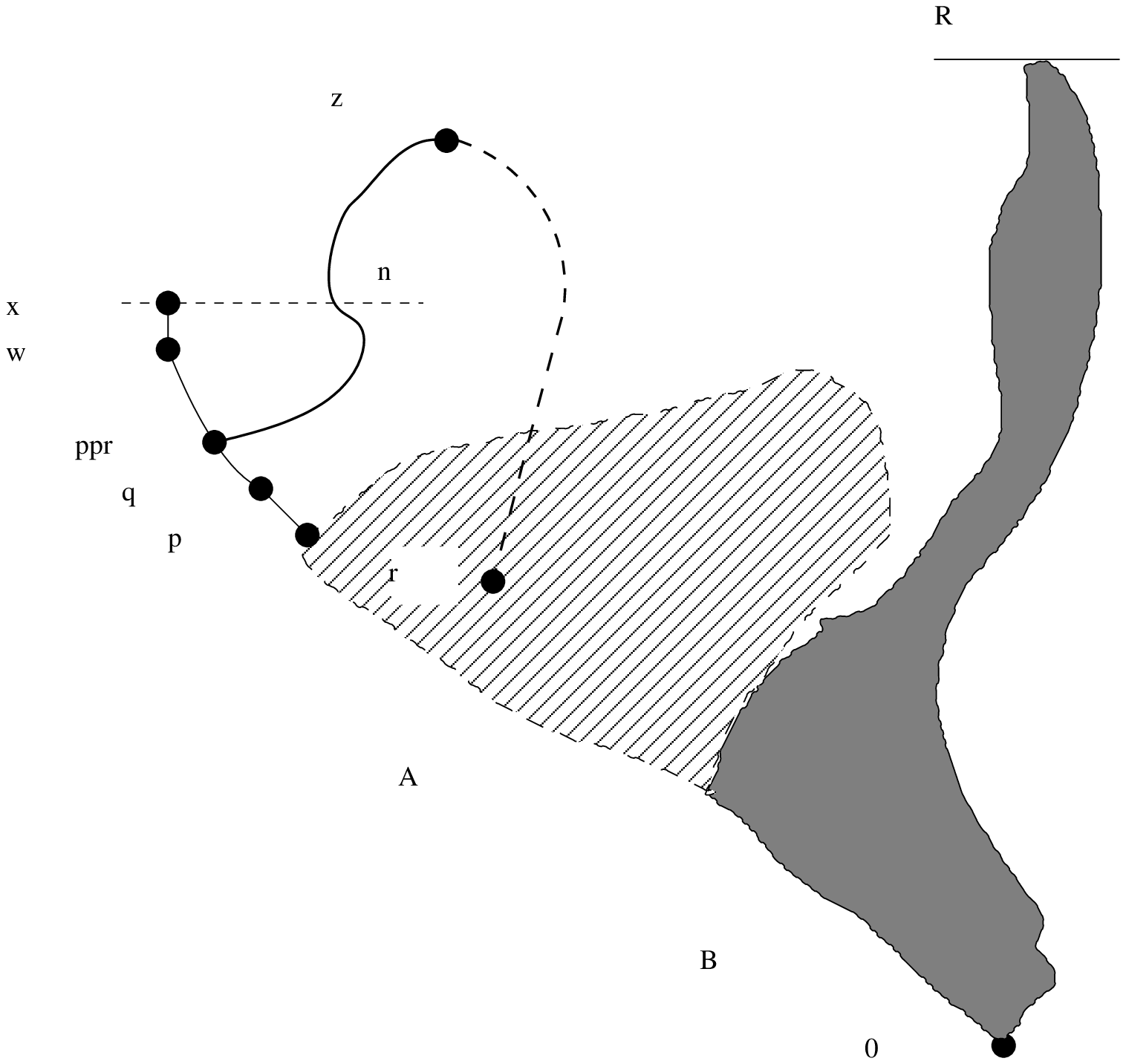}}
\caption{
\lbfg{intersection-lemma}
Illustration of the setup in Lemma~\ref{lem:intersect}.}
\end{figure}

We first describe informally the statement of the lemma,
{whose setup} is
illustrated in Figure \reffg{intersection-lemma}.
Suppose that $(\wvec, \xvec) \in D(n)$, and $\zerovec \conn \yvec$.
Let $(\pvec, \qvec)$ be an edge on an occupied path that starts
at $\zerovec$ and ends with the edge $(\wvec, \xvec)$. Assume
that $\qvec \nc R$. Then $\Ccal(\qvec)$
must intersect $\tCcal{\pvec, \qvec}$, otherwise a RW-connection
from $\xvec$ to $R$ in $\Ccal \cap U(n)$ could not occur.
Indeed, $\Ccal(\xvec)$
would have to intersect $\tCcal{\pvec, \qvec}$, but
the lemma gives a
more sophisticated version of the intersection requirements, which
allows us to have some control over the way the intersection occurs.
This is needed, because we will use the lemma recursively to
construct a set of paths realizing the intersections.
Assume that we are given a subgraph $A \cup B$ of
$\tCcal{\pvec, \qvec}$, that will represent a set of paths
already constructed, where $A$ will be a certain `preferred region.'
Assume that $A \cup B$ is disjoint from $\Ccal(\qvec)$, and
$\zerovec \in A \cup B$. Then there will be upwards occupied
paths from some vertex $\rvec \in A \cup B$ and some vertex
$\pvec' \in \ol{\qvec \xvec}$ to an intersection point
$\zvec$. It will be convenient, if we can also conclude that
$\rvec$ is in the preferred region $A$. For this reason, we will
also assume that any occupied path from $B$ to $\Ccal(\qvec)$
passes through $A$. Now we state the lemma precisely.}

\begin{lemma}
\label{lem:intersect}
Assume the event
$\{ (\wvec,\xvec) \in D(n),\, \zerovec \conn \yvec \}$. In addition,
assume the following:
\begin{itemize}
\item[(i)] $(\pvec,\qvec) \subset \Ccal$ and either
  $\qvec \conn \wvec$ or $(\pvec, \qvec) = (\wvec, \xvec)$;
\item[(ii)] $\qvec \nc R$;
\item[(iii)] $A$ and $B$ are subgraphs of $\tCcal{\pvec,\qvec}$
  with $\zerovec \in A \cup B$, and such that
  $(A \cup B) \cap \Ccal(\qvec) = \es$;
\item[(iv)] every occupied oriented path from $B$ to $\Ccal(\qvec)$
  passes through a vertex of $A$.
\end{itemize}
Then there exist $\pvec' \in \ol{\qvec \xvec}$,
$\rvec \in A$ and $\zvec$ with $|\pvec| < |\zvec| < R$, such that
\eq
\nnb
 \pvec' \conn \zvec \; \text{and} \; \rvec \conn \zvec \;
   \text{edge-disjointly, and edge-disjointly from} \;
   \ol{\pvec \xvec} \cup A \cup B.
\en
Here $\zvec$ may coincide with $\pvec'$ or $\rvec$.
\end{lemma}

\begin{proof}
We first show that $\Ccal(\qvec)$ and $\tCcal{\pvec,\qvec}$ must have
a common vertex $\vvec$. Fix a random walk path $\Gamma$ from $\xvec$ to $R$
in $U(n)$, showing that $(\wvec,\xvec) \in D(n)$. Note that $\Ccal$
(as a set of vertices) is the union
$\tCcal{\pvec,\qvec} \cup \Ccal(\qvec)$. Since
$\Gamma$ starts at $\xvec \in \Ccal(\qvec)$, but $\qvec \nc R$,
there is an edge $(\vvec,\vvec') \subset \Gamma$ such that
$\vvec \in \Ccal(\qvec)$ but $\vvec' \not\in \Ccal(\qvec)$,
and therefore $\vvec' \in \tCcal{\pvec,\qvec}$.
We need to have $|\vvec'| = |\vvec| - 1$ (otherwise
$\vvec' \in \Ccal(\qvec)$). We can rule out
$(\vvec',\vvec) = (\pvec,\qvec)$, since $\Gamma$
stays in $U(n)$, and $|\pvec| \le n-1$. It follows that
$\vvec \in \tCcal{\pvec,\qvec}$, and hence is in the intersection
$\Ccal(\qvec) \cap \tCcal{\pvec,\qvec}$.

Choose $\zvec \in \Ccal(\qvec) \cap \tCcal{\pvec,\qvec}$ with
$|\zvec|$ minimal. Since $\qvec \nc R$,
$|\pvec| < |\qvec| \le |\zvec| < R$.

We can find occupied oriented paths
$\ol{\qvec \zvec} \subset \Ccal(\qvec)$ and
$\ol{\zerovec \zvec} \subset \tCcal{\pvec,\qvec}$.
These two paths must be edge-disjoint by minimality
of $|\zvec|$. Let $\pvec'$ be the last visit of
$\ol{\qvec \zvec}$ to $\ol{\qvec \xvec}$, and let
$\rvec$ be the last visit of $\ol{\zerovec \zvec}$ to
$A \cup B$. Such a last visit exists, since we assumed
$\zerovec \in A \cup B$. Since $\zvec \not\in A \cup B$,
due to $(A \cup B) \cap \Ccal(\qvec) = \es$,
the last visit has to be in $A$ by assumption (iv).

The path $\ol{\pvec' \zvec}$ is edge-disjoint from
$\ol{\pvec \xvec}$, by the definition of $\pvec'$. It is
also edge-disjoint from $A \cup B$, by minimality of $|\zvec|$.
Likewise, the path $\ol{\rvec \zvec}$ is edge-disjoint from
$A \cup B$ by definition of $\rvec$. It is also edge-disjoint
from $\ol{\pvec' \zvec}$, by minimality of $|\zvec|$.
\end{proof}

\medbreak

\begin{rem}
Note that in the proof, we have first found a
vertex $\rvec \in A \cup B$, and assumption (iv) was only
used to show that we must have $\rvec \in A$. In fact, without
assumption (iv), we would get the statement of the Lemma
with $\rvec \in A \cup B$. The significance of being able to
ensure that $\rvec$ is in the smaller set $A$, as well as the
roles played by $A$ and $B$ will become apparent in Section~\ref{sec-AJ}.
\end{rem}

\subsection{The event $A_J(n,\wvec,\xvec,\yvec)$}
\label{sec-AJ}

In this section, we define the event $A_J(n,\wvec,\xvec,\yvec)$
and prove \refeq{Dnincl}.  The following lemma is key.

\begin{lemma}
\label{lem:paths}
Let $e = (\wvec,\xvec)$, and assume the event
$\{ e \in D(n),\, \zerovec \conn \yvec \}$. Then there exists $J \ge 0$,
such that the following vertices and paths (all edge-disjoint) exist:
\begin{itemize}
\item[(i)] vertices
  {$\uvec_0, \uvec_1, \dots, \uvec_J = \wvec$ such that
    $0 \le |\uvec_0| \le |\uvec_1| \le \dots \le |\uvec_J| = n - 1$;}
\item[(ii)] vertices
  {$\vvec_0, \vvec_1, \dots, \vvec_J = \xvec$}, and, if $J \ge 1$, vertices
  $\zvec_1, \dots, \zvec_J$ such that
\begin{gather}
  |\uvec_{i-1}|
  \le |\vvec_{i-1}|
  \le |\zvec_i|, \quad 1 \le i \le J; \label{e:uvz_i} \\
  |\uvec_{i-1}|
  < |\zvec_i|
  < R, \quad 1 \le i \le J; \label{e:vzR_i}
\end{gather}
\item[(iii)]
  $\zerovec \conn \uvec_0$ and $\uvec_{i-1} \conn \uvec_{i}$, $1 \le i \le J$;
\item[(iv)] $\uvec_{i-1} \conn \zvec_i$, $1 \le i \le J$;
\item[(v)]  $\vvec_{i-1}$ lies either on $\ol{\uvec_{i-1} \uvec_i}$
  or $\ol{\uvec_{i-1} \zvec_i}$, \, and
  $\vvec_i \conn \zvec_i$, \, $1 \le i \le J$.
\end{itemize}
In addition, at least one of the following holds: Case~(a)
$\vvec_0 \conn \yvec$; Case~(b) $\vvec_0 \conn R$ and there exists
$\vvec_*$ on $\ol{\zerovec \uvec_0}$ such that
$\vvec_* \conn \yvec$.
\end{lemma}

\begin{figure}
\psfrag{v0l0}{\ $\vvec_3$}
\psfrag{u0k0}{\ $\uvec_3$}
\psfrag{u1k1}{\ $\uvec_2$}
\psfrag{u2k2}{\ $\uvec_1$}
\psfrag{u3k3}{\ $\uvec_0$}
\psfrag{00}{\ $\zerovec$}
\psfrag{z1s1}{\ $\zvec_3$}
\psfrag{v1l1}{\ $\vvec_2$}
\psfrag{v2l2}{\ $\vvec_1$}
\psfrag{v3l3}{\ $\vvec_0$}
\psfrag{z2s2}{\ $\zvec_2$}
\psfrag{z3s3}{\ $\zvec_1$}
\psfrag{yN}{\ $\yvec$}
\psfrag{R}{\ $R$}
\psfrag{v4l4}{\ $\vvec_*$}
\psfrag{n}{\ $n$}
\center{\includegraphics[scale = 0.6]{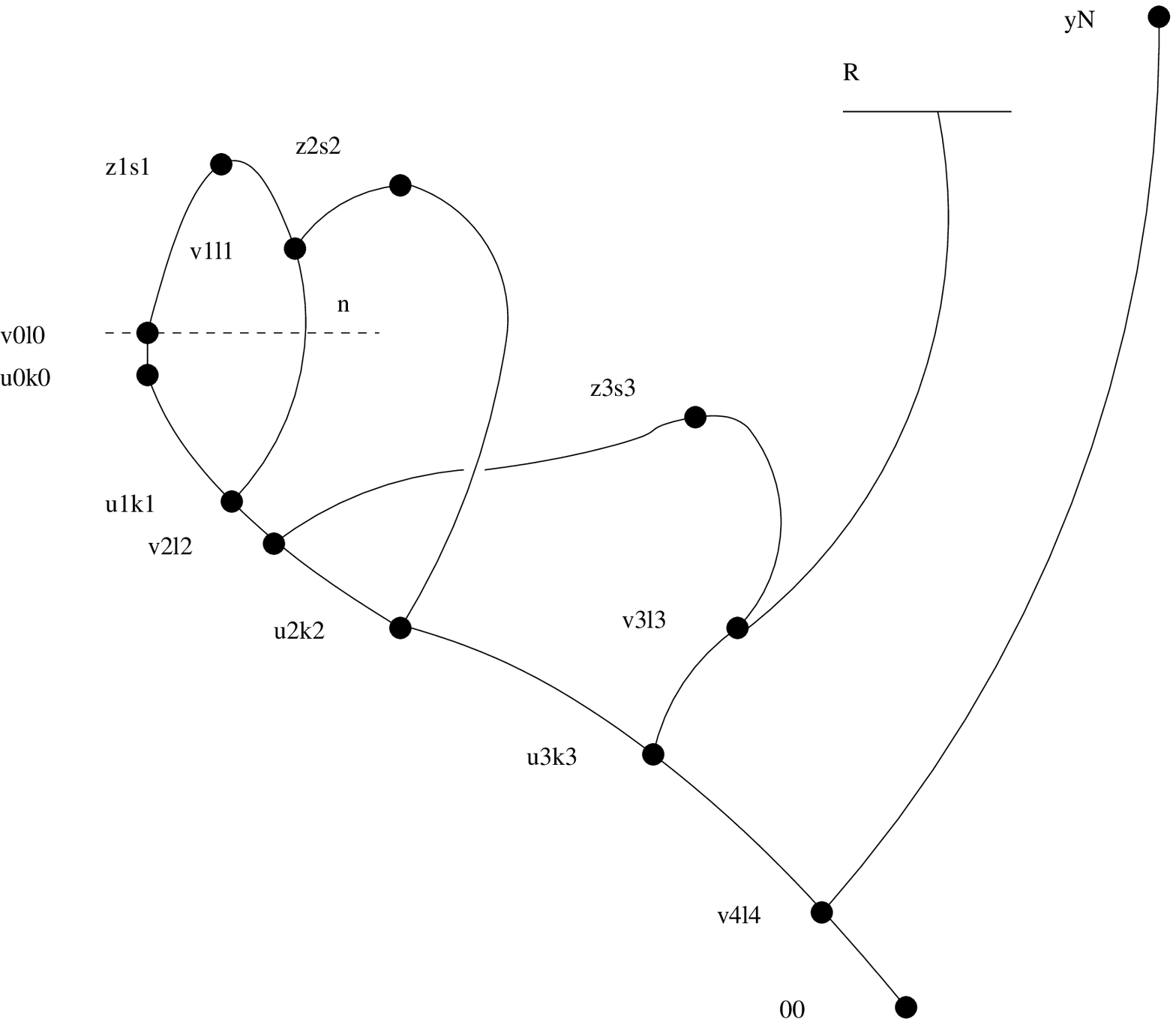}}
\caption{
\lbfg{paths}
The vertices and disjoint paths of $A_J(n,\wvec,\xvec,\yvec)$
for $J = 3$.  Here $\xvec = \vvec_3$ and $\wvec =\uvec_3$.}
\end{figure}

\begin{defn}
\label{def-AJ}
We denote by $A_J = A_J(n,\wvec,\xvec,\yvec)$ the event that the vertices and
disjoint paths listed in Lemma~\ref{lem:paths} exist,  and
$(\wvec,\xvec)$ is occupied.  See Figure~\reffg{paths}.
\end{defn}

The inclusion \refeq{Dnincl} then follows immediately from
Lemma~\ref{lem:paths}.

\medskip\noindent
\emph{Proof of Lemma~\ref{lem:paths}.}
Throughout the proof, we assume the event
$\{ e = (\wvec, \xvec)\in D(n),\, \zerovec \conn \yvec \}$.

We first show that if $\xvec \conn R$ then the lemma holds
with $J = 0$. Indeed, take $\uvec_0 = \wvec$ and $\vvec_0 = \xvec$.
Then $\zerovec \conn \uvec_0$, since $\uvec_0 \in \Ccal$.
Hence it is left to show that at least one of Cases (a) and (b)
holds. If $\vvec_0 = \xvec \conn \yvec$, then Case (a) holds.
If not, then since $\zerovec \conn \yvec$ we can find
$\vvec_* \in \ol{\zerovec \uvec_0}$ such that $\vvec_* \conn \yvec$
edge-disjointly from $\ol{\zerovec \uvec_0}$. The connection
$\ol{\vvec_* \yvec}$ has to be edge-disjoint from
$\ol{\wvec \xvec R}$, otherwise we are in Case (a). Hence
Case (b) holds.

For the rest of the proof, we assume $\xvec \nc R$.

\begin{figure}
\psfrag{x}{\fnts\ \ $\xvec$}
\psfrag{w}{\fnts\ \ $\wvec$}
\psfrag{p2}{\fnts\ \ $\pvec_2$}
\psfrag{q2}{\fnts\ \ $\qvec_2$}
\psfrag{u1r1}{\fnts\hskip-.5cm\ $\uvec_1 = \rvec_1$}
\psfrag{u0}{\fnts\ \ $\uvec_0$}
\psfrag{u0p0}{\fnts\hskip-.5cm\ $\uvec_0 = \pvec_0$}
\psfrag{v0r0}{\fnts\hskip-.5cm\ $\vvec_0 = \rvec_0$}
\psfrag{0}{\fnts\ \ $\zerovec$}
\psfrag{p1}{\fnts\ \ $\pvec_1$}
\psfrag{q1}{\fnts\ \ $\qvec_1$}
\psfrag{v1p1}{\fnts\hskip-.5cm\ $\vvec_1 = \pvec_1$}
\psfrag{v0}{\fnts\ \ $\vvec_0$}
\psfrag{z2}{\fnts\ \ $\zvec_2$}
\psfrag{z1}{\fnts\ \ $\zvec_1$}
\psfrag{R}{\fnts\ \ $R$}
\psfrag{n}{\fnts\ \ $n$}
\center{(a) \includegraphics[scale = 0.5]{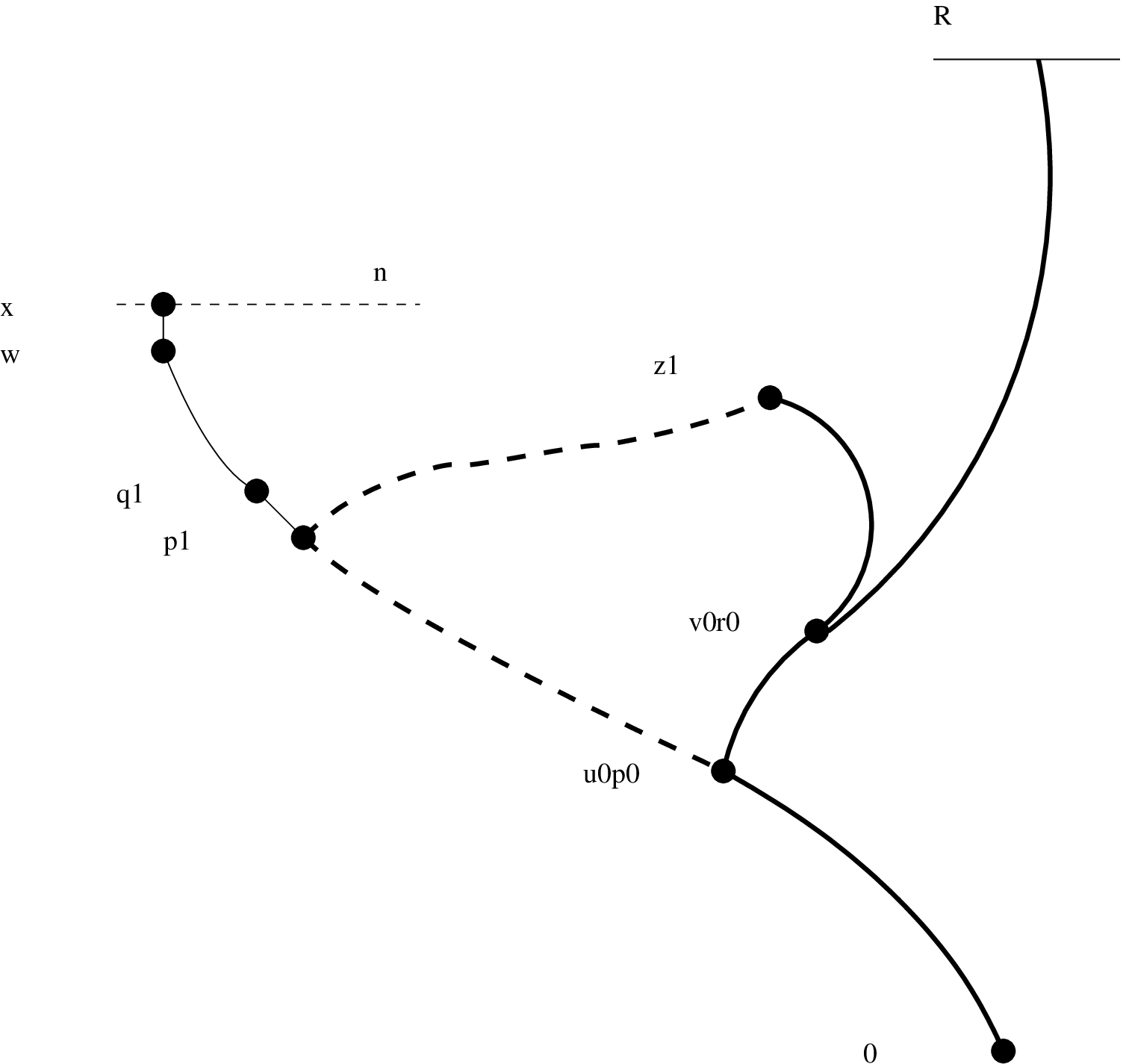} \hfill
(b) \includegraphics[scale = 0.5]{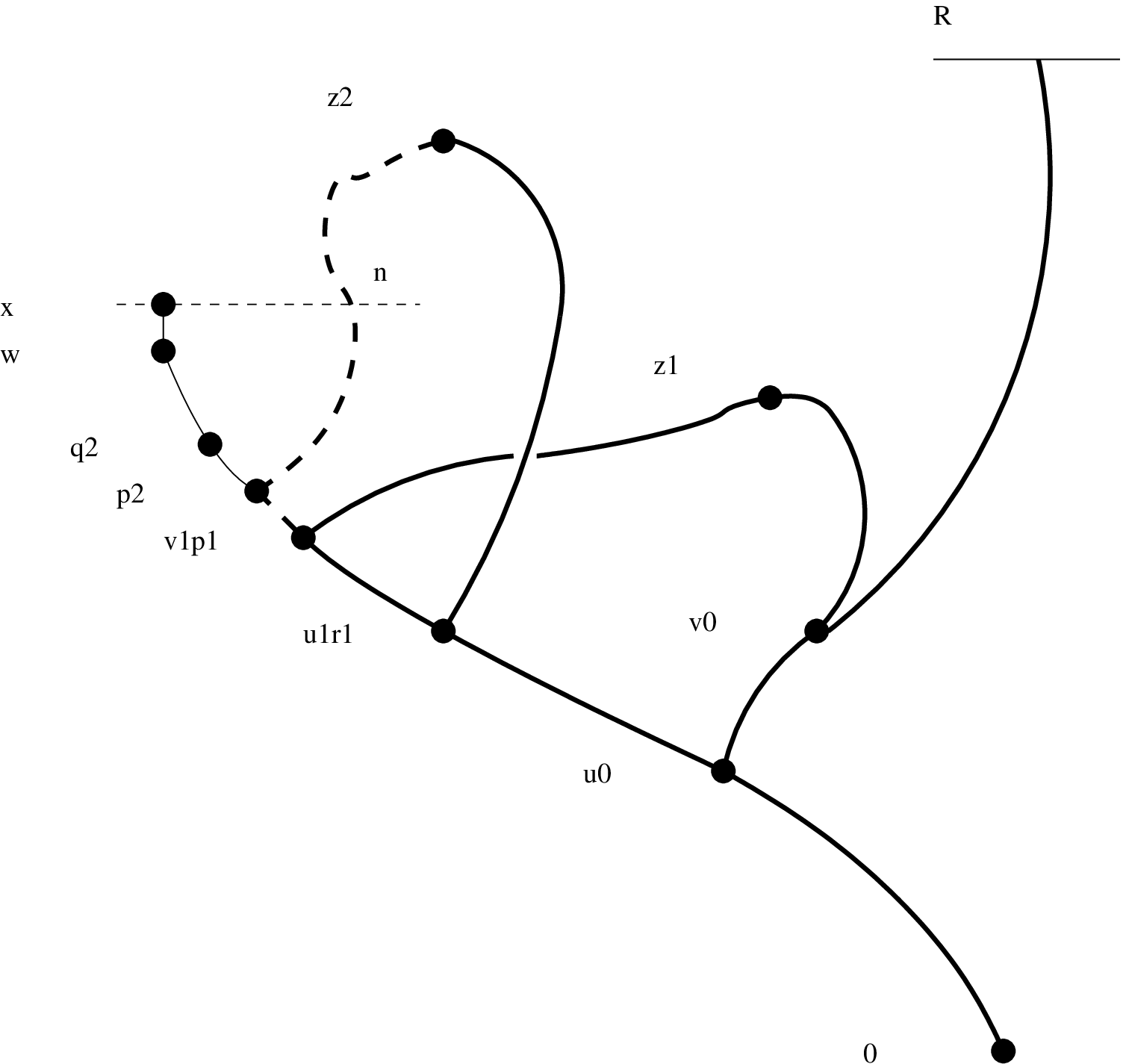}}
\caption{
\lbfg{recursion}
Assumptions of the recursion hypothesis for (a) $I = 1$; (b) $I = 2$.
The thick solid lines indicate the sets (a) $B_1$ and (b) $B_2$, and the
thick dashed lines the sets (a) $A_1$ and (b) $A_2$. The intersection lemma
is used to produce paths that join the thick dashed lines to the thin
solid lines.}
\end{figure}

We construct the paths claimed in the lemma recursively. Hence our
proof will be based on a recursion hypothesis whose statement involves
an integer $I \ge 0$, and which says that a subset of the paths claimed
in the lemma (depending on $I$) have already been constructed. 
In order to advance the recursion, the hypothesis also specifies graphs
$A_I$ and $B_I$ such that Lemma \ref{lem:intersect} can be applied with
$A = A_I$ and $B = B_I$.

The outline of the proof is the following.
Since the statement of the hypothesis for $I = 0$ is slightly
different than for $I \ge 1$, we state and verify the
hypothesis for $I = 0$ separately. This will show that the recursion
can be started. Since the general step of the recursion is
complex, we explain the first two steps of the recursion
($I = 1$ and $I = 2$) in some detail, before formulating the recursion
hypothesis precisely in the general case $I \ge 1$.
The recursion will lead to the proof of the lemma by the following
steps. 
We prove that if the hypothesis
holds for some value of $I \geq 0$, then either the conclusion of
Lemma~\ref{lem:paths} follows with $J = I + 1$, or else the hypothesis
also holds for $I + 1$. If, for some $i>0$,
the hypothesis holds for $I=0,1,\ldots, i$, then its statement will guarantee
the existence of vertices $\pvec_0, \pvec_1, \dots, \pvec_i$
with
\eq
\lbeq{pmon}
    |\pvec_0| <  |\pvec_1| < \cdots < | \pvec_i| < n.
\en
Consequently the hypothesis cannot hold for all $I=0,1,\ldots,n$,
and the implications just mentioned provide a proof of Lemma~\ref{lem:paths}.
We now carry out the details.

\medbreak

\noindent
\emph{{\bf (R) Recursion hypothesis for $I = 0$.}
There exists $\pvec_0, \qvec_0$ such that
\eqalign
\label{e:A0}
   \zerovec \conn \pvec_0, \qquad
   \pvec_0 \conn R,  \\
\label{e:rest0}
   \pvec_0 \conn \wvec \conn \xvec, \qquad
   \qvec_0 \nc R,
\enalign
where $(\pvec_0,\qvec_0)$ is the first edge in the path
$\ol{\pvec_0 \xvec}$. All paths stated are edge-disjoint.
Letting
\eqalign
  A_0
  &= \{ \ol{\zerovec \pvec_0}, \ol{\pvec_0 R} \}
  = \{ \text{paths in \eqref{e:A0}} \}, \nonumber \\
  B_0
  &= \es, \nonumber
\enalign
the hypotheses of Lemma~\ref{lem:intersect} are satisfied
with $\pvec = \pvec_0$, $\qvec = \qvec_0$, $A = A_0$ and
$B = B_0$.}

\medbreak

\noindent
{\bf Verification of (R) for $I = 0$.}
Since $\zerovec \conn \wvec$ and $\zerovec \conn R$, there exists
$\pvec_0$ such that
\eq
\nnb
   \zerovec \conn \pvec_0, \qquad
   \pvec_0 \conn \wvec \quad \text{and} \quad
   \pvec_0 \conn R \quad \text{disjointly.}
\en
Fix the paths $\ol{\zerovec \pvec_0}$, $\ol{\pvec_0 \wvec}$
and $\ol{\pvec_0 R}$, and let $(\pvec_0,\qvec_0)$ be the first
step of the path $\ol{\pvec_0 \xvec}$. If we select $\pvec_0$
so that $|\pvec_0|$ is maximal, then we have $\qvec_0 \nc R$.
We verify the hypotheses of Lemma~\ref{lem:intersect} with
these choices. First, (i), (ii) and $\zerovec \in A_0 \cup B_0$
are immediate. Also, $\Ccal(\qvec_0) \cap (A_0 \cup B_0)
= \Ccal(\qvec_0) \cap A_0 = \es$, since otherwise
$\qvec_0 \conn R$. Finally, (iv) is vacuous, since $B_0$ is empty.

\medbreak

Next, to illustrate the main idea of the proof, we explain the first
two steps of the recursion.

Since we have verified {\bf (R)} in the case $I = 0$, we can
apply Lemma~\ref{lem:intersect} with $\pvec = \pvec_0$,
$\qvec = \qvec_0$, $A = A_0$ and $B = B_0$. Lemma~\ref{lem:intersect}
shows that there exist $\pvec' \in \ol{\qvec_0 \xvec}$ and
$\rvec \in A_0 = \ol{\zerovec \pvec_0} \cup \ol{\pvec_0 R}$
and a vertex $\zvec$ such that $\pvec' \conn \zvec$ and
$\rvec \conn \zvec$. For reasons that will be explained
{in the third paragraph below},
we select $\pvec'$ with $|\pvec'|$ maximal such that the
conclusions of Lemma~\ref{lem:intersect} hold.
With this choice of $\pvec'$, we set $\pvec_1 = \pvec'$,
$\zvec_1 = \zvec$ and $\rvec_0 = \rvec$.
Note that $|\pvec_1|>|\pvec_0|$.
We define the
vertices $\uvec_0$ and $\vvec_0$ as follows.
Note that $\rvec_0 \in A_0$, which is the union of the
paths $\ol{\zerovec \pvec_0}$ and $\ol{\pvec_0 R}$.
If $\rvec_0 \in \ol{\pvec_0 R}$
then we set $\vvec_0 = \rvec_0$ and
$\uvec_0 = \pvec_0$,
and if $\rvec_0 \in \ol{\zerovec \pvec_0}$ then
we set $\vvec_0 = \pvec_0$, $\uvec_0 = \rvec_0$.
In either case, we have
$|\uvec_0| \le |\pvec_0| < |\zvec_1| < R$,
and hence \eqref{e:vzR_i} holds for $i = 1$.

The paths constructed so far are depicted in
Figure \reffg{recursion} (a). For the moment, the reader
should disregard $\qvec_1$, and the distinction between
thin, thick and dashed paths in the figure.
We either have $|\pvec_1| < |\xvec| = n$, as depicted in
Figure~\reffg{recursion}(a), or $\pvec_1 = \xvec$.

We first argue that in the case $\pvec_1 = \xvec$,
Lemma \ref{lem:paths} holds with $J = 1$. Indeed,
if $\pvec_1 = \xvec$, we set $\uvec_1 = \wvec$ and
$\vvec_1 = \xvec$. Then apart from the claim regarding
Cases (a) and (b), the vertices and paths required by
Lemma \ref{lem:paths} for $J = 1$ have been constructed.
(Note that the conclusion of Lemma \ref{lem:intersect}
guarantees that the newly constructed paths are
edge-disjoint from the old ones.) It is not difficult
to also show that either Case (a) or (b) holds, and
we leave the details of this to when we deal
with the general recursion step.

Next we explain how to continue the construction if
$|\pvec_1| < |\xvec| = n$. Let $\qvec_1$ denote the
first vertex on the path $\ol{\pvec_1 \xvec}$ following
$\pvec_1$. Let $B_1$ denote the union of the thick
solid lines in Figure~\reffg{recursion}(a), that is,
$B_1 = \ol{\zerovec \pvec_0} \cup \ol{\pvec_0 R} \cup
\ol{\rvec_0 \zvec_1} = A_0 \cup \ol{\rvec_0 \zvec_1}$.
Let $A_1$ denote the union of the dashed lines
in Figure~\reffg{recursion}(a), that is,
$A_1 = \ol{\pvec_0 \pvec_1} \cup \ol{\pvec_1 \zvec_1}$.
We want to apply Lemma~\ref{lem:intersect} with
$A = A_1$, $B = B_1$, etc. It is easy to verify
conditions (i)--(iii) of the lemma. The crucial
condition here is (iv), which allows us to conclude that
$\rvec \in A_1$, and hence the two new paths produced
by Lemma~\ref{lem:intersect} will connect the
dashed lines to the thin solid lines in
Figure~\reffg{recursion}(a). The reason condition
(iv) is satisfied is that we chose $|\pvec_1|$
to be maximal. Indeed, a glance at Figure~\reffg{recursion}(a)
suggests that if we had paths from $\ol{\qvec_1 \xvec}$ and
$B_1 \setminus A_1$ to a vertex $\zvec$ that are edge-disjoint
from $A_1 \cup B_1$, then that would contradict the
maximality of $|\pvec_1|$. (Recall the earlier
application of Lemma~\ref{lem:intersect} with $A = A_0$,
$B = B_0$, etc., and the choice of $\pvec_1$.)
We will verify the details of this when we deal with the
general case $I \ge 1$.

We can summarize the above discussion by saying that
Hypothesis {\bf (R)} for $I = 0$ should imply that
in the case $\pvec_1 \not= \xvec$ the following
statement holds.

\medbreak

\noindent
\emph{{\bf (R) Recursion hypothesis for $I = 1$.}
Vertices and paths (all edge-disjoint) with the following
properties exist:
\begin{itemize}
\item[(i)] $\pvec_1$ and $\qvec_1$ such that
\eq
\label{e:restIzzz}
  \pvec_1 \conn \wvec \conn \xvec, \qquad
  \qvec_1 \nc R,
\en
where $(\pvec_1,\qvec_1)$ is the first edge of the path
$\ol{\pvec_1 \xvec}$, and $|\pvec_1| > |\pvec_0|$;
\item[(ii)] $\uvec_0$, $\vvec_0$, $\zvec_1$, such that
\eq
\label{e:B1}
  \zerovec \conn \uvec_0,\, \uvec_0 \conn \zvec_1,\,
  \vvec_0 \conn R;
\en
\item[(iii)] $\uvec_0 \conn \pvec_1$;
\item[(iv)] $\vvec_0$ lies either on
    $\ol{\uvec_0 \pvec_1}$, in which case $\pvec_0 = \vvec_0$,
    or on $\ol{\uvec_0 \zvec_1}$, in which case
    $\pvec_0 = \uvec_0$;
\item[(v)] $\pvec_0 \conn \pvec_1 \conn \zvec_1$.
\end{itemize}
Letting
\eqalign
  A_1
  &= \{ \ol{\pvec_0 \pvec_1}, \ol{\pvec_1 \zvec_1} \}, \nonumber \\
  B_1
  &= A_0 \cup \{\ol{\rvec_0 \zvec_1}\}
  = \{ \text{paths in \eqref{e:B1}} \}
     \cup \{ \ol{\uvec_0 \pvec_0} \}
  , \nonumber
\enalign
the hypotheses of Lemma~\ref{lem:intersect} are satisfied with
$\pvec = \pvec_1$, $\qvec = \qvec_1$, $A = A_1$ and $B = B_1$.}

\medbreak

The next step of the construction is carried out similarly.
An application of Lemma~\ref{lem:intersect} gives the paths
shown in Figure~\reffg{recursion}(b). Again, we chose
$\pvec'$ so that $|\pvec'|$ is maximal, and set
$\pvec_2 = \pvec'$, $\zvec_2 = \zvec$ and $\rvec_1 = \rvec$
for this choice of $\pvec'$. We define $\uvec_1$ and
$\vvec_1$ depending on the location of $\rvec_1$,
similarly to the previous step.

If $\pvec_2 = \xvec$, we can conclude similarly to the previous
step that the lemma holds with $J = 2$. If $\pvec_2 \not= \xvec$,
as in Figure~\reffg{recursion}(b), we advance the induction
similarly to the previous step. This time, we use both the
choice of $\pvec_1$ and $\pvec_2$ to conclude the necessary
statement about $A_2$ and $B_2$.

\medbreak

Now we state the recursion hypothesis in general for $I \ge 1$.

\medbreak

\noindent
\emph{{\bf (R) Recursion hypothesis for $I \ge 1$.}
Vertices and paths (all edge-disjoint) with the following
properties exist:
\begin{itemize}
\item[(i)] $\pvec_I$ and $\qvec_I$ such that
\eq
\label{e:restI}
  \pvec_I \conn \wvec \conn \xvec, \qquad
  \qvec_I \nc R,
\en
where $(\pvec_I,\qvec_I)$ is the first edge of the path
$\ol{\pvec_I \xvec}$, and $|\pvec_I| > |\pvec_{I-1}|$;
\item[(ii)] $\uvec_i$, $0 \le i < I$; $\vvec_i$, $0 \le i < I$;
$\zvec_i$, $1 \le i \le I$, such that
\eqalign
\label{e:BI-1}
  &\text{Lemma~\ref{lem:paths} (iii) holds with $i$ restricted to
    $1 \le i < I$}, \\
\label{e:BI-2}
  &\text{Lemma~\ref{lem:paths} (iv) holds with $i$ restricted to
    $1 \le i \le I$}, \\
\label{e:BI-3}
  &\text{Lemma~\ref{lem:paths} (v) holds with $i$ restricted to
    $1 \le i < I$}, \\
\label{e:BI-4}
  &\vvec_0 \conn R;
\enalign
\item[(iii)] $\uvec_{I-1} \conn \pvec_I$;
\item[(iv)] $\vvec_{I-1}$ lies either on
    $\ol{\uvec_{I-1} \pvec_I}$, in which case $\pvec_{I-1} = \vvec_{I-1}$,
    or on $\ol{\uvec_{I-1} \zvec_I}$, in which case
    $\pvec_{I-1} = \uvec_{I-1}$;
\item[(v)] $\pvec_{I-1} \conn \pvec_I \conn \zvec_I$.
\end{itemize}
Letting
\eqalign
  A_I
  &= \{ \ol{\pvec_{I-1} \pvec_I}, \ol{\pvec_I \zvec_I} \}, \nonumber \\
  B_I
  &
  = B_{I-1} \cup A_{I-1} \cup \{ \ol{\rvec_{I-1} \zvec_{I}} \}
  = \{ \text{paths in \eqref{e:BI-1}--\eqref{e:BI-4}} \}
     \cup \{ \ol{\uvec_{I-1} \pvec_{I-1}} \}, \nonumber
\enalign the hypotheses of Lemma~\ref{lem:intersect} are satisfied
with $\pvec = \pvec_I$, $\qvec = \qvec_I$, $A = A_I$ and $B =
B_I$.}

\medbreak

Figure~\reffg{recursion} illustrates
those paths of Figure \reffg{paths} that have been constructed
at the stages $I = 1$ and $I = 2$. Note that $\pvec_I$
receives either the label $\uvec_I$ or $\vvec_I$. 
Hence $\pvec_i$ will always equal
either $\uvec_i$ or $\vvec_i$, depending on the location of
$\vvec_i$ (by part~(iv) of the hypothesis). 
Note also that \refeq{pmon} holds if {\bf (R)} holds for all $I=0,1,\ldots, i$.

\medbreak

\noindent
{\bf Consequence of (R): definition of $\pvec_{I+1}$,
$\uvec_I$, $\vvec_I$ and $\zvec_{I+1}$.}
We now assume that {\bf (R)} holds for some $I \ge 0$.
An application of Lemma \ref{lem:intersect} with the data given
in the hypothesis shows the existence of vertices $\pvec'$,
$\rvec$ and $\zvec$ with certain properties. We now choose
$\pvec'$ so that $|\pvec'|$ be maximal, and such that
the properties claimed in Lemma \ref{lem:intersect} hold.
We set $\pvec_{I+1} = \pvec'$, $\zvec_{I+1} = \zvec$ and
$\rvec_I = \rvec$ for this choice.

Note that $\rvec_I \in A_I$, which is a union of two paths
in both cases $I = 0$ and $I \ge 1$.
In the case $I = 0$, if $\rvec_0 \in \ol{\pvec_0 R}$
then we set $\vvec_0 = \rvec_0$ and
$\uvec_0 = \pvec_0$,
and if $\rvec_0 \in \ol{\zerovec \pvec_0}$ then
we set $\vvec_0 = \pvec_0$, $\uvec_0 = \rvec_0$.
Similarly, in the case $I \ge 1$, we set
$\vvec_I = \rvec_I$ and $\uvec_I = \pvec_I$ if
$\rvec_I \in \ol{\pvec_I \zvec_I}$, and we set
$\vvec_I = \pvec_I$, $\uvec_I = \rvec_I$
if $\rvec_I \in \ol{\pvec_{I-1} \pvec_I}$.
In both cases, it is clear that
$|\uvec_I| \le |\pvec_I| < |\zvec_{I+1}| < R$,
and hence \eqref{e:vzR_i} holds for $i = I+1$.

It follows immediately from these definitions, and from the
disjointness properties ensured by Lemma~\ref{lem:intersect},
that assumptions (ii)--(v) of {\bf (R)} now hold with
$I$ replaced by $I + 1$.

\medbreak

\noindent
{\bf Verification of Lemma~\ref{lem:paths} if
$\pvec_{I+1} = \xvec$.}
We show that if $\pvec_{I+1} = \xvec$, then Lemma~\ref{lem:paths}
holds with $J = I+1$. For this, we define
$\uvec_{I+1} = \wvec$ and $\vvec_{I+1} = \xvec$.
It is immediate from these definitions, from the disjointness properties
ensured by Lemma \ref{lem:intersect}, and from the already established
properties (ii)--(v) of hypothesis {\bf (R)} for $I + 1 = J$,
that (i)--(v) of Lemma \ref{lem:paths} hold.

It remains to show that either Case (a) or Case (b) holds.
Since $\zerovec \conn \yvec$, there exists
$\vvec_* \in \ol{\zerovec \uvec_0}$, such that
$\vvec_* \conn \yvec$ disjointly from
$\ol{\zerovec \uvec_0}$. If $\ol{\vvec_* \yvec}$ is not
disjoint from $\ol{\vvec_0 R}$, we are in Case (a), and we
can ignore $\vvec_*$. If $\ol{\vvec_* \yvec}$ intersects
$\ol{\uvec_0 \pvec_0}$ or $\ol{\uvec_0 \zvec_1}$, let
$\vvec_0'$ be the last such intersection. Note that
$\ol{\vvec_* \yvec}$ must be disjoint from all other
paths constructed, since those are subsets of
$\Ccal(\qvec_0)$, and $\qvec_0 \nc R$. Hence if the
intersection $\vvec_0'$ exists, we can replace $\vvec_0$
by $\vvec_0'$ and we are in Case (a). If the intersection
$\vvec_0'$ does not exist, we are in Case (b).
This verifies the claims of Lemma \ref{lem:paths}.
\medbreak

We are left to show that if $\pvec_{I+1} \not= \xvec$,
then {\bf (R)} must hold for $I+1$.

\medbreak

\noindent
{\bf Advancing the recursion $I \Longrightarrow I+1$
if $\pvec_{I+1} \not= \xvec$.}
Since $\pvec_{I+1} \in \ol{\qvec_I \xvec}$, but $\pvec_{I+1} \not= \xvec$,
we have $|\pvec_{I+1}| > |\pvec_I|$, and $\pvec_{I+1} \conn \wvec$,
showing (i) of hypothesis {\bf (R)}. We have already seen that
(ii)--(v) are guaranteed to hold.

We are left to show that the hypotheses of Lemma \ref{lem:intersect}
hold with the data given. (i), (ii) and $\zerovec \in A_{I+1} \cup B_{I+1}$
are clear from the definitions. By the definition of $\qvec_{I+1}$,
$A_{I+1} \cup B_{I+1}$ is a subgraph of $\tCcal{\qvec_{I+1}}$.

Assume, for a contradiction, that we have
$\zvec_* \in \Ccal(\qvec_{I+1}) \cap (A_{I+1} \cup B_{I+1})$.
Without loss of generality, assume that $\zvec_*$ is the first
visit of an occupied path $\ol{\qvec_{I+1} \zvec_*}$ to
$A_{I+1} \cup B_{I+1}$. In particular, $\ol{\qvec_{I+1} \zvec_*}$
is edge-disjoint from $A_{I+1} \cup B_{I+1}$. Observe that
\eq
\nnb
  A_{I+1} \cup B_{I+1}
  = A_{I+1} \cup A_I \cup B_I \cup \{ \ol{\rvec_I \zvec_{I+1}} \}.
\en
If we had $\zvec_* \in A_{I+1}$, then the disjoint paths
$\ol{\qvec_{I+1} \zvec_* \zvec_{I+1}}$ and
$\ol{\rvec_I \zvec_{I+1}}$ would satisfy the conclusions of
Lemma \ref{lem:intersect} for $\pvec = \pvec_I$, $\qvec = \qvec_I$,
etc. This contradicts the choice of $\pvec_{I+1}$ (the
maximality of $|\pvec_{I+1}|$), since
$|\qvec_{I+1}| > |\pvec_{I+1}|$.
If we had $\zvec_* \in \ol{\rvec_I \zvec_{I+1}}$, we get
a similar contradiction due to the paths
$\ol{\qvec_{I+1} \zvec_*}$ and $\ol{\rvec_I \zvec_*}$.
Finally, we can rule out $\zvec_* \in A_I \cup B_I$, since
$\Ccal(\qvec_{I+1}) \subset \Ccal(\qvec_I)$, and the latter is
disjoint from $A_I \cup B_I$.

We are left to show that every occupied path from $B_{I+1}$ to
$\Ccal(\qvec_{I+1})$ has to pass through $A_{I+1}$. Assume, for a
contradiction, that there exists $\zvec_* \in \Ccal(\qvec_{I+1})$,
and $\zvec_*' \in B_{I+1}$ such that $\zvec_*' \conn \zvec_*$
disjointly from $A_{I+1}$. By considering the last visit, we may
also assume that $\zvec_*'$ is the only vertex of
$\ol{\zvec_*' \zvec_*}$ in $A_{I+1} \cup B_{I+1}$.
We may also assume that $\ol{\qvec_{I+1} \zvec_*}$ and
$\ol{\zvec_*' \zvec_*}$ are edge-disjoint. We already
saw $\Ccal(\qvec_{I+1}) \cap (A_{I+1} \cup B_{I+1}) = \es$,
in particular, $\ol{\qvec_{I+1} \zvec_*}$ is edge-disjoint
from $A_{I+1} \cup B_{I+1}$. Observe that
\eq
\label{e:BI+1}
  B_{I+1}
  = A_I \cup B_I \cup \{ \ol{\rvec_I \zvec_{I+1}} \}
  = \bigcup_{i = 0}^I (A_i \cup \{ \ol{\rvec_i \zvec_{i+1}} \}).
\en
If we had $\zvec_*' \in \ol{\rvec_i \zvec_{i+1}}$, then
the paths $\ol{\qvec_{I+1} \zvec_*}$ and $\ol{\rvec_i \zvec_*' \zvec_*}$
would contradict the choice of $\pvec_{i+1}$.
Finally, if we had $\zvec_*' \in A_i$, then the paths
$\ol{\qvec_{I+1} \zvec_*}$ and $\ol{\zvec_*' \zvec_*}$ would
contradict the choice of $\pvec_{i+1}$.
This completes the verification of hypothesis {\bf (R)} for
$I + 1$.

\medbreak

This completes the proof of Lemma \ref{lem:paths}.
\qed

\subsection{A diagrammatic bound}
\label{sec-diag}

In this section, we use Lemma~\ref{lem:paths} and the BK inequality \cite{BK85}
to bound $\bP_{p_c} [A_J(n,\wvec,\xvec,\yvec)]$.
For this, we need the following preliminaries.

The \emph{critical survival probability} is defined by
\eq
\label{e:thetan}
  \theta_N
  = \bP_{p_c} ( \zerovec \conn N ).
\en
The {two papers \cite{HHS05a,HHS05b} show that
for $d > 4$ and $L \ge L_0(d)$, we have
$\theta_N \sim c N^{-1}$ as $N \to \infty$,
for some $c = c(d,L)=2+\Ocal(L^{-d})$.}
Moreover,
\eq
\label{e:survive}
  \theta_N
  \le \frac{K'}{N}, \quad \quad N \ge 0, \;\; L \ge L_0,
\en
with the constant $K'=5$ which is of course independent
of both $d$ and $L$ (see \cite[Eqn.~(1.11)]{HHS05a}).

To abbreviate the notation, when $\yvec_1 = (y_1,m_1)$ and $\yvec_2 = (y_2,m_2)$
we write
$\tau(\yvec_1,\yvec_2) = \tau_{m_2 - m_1} (y_2 - y_1)$.
We also introduce
\eqsplit
  U_1 (\uvec_0,\vvec_0,\uvec_1,\vvec_1,\zvec_1)
  &= \tau(\vvec_0,\uvec_1)\,
     \tau(\uvec_1,\vvec_1)\,
     \tau(\vvec_1,\zvec_1)\,
     \tau(\uvec_0,\zvec_1) \\
  U_2 (\uvec_0,\vvec_0,\uvec_1,\vvec_1,\zvec_1)
  &= \tau(\uvec_0,\uvec_1)\,
     \tau(\uvec_1,\vvec_1)\,
     \tau(\vvec_1,\zvec_1)\,
     \tau(\vvec_0,\zvec_1) \\
  U &= U_1 + U_2.
\ensplit
For $0 \le |\uvec_0| < n$ and $|\uvec_0| \le |\vvec_0| < R$
and $\yvec = (y,N)$, let
\eqsplit
  \vphi (\uvec_0,\vvec_0)
  &= \sum_{y \in \Zd}
     \tau(\zerovec,\uvec_0)\,
     \tau(\uvec_0,\vvec_0)\,
     \tau(\vvec_0,\yvec), \\
  \vphi_R (\uvec_0,\vvec_0)
  &= \sum_{y \in \Zd}
     \sum_{
     \vvec_* \in \Zd \times \Z_+
     }
     \tau(\zerovec,\vvec_*)\,
     \tau(\vvec_*,\uvec_0)\,
     \tau(\uvec_0,\vvec_0)\,
     \theta_{R - |\vvec_0|}\,
     \tau(\vvec_*,\yvec) \\
  \psi^{(0)} (\uvec_0,\vvec_0)
  &= \vphi (\uvec_0,\vvec_0) + \vphi_R (\uvec_0,\vvec_0).
\ensplit
For $I \ge 1$, $0 \le |\uvec_I| < n$ and $|\uvec_I| \le |\vvec_0| < R$, let
\eqsplit
\label{e:psiJdef}
  \psi^{(I)} (\uvec_I,\vvec_I)
  &= \sum_{\substack{\uvec_{I-1} \in \Zd \times Z_+ \\
     0 \le |\uvec_{I-1}| \le |\uvec_I|}} \ \
     \sum_{\substack{\zvec_{I} \in \Zd \times Z_+ \\
     |\vvec_I| < |\zvec_I| < R}} \ \
     \sum_{\substack{\vvec_{I-1} \in \Zd \times Z_+ \\
     |\uvec_{I-1}| \le |\vvec_{I-1}| \le |\zvec_I|}} \ \
     U(\uvec_{I-1},\vvec_{I-1},\uvec_I,\vvec_I,\zvec_I) \\
  &\qquad\qquad \times \psi^{(I-1)} (\uvec_{I-1},\vvec_{I-1}).
\ensplit

\begin{lemma}
\label{lem:AJs}
For $J \geq 0$,
\eq
\label{e:AJbound}
  \sum_{y \in \Zd} \bP_{p_c} \left[ A_J(n,\wvec,\xvec,\yvec) \right]
  \le \psi^{(J)}(\wvec,\xvec).
\en
\end{lemma}

\begin{proof}
Definition~\ref{def-AJ} guarantees that on the event
$A_J(n,\wvec,\xvec,\yvec)$ certain disjoint paths
exist. If we fix the vertices
$\uvec_0,\dots, \uvec_J$,
$\vvec_0,\dots, \vvec_J$ and
$\zvec_1,\dots,\zvec_J$, then the probability of the
existence of the disjoint paths is bounded by
the product of the probabilities of the existence of the
individual paths, by the BK inequality \cite{BK85}.
An individual path $\ol{\yvec_1 \yvec_2}$ contributes
a factor $\tau(\yvec_1,\yvec_2)$.
Now summing the bound over all the vertices but $\uvec_J$
and $\vvec_J$, gives an upper bound on
$\bP_{p_c} \left[ A_J(n,\wvec,\xvec,\yvec) \right]$.
Further summing over $y \in \Zd$ gives an upper bound for
the left-hand side of \eqref{e:AJbound}.

Now it is merely a matter of bookkeeping to check that we get
the expressions $\psi^{(J)}$. The terms $\vphi$ and
$\vphi_R$ correspond to Cases (a) and (b) of
Lemma~\ref{lem:paths}, respectively, and their sum
$\psi^{(0)}$ bounds the contribution of the paths
constructed when we initialized the recursion,
together with the path leading to $\yvec$.
When $J = 0$, and we take $\uvec_0 = \wvec$ and
$\vvec_0 = \xvec$, we get the bound
in \eqref{e:AJbound}, with $J = 0$.

When $J \ge 1$, the recursive definition of $\psi^{(I)}$
reflects the recursion of Lemma~\ref{lem:paths}.
The factor $U = U_1 + U_2$ gives the contribution of
the paths added in the $I$-th step: for $U_1$ these
are $\ol{\vvec_{I-1} \uvec_I}$, $\ol{\uvec_I \vvec_I}$,
$\ol{\vvec_I \zvec_I}$ and $\ol{\uvec_{I-1} \zvec_I}$
(when $\vvec_{I-1}$ lies on $\ol{\uvec_{I-1} \uvec_I}$),
and for $U_2$ they are
$\ol{\uvec_{I-1} \uvec_I}$, $\ol{\uvec_I \vvec_I}$,
$\ol{\vvec_I \zvec_I}$ and $\ol{\vvec_{I-1} \zvec_I}$
(when $\vvec_{I-1}$ lies on $\ol{\uvec_{I-1} \zvec_I}$).
Note that the path $\ol{\uvec_{I-1} \vvec_{I-1}}$ is
not present in $U$, since it is taken care of
inside $\psi^{(I-1)}$.
\end{proof}

\subsection{Estimation of diagrams}
\label{ssec:estimates}

It follows from \refeq{EDlim}, \refeq{Dnincl} and Lemma~\ref{lem:AJs} that
\eqalign
\label{e:Dnbound}
  \bE_\infty |D(n)|
  &\le \frac{1}{A}\limsup_{N \to \infty}  \sum_{w,x,y \in \Zd} \sum_{J = 0}^\infty
      \bP_{p_c} \left[ A_J(n,\wvec,\xvec,\yvec) \right]
  \nnb &
  \le \frac{1}{A}\limsup_{N \to \infty}
      \left[ \sum_{J = 0}^\infty \sum_{w,x \in \Zd} \psi^{(J)}(\wvec,\xvec) \right].
\enalign
{Fix $d>6$, $R \geq 1$, $0<a<1$ and $0 < n \leq \lfloor aR \rfloor$.
To prove \refeq{EDgoal} and hence Proposition~\ref{prop:piv},
it suffices to show
that there exist $c_2 = c_2(a)$ and a constant $0 < c_3 < \frac 12$} such that
\eq
\label{e:AJbd}
  \limsup_{N \to \infty}  \sum_{w,x \in \Zd}
    \psi^{(J)} (\wvec,\xvec)
  \le c_2 c_3^J, \quad J \ge 0,
\en
since \refeq{taun} and \eqref{e:Dnbound}--\eqref{e:AJbd} then imply that
\eq
\nnb
  \bE_\infty |D(n)|
  \le  \bar{K} c_2 \sum_{J=0}^\infty c_3^J
  = \bar{K} \frac{c_2}{1 - c_3} \leq 2\bar K c_2
  = {c_1(a)}.
\en

{We now state and prove two lemmas which imply \refeq{AJbd}.
Their proofs use the bound
\eq
\label{e:taunK}
  \tau_n
  \le \bar{K}, \quad n \ge 0,
\en
of \eqref{e:taun}, as well as \eqref{e:survive}.
It is in Lemma~\ref{lem:psiJs}, and only
there, that we need to assume $d > 6$ rather than $d>4$.}
The first lemma gives a bound on $\psi^{(0)}$.

\begin{lemma}
\label{lem:psi0}
{Let $d>4$, $R \geq 1$, $0<a<1$, $0 < n \leq \lfloor aR \rfloor$,}
$\wvec = (w,n-1)$ and $\xvec = (x,n)$.  Then
\eq
\label{e:psi0bnd}
  \limsup_{N \to \infty}
  \sum_{w,x \in \Zd} \psi^{(0)} (\wvec,\xvec)
  \le ( \bar{K}^3 + \bar{K}^4 K' a/(1-a) ).
\en
\end{lemma}

\begin{proof}
By definition and \refeq{taunK},
\eq
\nnb
  \sum_{w,x \in \Zd} \vphi (\wvec,\xvec)
  = \tau_{n-1} \tau_1 \tau_{N-n}
  \le \bar{K}^3.
\en
Similarly, writing $\vvec_* = (v_*,l_*)$,
\eq
\nnb
  \sum_{w,x \in \Zd} \vphi_R (\wvec,\xvec)
  = \sum_{l_* = 0}^{n-1}
    \tau_{l_*}\, \tau_{n-l_*-1}\, \tau_1\, \theta_{R-n}\, \tau_{N-l_*}
  \le \frac{\bar{K}^4 K' n}{R-n}.
\en
Since $n/(R-n) \le a/(1-a)$ because {$n\le\lfloor aR \rfloor$}, this gives \eqref{e:psi0bnd}.
\end{proof}

For $J \ge 1$, we use a somewhat stronger
formulation of the bound, in which $|\uvec_J|$ and $|\vvec_J|$ are
not restricted to the values $n-1$ and $n$. This will allow us
to prove a bound on $\psi^{(J)}$ by induction.

\begin{lemma}
\label{lem:psiJs}
{Let $d>6$, $R \geq 1$, $0<a<1$, $0 < n \leq \lfloor aR \rfloor$.}
Suppose that $0 \le k_J < n$, $k_J \le l_J < R$,
$\uvec_J = (u_J,k_J)$ and $\vvec_J = (v_J,l_J)$. Then
\eq
\label{e:psiJbnd}
  \limsup_{N \to \infty}
      \sum_{u_J,v_J \in \Zd} \psi^{(J)} (\uvec_J,\vvec_J)
  \le (2 \bar{K}^3 K^3 \beta)^J ( \bar{K}^3 + 3 \bar{K}^5 K' a/(1-a) ),
    \quad J \ge 1.
\en
\end{lemma}

\begin{proof}
We start by inserting
the definition of $\psi^{(J)}$ into the left-hand side of
\eqref{e:psiJbnd}.  With $\zvec_J = (z_J,s_J)$,
$\uvec_{J-1} = (u_{J-1},k_{J-1})$ and $\vvec_{J-1} = (v_{J-1},l_{J-1})$,
the left-hand side of \eqref{e:psiJbnd} equals
\eqsplit
\label{e:firststep0}
  &\limsup_{N \to \infty}  \sum_{u_J,v_J \in \Zd} \
  \sum_{z_J,u_{J-1},v_{J-1} \in \Zd} \
     \sum_{k_{J-1} = 0}^{k_J} \
     \sum_{s_J = l_J}^{R-1} \ \sum_{l_{J-1} = k_{J-1}}^{s_J} \
     U(\uvec_{J-1},\vvec_{J-1},\uvec_J,\vvec_J,\zvec_J) \\
  &\qquad\qquad \times \psi^{(J-1)} (\uvec_{J-1},\vvec_{J-1}).
\ensplit

The vertices $\uvec_J$, $\vvec_J$ and $\zvec_J$ only appear in the
factor $U$. We claim that
\eq
\label{e:bubblebound}
  \sum_{u_J,v_J,z_J \in \Zd}
     U(\uvec_{J-1},\vvec_{J-1},\uvec_J,\vvec_J,\zvec_J)
  \le 2 \bar{K}^3 K \beta (s_J - k_{J-1} + 1)^{-d/2}.
\en
To see this, note that  $s_J = |\zvec_J| > |\uvec_{J-1}| = k_{J-1}$,
by \eqref{e:vzR_i}.
For the $U_1$ term, we use \eqref{e:taubound} to bound
$\tau(\uvec_{J-1},\zvec_J)$
by $K \beta (s_J - k_{J-1} + 1)^{-d/2}$. Then the sums over
$z_J$, $v_J$ and $u_J$ contribute the factor $\bar{K}^3$,
by using \eqref{e:taunK} for the other three factors in $U_1$.
For the $U_2$ term,
we apply \eqref{e:taubound} and $\tau_n \le \bar{K}$ to see that
\eq
\label{e:space-conv}
  \sup_{x \in \Zd} \sum_{y \in \Zd} \tau_n (y) \tau_m (x - y)
  \le K \beta (n + m + 1)^{-d/2}, \quad n + m \ge 1.
\en
An application of \eqref{e:space-conv}
to the convolution of $\tau(\uvec_{J-1},\uvec_J)$,
$\tau(\uvec_J,\vvec_J)$ and $\tau(\vvec_J,\zvec_J)$, together with
\eqref{e:taunK}, yields an upper bound of the same form.
This proves \eqref{e:bubblebound}.
Inserting \eqref{e:bubblebound} into \eqref{e:firststep0}
and rearranging, we get
\eqalign
  \refeq{firststep0} &\leq 2 \bar{K}^3 K \beta \sum_{k_{J-1} = 0}^{k_J} \
     \sum_{s_J = l_J}^{R-1} \ (s_J - k_{J-1} + 1)^{-d/2} \nonumber\\
  &\qquad\qquad \times \sum_{l_{J-1} = k_{J-1}}^{s_J}
      \limsup_{N \to \infty}
     \sum_{u_{J-1},v_{J-1} \in \Zd} \
     \psi^{(J-1)} (\uvec_{J-1},\vvec_{J-1})  .
\label{e:one-step}
\enalign

Now we prove \eqref{e:psiJbnd} by induction on $J$.
To start the induction, we verify \eqref{e:psiJbnd} for $J = 1$.
This is most of the work; advancing the induction is easy.
When $J = 1$, the $\limsup$ in \eqref{e:one-step}
consists of two terms, corresponding
to $\vphi$ and $\vphi_R$. The $\vphi$-term is bounded by
\eqsplit
\label{e:phicase}
  &\limsup_{N \to \infty}
    \sum_{u_0,v_0 \in \Zd} \
    \sum_{y \in \Zd} \ \tau(\zerovec,\uvec_0)\, \tau(\uvec_0,\vvec_0)\,
    \tau(\vvec_0,\yvec)
  = \limsup_{N \to \infty}\ \tau_{k_0} \tau_{l_0 - k_0} \tau_{N - l_0}
  \le \bar{K}^3.
\ensplit
Inserting this into \eqref{e:one-step}, and assuming $d>6$,
we see that the $\varphi$ contribution to \eqref{e:one-step}
is bounded by
\eq
\lbeq{wxyz}
  2 \bar{K}^3 K \beta \bar{K}^3 \sum_{k_0 = 0}^{k_1} \
     \sum_{s_1 = l_1}^{R-1} \ (s_1 - k_0 + 1)^{(2-d)/2}
  \le (2 \bar{K}^3 K^2 \beta) (\bar{K}^3).
\en

The $\vphi_R$ term is bounded as follows. First, the
$\limsup$ is bounded by
\eqalign
  &\limsup_{N \to \infty}
    \sum_{u_0,v_0 \in \Zd} \
    \sum_{l_* = 0}^{k_0} \
    \sum_{y,v_* \in \Zd} \ \tau(\zerovec,\vvec_*)\, \tau(\vvec_*,\uvec_0)\,
    \tau(\uvec_0,\vvec_0)\, \theta_{R-l_0}\,
    \tau(\vvec_*,\yvec) \nonumber \\
  &\qquad
   \le \bar{K}^4 \sum_{l_* = 0}^{k_0} \theta_{R-l_0}
  \le \bar{K}^4 (k_0 + 1) \frac{K'}{R-l_0}
  \le \bar{K}^4 K' n \frac1{R - l_0}.
\enalign
We insert this bound into   \eqref{e:one-step} to obtain
\eq
\label{e:Rcase}
  (2 \bar{K}^3 K \beta) (\bar{K}^4 K') n \sum_{k_0 = 0}^{k_1} \
     \sum_{s_1 = l_1}^{R-1} \ (s_1 - k_0 + 1)^{-d/2}
     \sum_{l_0 = k_0}^{s_1} \frac1{R - l_0}.
\en
We split the sum over $s_1$ into the cases: (1) $s_1 < n + (R-n)/2$;
(2) $s_1 \ge n + (R-n)/2$. In case~(1), we have
\eq
\nnb
  \frac1{R - l_0}
  \le \frac1{R - s_1}
  \le \frac{2}{R - n}.
\en
Inserting this into \eqref{e:Rcase}, the contribution of case (1)
to the expression in \eqref{e:Rcase} is bounded by
\eqsplit
\label{e:case1}
   &(2 \bar{K}^3 K \beta ) (2 \bar{K}^4 K') \frac{n}{R-n} \sum_{k_0 = 0}^{k_1} \
     \sum_{s_1 = l_1}^{n + (R-n)/2} \ (s_1 - k_0 + 1)^{(2-d)/2} \\
   &\qquad \le (2 \bar{K}^3 K^2 \beta) (2 \bar{K}^4 K') \frac{n}{R-n}
   \le (2 \bar{K}^3 K^2 \beta) (2 \bar{K}^4 K') \frac{a}{1-a}.
\ensplit
In case~(2), since $n \geq k_1\geq k_0$ we have
\eq
\nnb
  (s_1 - k_0 + 1)^{-d/2}
  \le K (R - k_0 + 1)^{-d/2},
\en
and the sum over $l_0$ in \eqref{e:Rcase} is bounded by
$\log(R - k_0 + 1) \le \bar{K} (R - k_0 + 1)^\delta$ for
some fixed exponent $\delta$
(e.g., $\delta = 1/4$ suffices).
Therefore the contribution of case (2)
to the expression in \eqref{e:Rcase} is bounded by
\eqsplit
\label{e:case2}
  &(2 \bar{K}^3 K^2  \beta) (\bar{K}^5 K') n \frac{R - n}{2}
     \sum_{k_0 = 0}^{k_1} \
     (R - k_0 + 1)^{(2 \delta - d)/2} \\
  &\qquad \le (2 \bar{K}^3 K^3  \beta) (\bar{K}^5 K') n \frac{R - n}{2}
     (R - n)^{(2 \delta + 2 - d)/2} \\
  &\qquad \le (2 \bar{K}^3 K^3  \beta) (\bar{K}^5 K') \frac{n}{R-n}
     (R - n)^{(2 \delta + 6 - d)/2} \\
  &\qquad \le (2 \bar{K}^3 K^3  \beta) (\bar{K}^5 K') \frac{a}{1-a}.
\ensplit

Putting \eqref{e:case1} and \eqref{e:case2} together, we get
that \eqref{e:Rcase} is bounded by
$(2 \bar{K}^3 K^3 \beta) (3 \bar{K}^5 K' a/(1-a))$.
Together with \eqref{e:wxyz} this proves the $J = 1$ case
of \eqref{e:psiJbnd}.

To advance the induction,
we assume now that \eqref{e:psiJbnd} holds for an integer $J = M - 1 \ge 1$,
and prove that it holds for $J = M$.
Using $d > 6$, we insert the bound \eqref{e:psiJbnd}
into \eqref{e:one-step} to get that the right-hand side of
\eqref{e:one-step} is bounded by
\eqalign
  & (2 \bar{K}^3 K \beta) (2 \bar{K}^3 K^3 \beta)^{M-1}
     (\bar{K}^3 + 3 \bar{K}^5 K' a/(1-a))
     \sum_{k_{M-1} = 0}^{k_M} \
     \sum_{s_M = l_M}^{R-1} \ (s_M - k_{M-1} + 1)^{(2-d)/2} \nonumber \\
  &\qquad \le (2 \bar{K}^3 K^3 \beta)^M
     (\bar{K}^3 + 3 \bar{K}^5 K' a/(1-a)).
\enalign
This completes the proof of \eqref{e:psiJbnd}.
\end{proof}

\medbreak

\noindent
\emph{Proof of \refeq{AJbd}.}
It follows immediately from Lemmas \ref{lem:psi0}--\ref{lem:psiJs}
that \eqref{e:AJbd} holds with $c_2 = (\bar{K}^3 + 3 \bar{K}^5 K' a/(1-a))$
and $c_3 = 2 \bar{K}^3 K^3 \beta$.
Recall that the constant {$K'=5$} of \refeq{survive} is independent of $d$ and $L$.
Choosing $\beta$ small
ensures that {$0 < c_3 < \frac 12$}. This proves \refeq{AJbd},
and thus completes the proof of Proposition~\ref{prop:piv}.
\qed

\section*{Acknowledgements}
The work of MTB, AAJ and GS was supported in part by
NSERC of Canada.
The work of TK was supported in part by the Ministry of Education, Culture,
Sports, Science and Technology of Japan, Grant-in-Aid 18654018 (Houga).
{We thank an anonymous referee for suggesting several improvements to
the exposition.}

\end{document}